\documentclass[12pt]{report}

\usepackage{thesis}

\usepackage{amsmath}
\usepackage{url}
\usepackage{subfigure}
\usepackage[pdftex]{graphicx}
\usepackage{multirow}
\usepackage{chapterbib}
\usepackage{amssymb}
\usepackage[titletoc,toc,title]{appendix}
\usepackage{titlesec}
\usepackage{listings,color}
\usepackage{cite}
\usepackage{tabularx}
\usepackage[makeroom]{cancel}

\usepackage{titletoc}

\allowdisplaybreaks
\usepackage{currfile}
\title{Predicting Flow Reversals in a Computational Fluid Dynamics Simulated Thermosyphon using Data Assimilation}
\author{Andrew J. Reagan}
\defensedate{November 15, 2013}
\thesis 			    \masterscience			\apma				      \jangrad				  \def\gradyear{2014}

\advisor{Christopher M. Danforth, Ph.D.}
\readerone{Peter Sheridan Dodds, Ph.D.}
\chair{Darren Hitt, Ph.D.}

\newcommand{\pdiff}[2]{\frac{\partial #1}{\partial #2}}

\newcommand{\diff}[2]{\frac{{\rm d}#1}{{\rm d}#2}}

\newcommand{\rhoref}{\rho_{\text{ref}}}
\newcommand{\dphi}{\text{d}\phi}

\begin{document}

\titlecontents{section}
	      [3.8em]
	      {}
	      {\contentslabel{2.3em}}
	      {\hspace*{-2.3em}}
	      {\titlerule*[1pc]{.}\contentspage}

\maketitle
\makeacceptance

\begin{abstract}
A thermal convection loop is a circular chamber filled with water, heated on the bottom half and cooled on the top half.
With sufficiently large forcing of heat, the direction of fluid flow in the loop oscillates chaotically, forming an analog to the Earth's weather.
As is the case for state-of-the-art weather models, we only observe the statistics over a small region of state space, making prediction difficult.
To overcome this challenge, data assimilation methods, and specifically ensemble methods, use the computational model itself to estimate the uncertainty of the model to optimally combine these observations into an initial condition for predicting the future state.
First, we build and verify four distinct DA methods.
Then, a computational fluid dynamics simulation of the loop and a reduced order model are both used by these DA methods to predict flow reversals.
The results contribute to a testbed for algorithm development.

\end{abstract}

\begin{dedication}
{\it in dedication to }
\\
\vskip 2em
my parents, Donna and Kevin Reagan, for their unwavering support. \end{dedication}

\begin{acknowledgements}
I would like to thank Professors Chris Danforth and Peter Dodds for their outstanding advising, as well as Professors Darren Hitt and Yves Dubief for their valuable feedback.
This work was made possible by funding from the Mathematics and Climate Research Network and Vermont NASA EPSCoR program.
I would like to thank my fellow students for all of their help and toleration through the past two years.
And finally I would like to thank my girlfriend Sam Spisiak, whose love makes this thesis worth writing.
\end{acknowledgements}

\tableofcontents
\listoffigures
\listoftables

\mainmatter
\sloppy

\titlecontents{section}[3.5em]{\vspace{2px}}{}{}{\titlerule*[1pc]{.}\contentspage}
\titlecontents{subsection}[5em]{}{}{}{}

\chapter{Introduction}

\begin{quote}
In this chapter we explore the current state of numerical weather prediction, in particular data assimilation, along with an introduction to computational fluid dynamics and reduced order experiments.
\end{quote}

\section{Introduction}

Prediction of the future state of complex systems is integral to the functioning of our society.
Some of these systems include weather \shortcite{weather-violence2013}, health \shortcite{ginsberg2008detecting}, the economy \shortcite{sornette2006predictability}, marketing \shortcite{asur2010predicting} and engineering \shortcite{savely 1972}.
For weather in particular, this prediction is made using supercomputers across the world in the form of numerical weather model integrations taking our current best guess of the weather into the future.
The accuracy of these predictions depend on the accuracy of the models themselves, and the quality of our knowledge of the current state of the atmosphere.

Model accuracy has improved with better meteorological understanding of weather processes and advances in computing technology.
To solve the initial value problem, techniques developed over the past 50 years are now broadly known as {\em data assimilation}.
Formally, data assimilation is the process of using all available information, including short-range model forecasts and physical observations, to estimate the current state of a system as accurately as possible \shortcite{yang2006}.

We employ a toy climate experiment as a testbed for improving numerical weather prediction algorithms, focusing specifically on data assimilation methods.
This approach is akin to the historical development of current methodologies, and provides a tractable system for rigorous analysis.
The experiment is a thermal convection loop, which by design simplifies our problem into the prediction of convection.
The dynamics of thermal convection loops have been explored under both periodic \shortcite{keller1966} and chaotic \shortcite{welander1967,creveling1975stability,gorman1984,gorman1986,ehrhard1990dynamical,yuen1999,jiang2003,burroughs2005reduced,desrayaud2006numerical,yang2006,ridouane2010} regimes.
A full characterization of the computational behaivor of a loop under flux boundary conditions by Louisos et. al. describes four regimes: chaotic convection with reversals, high Ra aperiodic stable convection, steady stable convection, and conduction/quasi-conduction \shortcite{louisos2013}.
For the remainder of this work, we focus on the chaotic flow regime.

Computational simulations of the thermal convection loop are performed with the open-source finite volume C++ library OpenFOAM \shortcite{jasak2007}.
The open-source nature of this software enables its integration with the data assimilation framework that this work provides.

\section{History of NWP}

The importance of Vilhelm Bjerknes in early developments in NWP is described by Thompson's ``Charney and the revival of NWP'' \shortcite{thompson1990}:

\begin{quote}
It was not until 1904 that Vilhelm Bjerknes - in a remarkable manifesto and testament of deterministic faith - stated the central problem of NWP.
This was the first explicit, coherent, recognition that the future state of the atmosphere is, \emph{in principle}, completely determined by its detailed initial state and known boundary conditions, together with Newton's equations of motion, the Boyle-Charles-Dalton equation of state, the equation of mass continuity, and the thermodynamic energy equation.
Bjerknes went further: he outlined an ambitious, but logical program of observation, graphical analysis of meterological data and graphical solution of the governing equations.
He succeeded in persuading the Norwegians to support an expanded network of surface observation stations, founded the famous Bergen School of synoptic and dynamic meteorology, and ushered in the famous polar front theory of cyclone formation.
Beyond providing a clear goal and a sound physical approach to dynamical weather prediction, V. Bjerknes instilled his ideas in the minds of his students and their students in Bergen and in Oslo, three of of whom were later to write important chapters in the development of NWP in the US (Rossby, Eliassen, and Fj\"{o}rtoft).
\end{quote}

It then wasn't until 1922 that Lewis Richardson suggested a practical way to solve these equations.
Richardson used a horizontal grid of about 200km, and 4 vertical layers in a computational mesh of Germany and was able to solve this system by hand \shortcite{richardson1965}.
Using the observations available, he computed the time derivative of pressure in the central box, predicting a change of 146 hPa, whereas in reality there was almost no change.
This discrepancy was due mostly to the fact that his initial conditions (ICs) were not balanced, and therefore included fast-moving gravity waves which masked the true climate signal \shortcite{kalnay2003}.
Regardless, had the integration had continued, it would have ``exploded,'' because together his mesh and timestep did not satisfy the Courant-Friedrichs-Lewy (CFL) condition, defined as follows.
Solving a partial differential equation using finite differences, the CFL requires the Courant number to be less than one for a convergent solution.
The dimensionless Courant number relates the movement of a quantity to the size of the mesh, written for a one-dimensional problem as $C = u \Delta t/\Delta x $ for $u$ the velocity, $\Delta t$ the timestep, and $\Delta x$ the mesh size.

Early on, the problem of fast traveling waves from unbalanced IC was solved by filtering the equations of motion based on the quasi-geostrophic (slow time) balance.
Although this made running models feasible on the computers available during WW2, the incorporation of the full equations proved necessary for more accurate forecasts.

The shortcomings of our understanding of the subgrid-scale phenomena do impose a limit on our ability to predict the weather, but this is not the upper limit.
The upper limit of predictability exists for even a perfect model, as Edward Lorenz showed in with a three variable model found to exhibit sensitive dependence on initial conditions \shortcite{lorenz1963}.

Numerical weather prediction has since constantly pushed the boundary of computational power, and required a better understanding of atmospheric processes to make better predictions.
In this direction, we make use of advanced computational resources in the context of a simplified atmospheric experiment to improve prediction skill.

\section{IC determination: data assimilation}

Areas as disparate as quadcopter stabilization \shortcite{achtelik2009visual} to the tracking of ballistic missle re-entry \shortcite{siouris1997tracking} use data assimilation.
The purpose of data assimilation in weather prediction is defined by Talagrand as ``using all the available information, to determine as accurately as possible the state of the atmospheric (or oceanic) flow.'' \shortcite{talagrand1997assimilation}
One such data assimilation algorithm, the Kalman filter, was originally implemented in the navigation system of Apollo program \shortcite{kalman1961new,savely1972}.

Data assimilation algorithms consist of a 3-part cycle: predict, observe, and assimilate.
Formally, the data assimilation problem is solved by minimizing the initial condition error in the presence of specific constraints.
The prediction step involves making a prediction of the future state of the system, as well as the error of the model, in some capacity.
Observing systems come in many flavors: rawindsomes and satellite irradiance for the atmosphere, temperature and velocity reconstruction from sensors in experiments, and sampling the market in finance.
Assimilation is the combination of these observations and the predictive model in such a way that minimizes the error of the initial condition state, which we denote the analysis.

The difficulties of this process are multifaceted, and are addressed by parameter adjustments of existing algorithms or new approaches altogether.
The first problem is ensuring the analysis respects the core balances of the model.
In the atmosphere this is the balance of Coriolis force and pressure gradient (geostrophic balance), and in other systems this can arise in the conservation of physical (or otherwise) quantities.

In real applications, observations are made with irregular spacing in time and often can be made much more frequently than the length of data assimilation cycle.
However the incorporation of more observations will not always improve forecasts: in the case where these observations can only be incorporated at the assimilation time, the forecasts can degrade in quality \shortcite{kalnay20074}.
Algorithms which are capable of assimilating observations at the time they are made have a distinct advantage in this respect.
In addition, the observations are not always direct measurements of the variables of the model.
This can create non-linearities in the observation operator, the function that transforms observations in the physical space into model space.
In Chapter 2 we will see where a non-linear operator complicates the minimization of the analysis error.

As George Box reminds us, ``all models are wrong'' \shortcite{box1978statistics}.
The incorporation of model bias into the assimilation becomes important with imperfect models, and although this issue is well-studied, it remains a concern \shortcite{allgaier2012empirical}.

The final difficulty we need to address is the computational cost of the algorithms.
The complexity of models creates issues for both the computation speed and storage requirements, an example of each follows.
The direct computation of the model error covariance, integrating a linear tangent model which propagates errors in the model forward in time, amounts to $n^2$ model integrations where $n$ is the number of variables in the model.
For large systems, e.g. $n = 10^{10}$ for global forecast models, this is well beyond reach.
Storing a covariance matrix at the full rank of model variables is not possible even for a model of the experiment we will address, which has 240,000 model variables.
Workarounds for these specific problems include approximating the covariance matrix (e.g. with a static matrix, or updating it using an ensemble of model runs), and reformulating the assimilation problem such that the covariance is only computed in observation space, which is typically orders of magnitude smaller.

\section{NWP and fluid dynamics: computational models}

As early as 1823, the basic equations governing flow of a viscous, heat-conducting fluid were derived from Newton's law of motion \shortcite{navier,stokes}.
These are the Navier-Stokes equations.
The equations of geophysical fluid dynamics describe the flow in Earth's atmosphere, differing from the Navier-Stokes equations by incorporating the Coriolis effect and neglecting the vertical component of advection.
In both cases, the principle equations are the momentum equations, continuity equation and energy equation.

As stated, our interest here is not the accuracy of computational fluid dynamic simulations themselves, but rather their incorporation into the prediction problem.
For this reason, we limit the depth of exploration of available CFD models.

We chose to perform simulations in the open-source C++ library OpenFOAM.
Because the code is open-source, we can modify it as necessary for data assimilation, and of particular importance becomes the discretized grid.

There are two main directions in which we do explore the computational model: meshing and solving.
To solve the equations of motion numerically, we first must discretize them in space by creating the computational grid, which we refer to as a mesh.
Properly designing and characterizing the mesh is the basis of an accurate solution.
Still, without the use of (notoriously unreliable) turbulence models, and without a grid small enough to resolve all of the scales of turbulence, we inevitably fail to capture the effects of sub-grid scale phenomena.
In weather models, where grid sizes for global models are in the 1-5km range horizontally, this is a known problem and is an open area of research.
In addition to the laminar equations, to account for turbulence, we explore include the $k-\epsilon$ model, $k-\omega$ model, and large eddy simulations (in 3-dimensions only).

Secondly, the choice of solving technique is a source of numerical and algorithmic uncertainty.
The schemes which OpenFOAM has available are dissapative and are limited in scope.
That is, small perturbations experience damping (such as the analysis perturbations to an ensemble members in ensemble Kalman filters).
Without consideration this problem, the ensemble does not adequately sample the uncertainty of the model to obtain an approximation of the model error growth.

We investigate these meshing and solving techniques in greater detail in Chapter 3.

\section{Toy climate experiments}

\begin{figure}[t!]
  \centering
  \includegraphics[width=0.49\textwidth]{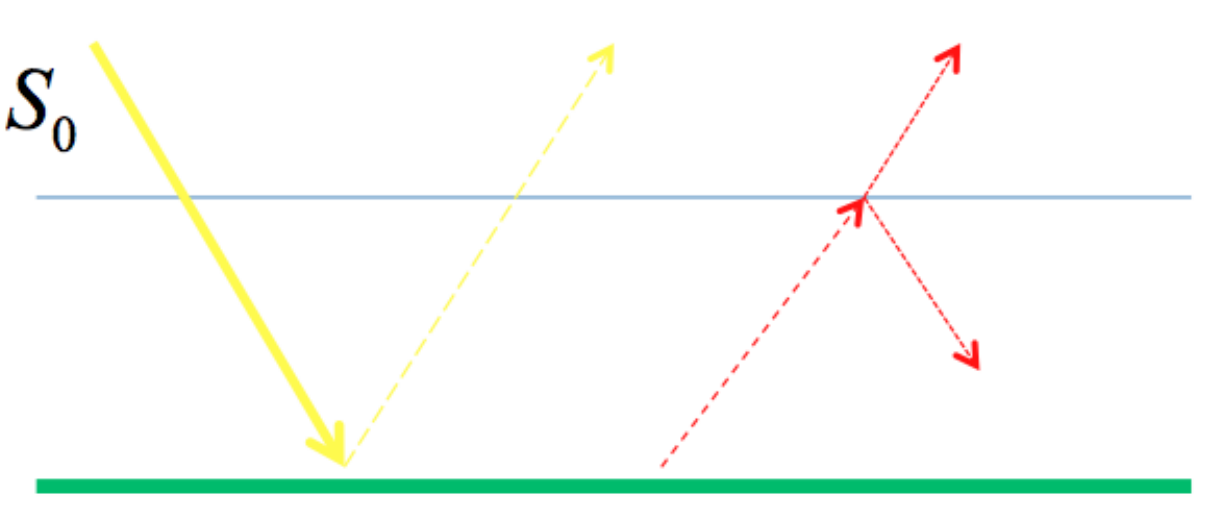}
  \caption[A schematic of the Energy Balance Model (EBM), a conceptual climate model based on the conservation of energy]{
    A schematic of the Energy Balance Model (EBM), a conceptual climate model based on the conservation of energy.
    $S_0$, represented by the yellow arrow, represents energy coming from the sun absorbed by Earth's surface (horizontal green bar).
    The red, dashed arrow represents energy radiating out of Earth's surface, of which some is reflected back by the atmosphere (horizontal gray line).
    Since energy is conserved in the EBM, the parameterizations of surface albedo (amount reflected by Earth's surface), atmospheric absorbivity, and radiation rate determine the long term behaivor of internal energy (temperature).
    This schematic was provided in a lecture given by Dr. Chris Jones {\protect \shortcite{MCRN501}}.
  }
  \label{fig:EBM}
\end{figure}

Mathematical and computational models span from simple to kitchen-sink, depending on the phenomenon of interest.
Simple, or ``toy'', models and experiments are important to gain theoretical understanding of the key processes at play.
An example of such a conceptual model used to study climate is the energy balance model, Figure \ref{fig:EBM}.

With a simple model, it is possible to capture the key properties of simple concepts (e.g. the greenhouse effect in the above energy balance model).
For this reason they have been, and continue to be, employed throughout scientific disciplines.
The ``kitchen-sink'' models, like the Earth System Models (ESMs) used to study climate, incorporate all of the known physical, chemical, and biological processes for the fullest picture of the system.
In the ESM case, the simulations represent our computational laboratory.

We study a model that is more complex than the EBM, but is still analytically tractable, known as the Lorenz 1963 system (Lorenz 63 for short).
In 1962, Barry Saltzmann attempted to model convection in a Rayleigh-B\'{e}rnard cell  by reducing the equations of motion into their core processes \shortcite{saltzman1962finite}.
Then in 1963 Edward Lorenz reduced this system ever further to 3 equations, leading to his landmark discovery of deterministic non-periodic flow \shortcite{lorenz1963}.
These equations have since been the subject of intense study and have changed the way we view prediction and determinism, remaining the simple system of choice for examining nonlinear behaivor today \cite{kalnay20074}.

Thermal convection loops are one of the few known experimental implementations of Lorenz's 1963 system of equations.
Put simply, a thermal convection loop is a hula-hoop shaped tube filled with water, heated on the bottom half and cooled on the top half.
As a non-mechanical heat pump, thermal convection loops have found applications for heat transfer in solar hot water heaters \shortcite{belessiotis2002analytical}, cooling systems for computers \shortcite{beitelmal2002two}, roads and railways that cross permafrost \shortcite{cherry2006next}, nuclear power plants \shortcite{detman1968thermos,beine1992integrated,kwant1992prism}, and other industrial applications.
Our particular experiment consists of a plastic tube, heated on the bottom by a constant temperature water bath, and cooled on top by room temperature air.

\begin{figure}[t!]
  \centering
  \includegraphics[width=0.49\textwidth]{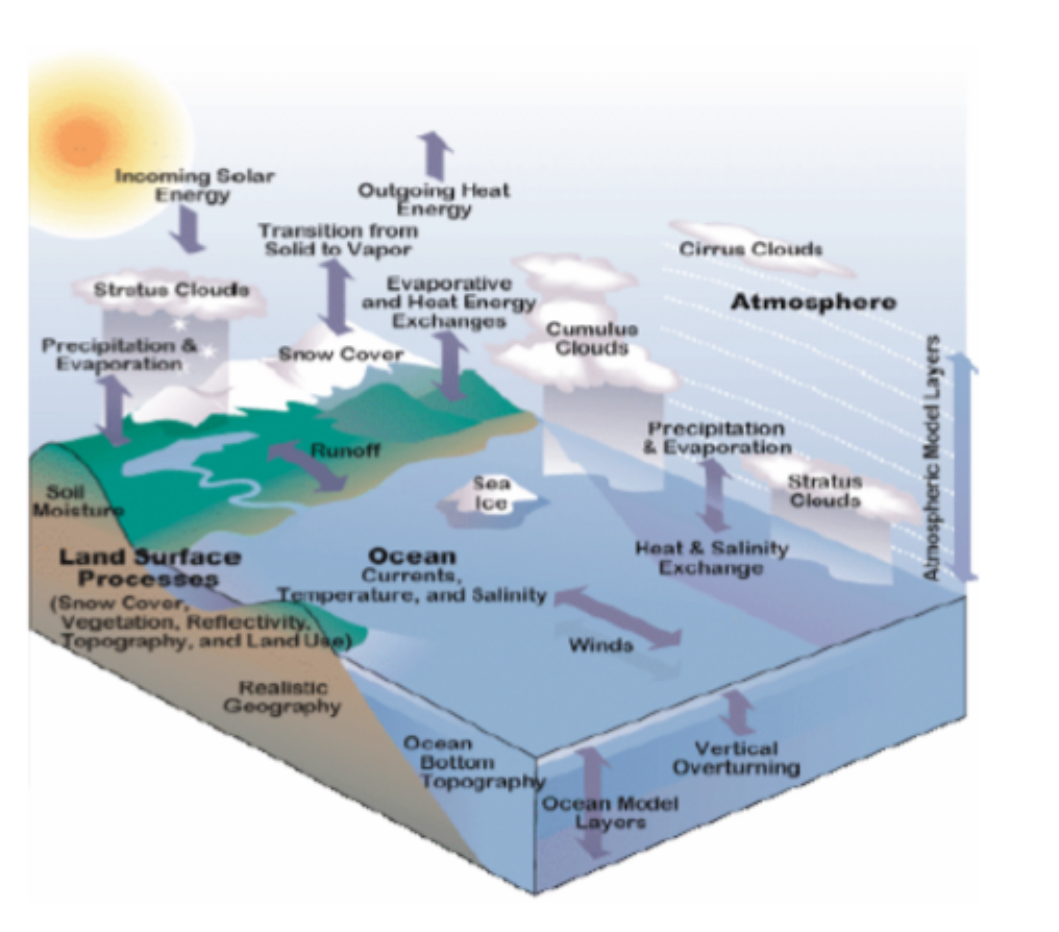}
  \caption[The ``kitchen-sink'' Earth System Model captures as many known processes as computationally feasible]{
    The ``kitchen-sink'' Earth System Model captures as many known processes as computationally feasible.
    Models of this size are most useful in making global climate predictions and global weather predictions.
    By incorporating physically tuned parameterizations of subgrid-scale processes, these models are able to make accurate predictions of Earth's climate and weather.
    Due to the size and complexity of the model, the physical processes driving the model dynamics are, however, difficult to study.
    This schematic describes the Community Earth System Model developed by NCAR {\protect \shortcite{ESMpicture}}.
  }
  \label{fig:ESM}
\end{figure}

Varying the forcing of temperature applied to the bottom half of the loop, there are three experimentally accessible flow regimes of interest to us: conduction, steady convection, and unsteady convection.
The latter regime exhibits the much studied phenomenon of chaos in its changes of rotational direction.
This regime is analytically equivalent the Lorenz 1963 system with the canonical parameter values listed in Table \ref{table:lorenz63params}.
Edward Lorenz discovered and studied chaos most of his life, and is said to have jotted this simple description on a visit to the University of Maryland \cite{danforth2013blog}:
\begin{quote}
Chaos: When the present determines the future, but the approximate present does not approximately determine the future.
\end{quote}

In chapter 3, we present a derivation of Ehrhard-M\"{u}ller Equations governing the flow in the loop following the work of Kameron Harris (2011).
These equations are akin to the Lorenz 63 system, as we shall see, and the experimental implementation is described in more detail in the first section of Chapter 3.

In chapter 4 we present results for the predictive skill tests and summarize the successes of this work and our contributions to the field.
We find that below a critical threshold of observations, predictions are useless, where useless forecasts are defined as those with RMSE greater than 70\% of climatological variance.
The skill of forecasts is also sensitive to the choice of covariance localization, which we discuss in more detail in chapter 4.
And finally, we present ``foamLab,'' a test-bed for DA which links models, DA techniques, and localization in a general framework with the ability to extend these results and test experimental methods for DA within different models.

\newcommand{\mbe}{\mathbf{\epsilon}}
\newcommand{\mbx}{\mathbf{x}}
\newcommand{\mby}{\mathbf{y}}
\newcommand{\mbd}{\mathbf{d}}
\newcommand{\mbB}{\mathbf{B}}
\newcommand{\mbW}{\mathbf{W}}
\newcommand{\mbR}{\mathbf{R}}
\newcommand{\mbH}{\mathbf{H}}
\newcommand{\mbK}{\mathbf{K}}
\newcommand{\mbP}{\mathbf{P}}
\newcommand{\mbZ}{\mathbf{Z}}
\newcommand{\mbw}{\mathbf{w}}
\newcommand{\mbX}{\mathbf{X}}
\newcommand{\mbY}{\mathbf{Y}}
\newcommand{\mb}{\mathbf{}}
\newcommand{\expv}[1]{E \left [ #1 \right ]}

\chapter{Data Assimilation}
\chaptermark{Data Assimilation}
The goal of this chapter is to introduce and detail the practical implementation of state-of-the-art data assimilation algorithms.
We derive the general form for the problem, and then talk about how each filter works.
Starting with the simplest examples, we formulate the most advanced 4-dimensional Variational Filter and 4-dimensional Local Ensemble Transform Kalman Filter to implement with a computational fluid dynamics model.
For a selection of these algorithms, we illustrate some of their properties with a MATLAB implementation of the Lorenz 1963 three-variable model: \begin{align*}
\frac{dx}{dt} &= \sigma (y-x)\\
\frac{dy}{dt} &= \rho x - y -xz \\
\frac{dz}{dt} &= xy -  \beta z\end{align*}

\begin{figure}[h!]
  \centering
  \includegraphics[width=0.49\textwidth]{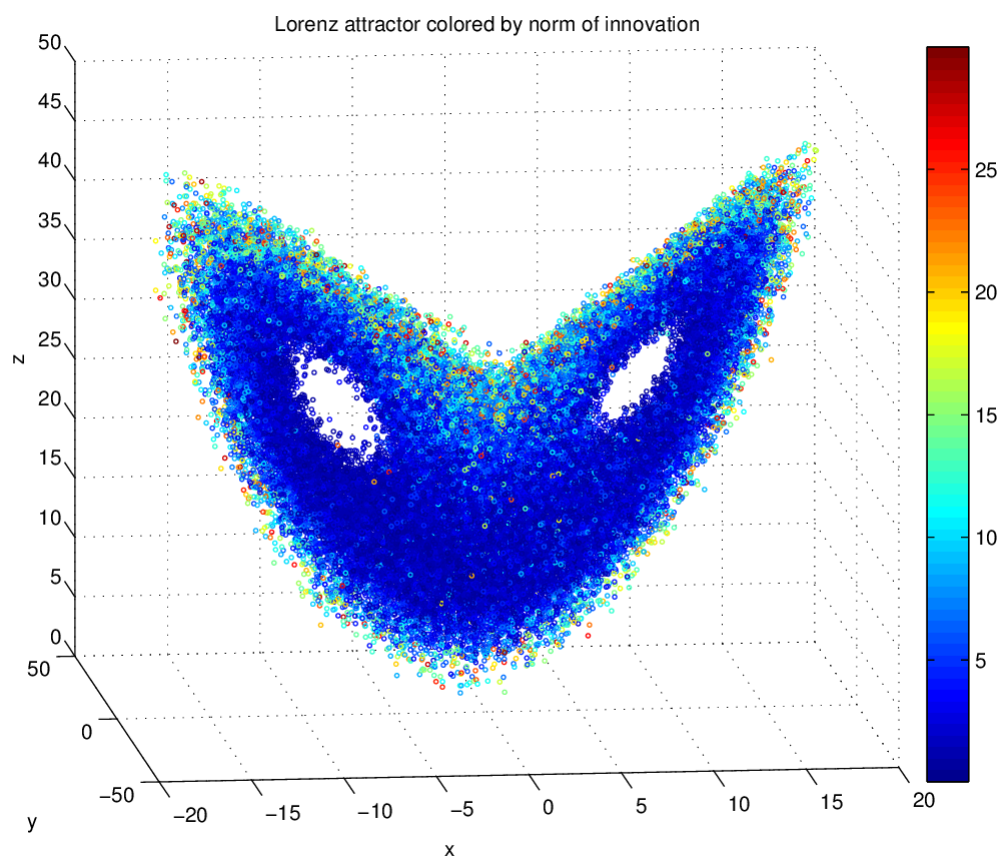}
  \caption[The Lorenz 63 attractor colored by the norm of the analysis increment]{
    The Lorenz 63 attractor colored by the norm of the analysis increment.
    Larger numbers represent where forecasts are poor, and the data assimilation algorithm puts more weight on the observation when initializing a new forecast.
    This figure was bitmapped for inclusion on the arXiv, for a full version visit {\protect \url{http://andyreagan.com}}.
  }
  \label{fig:lorenzattractor}
\end{figure}

The cannonical choice of $\sigma = 10, \beta = 8/3$ and $\rho = 28$, produce the well known butterfly attractor, is used for all examples here (Figure \ref{fig:lorenzattractor}).

\section{Basic derivation, OI}

We begin by deriving the equations for optimal interpolation (OI) and for the Kalman Filter with one scalar observation.
The goal of our derivation is to highlight which statistical information is useful, which is necessary, and what assumptions are generally made in the formulation of more complex algorithms.

Beginning with the one variable $U$, which could represent the horizontal velocity at one grid-point, say we make two independent observations of $U$: $U_1$ and $U_2$.
If the true state of $U$ is given by $U_t$, then the goal is to produce an analysis of $U_{1,2}$ which we call $U_a$, that is our best guess of $U_t$.

Denote the errors of $U_{1,2}$ by $\epsilon _{1,2}$ respectively, such that
$$ U_{1,2} = U_t + \epsilon _{1,2} .$$

For starters, we assume that we do not know $\epsilon _{1,2}$.
Next, assume that the errors vary around the truth, such that they are unbiased.
This is written as
\begin{equation*} \expv{U_{1,2}-U_t} = \expv{\epsilon _{1,2} } = 0. \end{equation*}

As of yet, our best guess of $U_t$ would be the average $(U_1 + U_2)/2$, and for any hope of knowing $U_t$ well we would need many unbiased observations at the same time.
This is unrealistic, and to move forward we assume that we know the variance of each $U_{1,2}$.
With this knowledge, which we write
\begin{equation} \expv{\epsilon _{1,2} ^2} = \sigma _{1,2} ^2 .\end{equation}
Further assume that these errors are uncorrelated
\begin{equation} \expv{\epsilon _1 \epsilon _1} = 0, \end{equation}
which in the case of observing the atmosphere, is hopeful.

We next attempt a linear fit of these observations to the truth, and equipped with the knowledge of observation variances we can outperform an average of the observations.
For $U_a$, our best guess of $U_t$, we write this linear fit as
\begin{equation} U_a = \alpha _1 U_1 + \alpha _2 U_2 \label{eq:analproblem}. \end{equation}
If $U_a$ is to be our best guess, we desire it to unbiased
\begin{equation*} \expv{U_a} = \expv{U_t} \Rightarrow \alpha _1 + \alpha _2 = 1\end{equation*}
and to have minimal error variance $\sigma _a ^2$.

At this point, note all of the assumptions that we have made: independent observations, unbiased instruments, knowledge of error variance, and uncorrelated observation errors.
Of these, knowing the error variances (and in multiple dimensions, error covariance matrices) will become our focus, and better approximations of these errors have been largely responsible for recent improvements in forecast skill.

Solving \ref{eq:analproblem}, by substituting $\alpha _2 = 1 - \alpha _1$, we have
\begin{align*} \sigma _a ^2 &= \expv{ \left ( \alpha_1 (U_1-U_t)+(1-\alpha_1)(U_2 - U_t) \right ) ^2 } \\
&= \alpha_1^2  \expv{ \epsilon _1 ^2 } + (1-\alpha _1^2) \expv{ \epsilon _2 ^2} + 2\alpha _1 (1- \alpha _1) \expv{\epsilon _1 \epsilon _2} \\
&=\alpha _1 ^2 \sigma_1 ^2 + (1-\alpha _1 ^2) \sigma _2^2\end{align*}
To minimize $\sigma _a^2$ with respect to $\alpha _1$ we compute $\frac{\partial }{\partial \alpha _1}$ of the above and solve for 0 to obtain
$$ \alpha _1 = \frac{\sigma _2 ^2}{2(\sigma_2^2 - \sigma _1^2)} $$
and similarly for $\alpha _2$.

We can also state the relationship for $\sigma _a ^2$ as a function of $\sigma_{1,2}^2$ by substituting our form for $\alpha _1$ into the second step of \ref{eq:analproblem} to obtain:
\begin{equation} \frac{1}{\sigma _a ^2} = \frac{1}{\sigma _1^2} + \frac{1}{\sigma_2^2} .\label{eq:analprec} \end{equation}
This last formula says that given accurate statistical knowledge of $U_{1,2}$, the precision of the analysis is the sum of the precision of the measurements.

\begin{figure}[t!]
  \centering
      \includegraphics[width=0.75\textwidth]{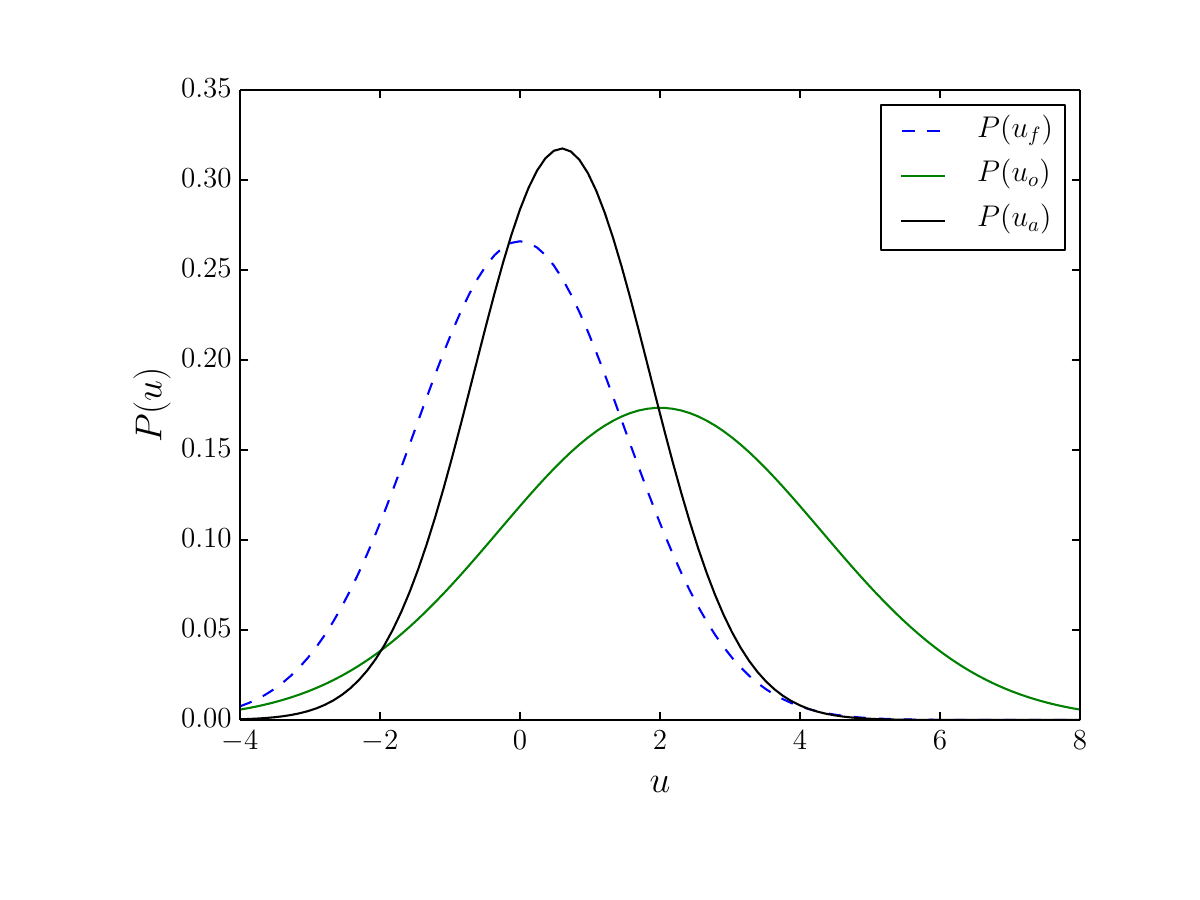}
  \caption[A 1-dimensional example of this data assimilation cycle]{
    A 1-dimensional example of this data assimilation cycle.
    First, a forecast is made with mean 0 and variance 1, labeled $u_f$.
    An observation denoted $u_o$ is the made with mean 2 and variance 2.
    From this information and Equation \ref{eq:1dda} we compute the analysis $u_f$ with mean 2/3, noting the greater reliance on the lower variance input $u_f$.
  }
  \label{fig:DAexample1D}
\end{figure}

We can now formulate the simplest data assimilation cycle by using our forecast model to generate a background forecast $U_b$ and setting $U_1 = U_b$.
We then assimilate the observations $U_o$ by setting $U_2 = U_o$ and solving for the analysis velocity $U_a$ as
\begin{equation} U_a = U_b + W (U_o - U_b) .\end{equation}

The difference $U_o - U_b$ is commonly reffered to as the ``innovation'' and represents the new knowledge from the observations with respect to the background forecast.
The optimal weight $W$ is then given by \ref{eq:analprec} as
\begin{equation} W = \frac{\sigma _b ^2}{\sigma_b ^2 - \sigma _o ^2} = \sigma_b ^2 (\sigma_b^2 - \sigma_o^2) ^{-1} \label{eq:1dda}.\end{equation}

The analysis error variance and precision are defined as before, and following Kalnay (2003) for convenience we can rewrite this in terms of $W$ as:
\begin{equation} \sigma _a ^2 = (1-W ) \sigma_b^2 .\end{equation}

We now generalize the above description for models of more than one variable to obtain the full equations for OI.
This was originally done by Aerology et al (1960), and independently by Gandin (1965). The following notation will be used for the OI and the other multivariate filters.
Consider a field of grid points, where there may be $\mbx_i = (p,p_{\text{rgh}},T,u_x,u_y,u_z)$ variables at each point for pressure, hydrostatic pressure, temperature, and three dimensions of velocity.
We combine all the variables into one list $\mbx$ of length $n$ where $n$ is the product of the number of grid points and the number of model variables.
This list is ordered by grid point and then by variable (e.g. $\mbx = (p_1,T_1,p_2,T_2,\ldots,p_{n/2},T_{n/2})$ for a field with just variables pressure $p$ and temperature $T$).
Denote the observations as $\mby_0$, where the ordering is the same as $\mbx$ but the dimension is much smaller (e.g. by a factor of $10^3$ for NWP) and the observations are not strictly at the locations of $\mbx$ or even direct measurements of the variables in $\mbx$.
In our experiment, temperature measurements are in the model variable $T$, but are understood as the average of the temperature of the grid points surrounding it within some chosen radius.
Examples of indirect and spatially heterogeneous observations in NWP are radar reflectivities, Doppler shifts, satellite radiances, and global positioning system atmospheric refractivities \shortcite{kalnay2003}.

We can now write the analysis with similar notation as for the scalar case with the addition of the observation operator $H$ which transforms model variables into observation space:
\begin{align} \mbx_a &= \mbx_b + \mbW \mbd \\
\mbd &=  \mby_0 - H (\mbx _b ) \\
\mbe _a  &= \mbx_t - \mbx_a . \end{align}

Note that $H$ is not required to be linear, but here we write $H$ as an $n_{\text{obs}} \times n$ matrix for a linear transformation.

Again, our goal is to obtain the weight matrix $\mbW$ that minimizes the analysis error $\mbe_a$. In Kalman Filtering, the weight matrix $\mbW$ is called the gain matrix $\mbK$.
First, we assume that the background and observations are unbiased.
\begin{align*} \expv{ \mbe_b} &= \expv{ \mbx _b} - \expv{\mbx _t } = 0\\
\expv{ \mbe_o }  &= \expv{ \mby_o } - \expv{ H(\mbx_t)} = 0 .\end{align*}

This a reasonable assumption, and if it does not hold true, we can account for this as a part of the analysis cycle \shortcite{dee1998data}.

Define the forecast error covariances as
\begin{align} \mathbf{P}_a &= \mathbf{A} = \expv{\mbe_a \mbe_a ^T}\\
\mathbf{P}_b &= \mbB = \expv{ \mbe_b \mbe_b ^T}\\
\mathbf{P}_o &= \mbR = \expv{ \mbe_o \mbe_o ^T} . \end{align}

Linearize the observation operator
\begin{equation} H(\mbx + \partial \mbx ) = H(\mbx) + \mbH \partial \mbx \end{equation}
where $h_{ij} = \partial H_i / \partial x_j$ are the entries of $\mbH$.
Since the background is typically assumed to be a reasonably good forecast, it is convenient to write our equations in terms of changes to $\mbx _b$.
The observational increment $\mbd$ is thus
\begin{align} \mbd &= \mby_o - H(\mbx_b) = \mby_o - H(\mbx_t + (\mbx_b-\mbx_t)) \\
&= \mby_o - H(\mbx _t) \mbH (\mbx_b -\mbx_t) = \mbe_o - \mbH \mbe_b . \end{align}

To solve for the weight matrix $\mbW$ assume that we know $\mbB$ and $\mbR$, and that their errors are uncorrelated:
\begin{equation} \expv{ \mbe_o \mbe_b ^T } = \expv{\mbe_b\mbe_o ^T} = 0 . \end{equation}

We now use the method of best linear unbiased estimation (BLUE) to derive $\mbW$ in the equation $\mbx_a - \mbx _b = \mbW \mbd$, which approximates $\mbx_a - \mbx_b = \mbW \mbd - \mbe_a$.
Since $\mbd = \mbe_o - \mbH \mbe_b$ and the well known BLUE solution to this simple equation, the $\mbW$ that minimizes $\mbe_a \mbe_a ^T$ is given by
\begin{align*} \mbW &= E\left [ ( \mbx _t - \mbx _b ) ( \mby _o - \mbH \mbx _b )^T \right] \left (E \left [ (\mby _o - \mbH \mbx _b)(\mby - \mbH \mbx _b)^T    \right ] \right ) ^{-1}\\
&= E\left [ (- \mbe _b ) (\mbe _o  - \mbH \mbe _b ) ^T \right ] \left ( E \left [ (\mbe _o - \mbH \mbe _b )( \mbe _o - \mbH \mbe _b ) ^T \right ] \right ) ^{-1} \\
&= E \left [ - \mbe _b \mbe _o ^T + \mbe _b \mbe _b ^T \mbH ^T \right ] \left ( E \left [ \mbe _o \mbe_o ^T - \mbH \mbe _b \mbe_o ^T - \mbe _o \mbe _b ^T \mbH ^T + \mbH \mbe _b \mbe _b ^T \mbH ^T \right ] \right ) ^{-1}\\
&= \left ( \cancel{E\left [ - \mbe _b \mbe _o ^T \right ]} + E \left [ \mbe _b \mbe _b ^T \right ] \mbH ^T \right ) \cdot \\
& ~~~~~~~\left ( E \left [ \mbe _o \mbe _o ^T \right ] - \cancel{\mbH E\left [ \mbe _b \mbe _o ^T \right ]} - \cancel{E \left [ \mbe _o \mbe _b ^T \right ] \mbH ^T } + \mbH E \left [ \mbe _b \mbe _b ^T \right ] \mbH ^T \right ) ^{-1}\\
&= \mbB \mbH ^T \left ( \mbR + \mbH \mbB \mbH ^T \right ) ^{-1}  \end{align*}

Now we solve for the analysis error covariance $\mathbf{P} _a $ using $\mbW$ from above and $\mbe _a = \mbe _b + \mbW \mbd$.
This is (with the final bit of sneakiness using $\mbW\mbR-\mbB\mbH^T=-\mbW\mbH\mbB\mbH^T$):
\begin{align*} \mathbf{P} _a &= \expv{\mbe_a \mbe_a^T} = \expv{(\mbe_b+ \mbW \mbd)(\mbe_b + \mbW \mbd)^T}\\
&= \expv{(\mbe_b + \mbW (\mbe _o - \mbH \mbe _b ) ) (\mbe_b + \mbW (\mbe _o - \mbH \mbe _b ) )^T}\\
&= \expv{(\mbe_b + \mbW (\mbe _o - \mbH \mbe _b ) ) (\mbe_b^T + (\mbe _O^T -  \mbe _b^T \mbH ^T )\mbW ^T )}\\
&= E \left [ \mbe_b\mbe_b^T + \mbe_b (\mbe _o^T -  \mbe _b^T \mbH ^T )\mbW ^T + \mbW (\mbe _o - \mbH \mbe _b )	\mbe _b ^T + \right.\\
&~~~~~~\left. \mbW (\mbe _o - \mbH \mbe _b )  (\mbe _O^T -  \mbe _b^T \mbH ^T )\mbW ^T \right ] \\
&= \expv{\mbe_b\mbe_b^T}  + \expv{\mbe_b (\mbe _o^T -  \mbe _b^T \mbH ^T )} \mbW ^T + \mbW \expv{ (\mbe _o - \mbH \mbe _b )  \mbe _b ^T} + \\
&~~~~~~\mbW \expv{(\mbe _o - \mbH \mbe _b )  (\mbe _O^T -  \mbe _b^T \mbH ^T )}\mbW ^T\\
&= \mbB - \mbB \mbH ^T \mbW^T - \mbW \mbH \mbB + \\
&~~~~~~\mbW \expv{\mbe_o\mbe_o^T - \cancel{\mbe_o\mbe_b^T}\mbH ^T - \cancel{\mbH \mbe _b \mbe _o ^T } + \mbH \mbe_b \mbe_b^T \mbH ^T }\mbW ^T\\
&= \mbB - \mbB \mbH ^T \mbW^T - \mbW \mbH \mbB +\mbW \mbR \mbW+  \mbW \mbH \mbB \mbH ^T \mbW ^T\\
&= \mbB - \mbW\mbH\mbB + (\mbW \mbR - \mbB \mbH^T ) \mbW^T + \mbW \mbH \mbB \mbH ^T \mbW ^T\\
&=(\mathbf{I} - \mbW \mbH ) \mbB \end{align*}

Note that in this formulation, we have assumed that $\mbR$ (the observation error covariance) accounts for both the instrumental error and the ``representativeness'' error, a term coined by Andrew Lorenc to describe the fact that the model does not account for subgrid-scale variability \shortcite{lorenc1981global}.
To summarize, the equations for OI are
\begin{align} \mbW &= \mbB \mbH ^T \left ( \mbR + \mbH \mbB \mbH ^T \right ) ^{-1}  \\
\mathbf{P} _a &=(\mathbf{I} - \mbW \mbH ) \mbB \end{align}

As a final note, it is also possible to derive the ``precision'' of the analysis $\mathbf{P} _a ^{-1}$ but in practice it is impractical to compute the inverse of the full matrix $\mathbf{P}_a$, so we omit this derivation.
In general is also not possible to solve for $\mbW$ for the whole model, and localization of the solution is achieved by considering each grid point and the observations within some set radius of this.
For different models, this equates to different assumptions on the correlations of observations, which we do not consider here.

\section{3D-Var}

Simply put, 3D-Var is the variational (cost-function) approach to finding the analysis.
It has been shown that 3D-var solves the same statistical problem as OI \shortcite{lorenc1986analysis}.
The usefulness of the variational approach comes from the computational efficiency, when solved with an iterative method.
This is why the 3D-Var analysis scheme became the most widely used data assimilation technique in operational weather forecasting, only to be succeeded by 4D-Var today.

The formulation of the cost function can be done in a Bayesian framework \shortcite{purser1984new} or using maximum likelihood approach (below).
First, it is not difficult to see that solving for $\partial J/ \partial U = 0$ for $U = U_a$ where we define $J(U)$ as
\begin{equation} J(U) = \frac{1}{2} \left [ \frac{(U-U_1)^2}{\sigma_1^2}+\frac{(U-U_2)^2}{\sigma_2^2} \right ] \end{equation}
gives the same solution as obtained through least squares in our scalar assimilation.
The formulation of this cost function obtained through the ML approach is simple \shortcite{edwards1984likelihood}.
Consider the likelihood of true value $U$ given an observation $U_1$ with standard deviation $\sigma_1$:
\begin{equation} L_{\sigma_1}(U|U_1)=p_{\sigma_1}(U_1|U)=\frac{1}{\sqrt{2\pi}\sigma_1}e^{-\frac{(T_1-T)^2}{2\sigma_1^2}} .\end{equation}
The most likely value of $U$ given two observations $U_{1,2}$ is the value which maximizes the joint likelihood (the product of the likelihoods).
Since the logarithm is increasing, the most likely value must also maximize the logarithm of the joint likelihood:
\begin{equation} \max_U \ln L_{\sigma_1,\sigma_2} (U|U_{1},U_{2}) = \max _U \left [ a - \frac{(U-U_1)^2}{\sigma_1^2}-\frac{(U-U_2)^2}{\sigma_2^2} \right ] \end{equation}
and this is equivalent to minimizing the cost function $J$ above.

Jumping straight the multi-dimensional analog of this, a result of Andrew Lorenc \shortcite{lorenc1986analysis}, we define the multivariate 3D-Var as finding the $\mbx _a$ that minimizes the cost function
\begin{equation} J(\mbx) = (\mbx - \mbx_b) ^T \mbB ^{-1} (\mbx - \mbx_b) + (\mby_o + H(\mbx))^T\mbR (\mby_o - H(\mbx)) .\end{equation}

The formal minimization of $J$ can be obtained by solving for $\nabla_x J(\mbx_a)$, and the result from Eugenia Kalnay is \shortcite{kalnay2003}:
\begin{equation} \mbx_a = \mbx_b+(\mbB^{-1} + \mbH^T\mbR^{-1}\mbH)^{-1}\mbH^T\mbR^{-1}(\mby_o-H(\mbx_b)) .\end{equation}

The principle difference between the Physical Space Analysis Scheme (PSAS) and 3D-Var is that the cost function $J$ is solved in observation space, because the formulation is so similar we do not derive it here.
If there are less observations than there are model variables, the PSAS is a more efficient method.

Finally we note for 3D-Var that operationally $\mbB$ can be estimated with some flow dependence with the National Meteorological Center (NMC) method, i.e.
\begin{equation} \mbB \approx \alpha \expv{(\mbx_f(48h)-\mbx_f(24h))(\mbx_f(48h)-\mbx_f(24h))^T} \end{equation}
with on the order of 50 different short range forecasts.
Up to this point, we have assumed knowledge of a static model error covariance matrix $\mbB$, or estimated it using the above NMC method.
In the following sections, we explore advanced methods that compute or estimate a changing model error covariance as part of the data assimilation cycle.

\section{Extended Kalman Filter}

The Extended Kalman Filter (EKF) is the ``gold standard'' of data assimilation methods \shortcite{kalnay2003}.
It is the Kalman Filter ``extended'' for application to a nonlinear model.
It works by updating our knowledge of error covariance matrix for each assimilation window, using a Tangent Linear Model (TLM), whereas in a traditional Kalman Filter the model itself can update the error covariance.
The TLM is precisely the model (written as a matrix) that transforms a perturbation at time $t$ to a perturbation at time $t+\Delta t$, analytically equivalent to the Jacobian of the model.
Using the notation of Kalnay \cite{kalnay2003}, this amounts to making a forecast with the nonlinear model $M$, and updating the error covariance matrix $\mbP$ with the TLM $L$, and adjoint model $L^T$

\begin{align*} \mbx^f (t_i) &= M _{i-1} [\mbx ^a (t_{i-1} ) ]\\
\mbP^f (t_i ) &= L_{i-1} \mbP^a (t_{i-1} ) L^T _{i-1} + \mathbf{Q} (t_{i-1} ) \end{align*}

where $\mathbf{Q}$ is the noise covariance matrix (model error).
In the experiments with Lorenz 63 presented in this section, $\mathbf{Q} = 0$ since our model is perfect.
In NWP, $\mathbf{Q}$ must be approximated, e.g. using statistical moments on the analysis increments \shortcite{danforth2007estimating,li2009accounting}.

The analysis step is then written as (for $H$ the observation operator):
\begin{align} \mbx^a (t_i ) &= \mbx^f (t_i) + \mbK_i \mbd_i\\
\mbP^a (t_i) &= (\mathbf{I} - \mbK_i \mbH_i )\mbP^f (t_i) \end{align}
where
\[ \mbd_i = \mby_i^o - \mbH[x^f (t_i) ] \]
is the innovation. The Kalman gain matrix is computed to minimize the analysis error covariance $P^a _i$ as
\[ \mbK_i = \mbP^f (t_i) \mbH_i ^T [ \mbR_i + \mbH_i \mbP^f (t_i) \mbH^T ] ^{-1} \]
where $\mbR_i$ is the observation error covariance.

Since we are making observations of the truth with random normal errors of standard deviation $\mbe$, the observational error covariance matrix $\mbR$ is a diagonal matrix with $\epsilon$ along the diagonal.

The most difficult, and most computationally expensive, part of the EKF is deriving and integrating the TLM.
For this reason, the EKF is not used operationally, and later we will turn to statistical approximations of the EKF using ensembles of model forecasts.
Next, we go into more detail regarding the TLM.

\subsection{Coding the TLM}

The TLM is the model which advances an initial perturbation $\delta \mbx_{i}$ at timestep $i$ to a final perturbation $\delta \mbx_{i+1}$ at timestep $i+1$.
The dynamical system we are interested in, Lorenz '63, is given as a system of ODE's:
\[ \frac{d\mbx}{dt} = F(\mbx) .\]
We integrate this system using a numerical scheme of our choice (in the given examples we use a second-order Runge-Kutta method), to obtain a model $M$ discretized in time.
\[ \mbx(t) = M[ \mbx(t_0) ] .\]
Introducing a small perturbation $\mby$, we can approximate our model $M$ applied to $\mbx(t_0) + \mby(t_0)$ with a Taylor series around $\mbx(t_0)$:
\begin{align*} M[ \mbx(t_0) + \mby(t_0) ] &= M [ \mbx(t_0) ] + \frac{\partial M}{\partial \mbx} \mby(t_0) + O [ \mby(t_0) ^2 ]\\ &\approx \mbx(t) + \frac{\partial M}{\partial \mbx} \mby(t_0) .\end{align*}
We can then solve for the linear evolution of the small perturbation $\mby(t_0)$ as
\begin{equation} \frac{d\mby }{dt } = \mathbf{J} \mby \label{eq:ODETLM} \end{equation}
where $\mathbf{J} = \partial F / \partial \mbx$ is the Jacobian of $F$.
We can solve the above system of linear ordinary differential equations using the same numerical scheme as we did for the nonlinear model.

\begin{figure}[h!]
  \centering
  \includegraphics[width=0.89\textwidth]{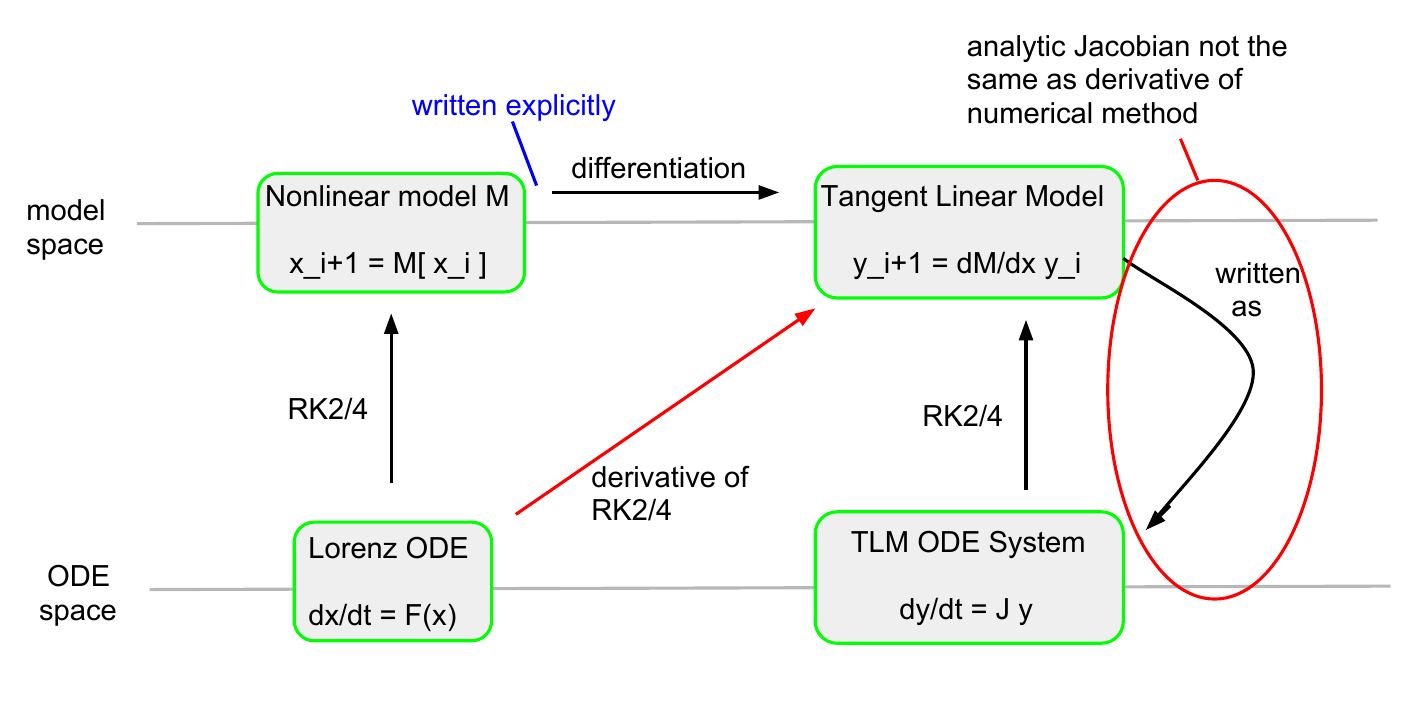}
  \caption[An explanation of how and why the best way to obtain a TLM is with a differentiated numerical scheme]{
    An explanation of how and why the best way to obtain a TLM is with a differentiated numerical scheme.
  }
  \label{fig:TLMscheme}
\end{figure}

One problem with solving the system of equations given by Equation \ref{eq:ODETLM} is that the Jacobian matrix of discretized code is not necessarily identical to the discretization of the Jacobian operator for the analytic system.
This is a problem because we need to have the TLM of our model $M$, which is the time-space discretization of the solution to $d\mbx/dt = F(\mbx)$.
We can apply our numerical method to the $d\mbx/dt = F(\mbx)$ to obtain $M$ explicitly, and then take the Jacobian of the result.
This method is, however, prohibitively costly, since Runge-Kutta methods are implicit.
It is therefore desirable to take the derivative of the numerical scheme directly, and apply this differentiated numerical scheme to the system of equations $F(\mbx)$ to obtain the TLM.
A schematic of this scenario is illustrated in Figure \ref{fig:TLMscheme}.
To that the derivative of numerical code for implementing the EKF on models larger than 3 dimensions (i.e. global weather models written in Fortan), automatic code differentiation is used \shortcite{autodiff1981}.

\begin{figure}[h!]
  \centering
  \includegraphics[width=0.79\textwidth]{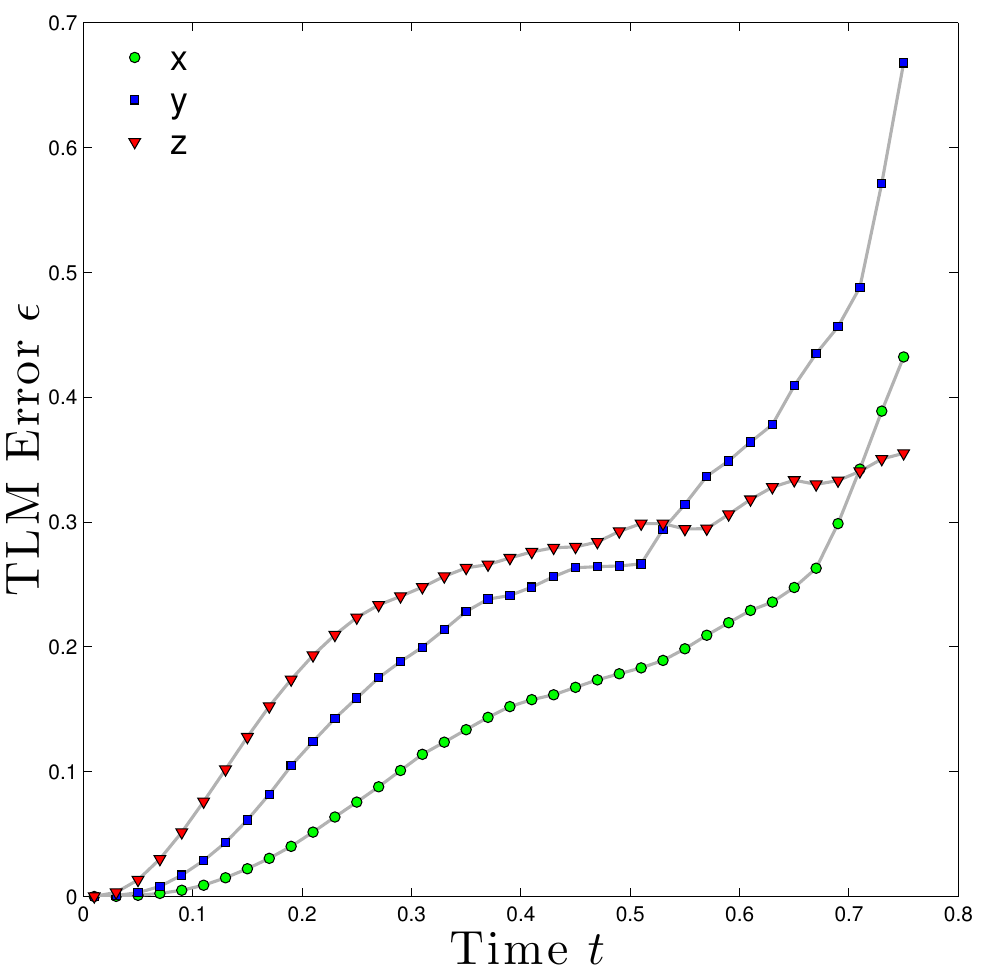}
  \caption[The future error predicted by the TLM is compared to the error growth in Lorenz 63 system for an initial perturbation with standard deviation of 0.1, averaged over 1000 TLM integrations]{
    The future error predicted by the TLM is compared to the error growth in Lorenz 63 system for an initial perturbation with standard deviation of 0.1, averaged over 1000 TLM integrations.
    The $\epsilon$ is not the error predicted by the TLM, but rather the error of the TLM in predicting the error growth.
  }
  \label{fig:TLMverification}
\end{figure}

To verify our implementation of the TLM, we propagate a small error in the Lorenz 63 system and plot the difference between that error and the TLM predicted error, for each variable (Figure \ref{fig:TLMverification}).

\section{Ensemble Kalman Filters}

The computational cost of the EKF is mitigated through the approximation of the error covariance matrix $\mbP_f$ from the model itself, without the use of a TLM.

One such approach is the use of a forecast ensemble, where a collection of models (ensemble members) are used to statistically sample model error propagation.
With ensemble members spanning the model analysis error space, the forecasts of these ensemble members are then used to estimate the model forecast error covariance.
In the limit of infinite ensemble size, ensemble methods are equivalent to the EKF \shortcite{evensen2003ensemble}.
There are two prevailing methodologies for maintaining independent ensemble members: (1) adding errors to the observation used by each forecast and assimilating each ensemble individually or (2) distributing ensemble initial conditions from the analysis prediction.
In the first method, the random error added to the observation can be truly random \shortcite{harris2011predicting} or be drawn from a distribution related to the analysis covariance.
The second method distributes ensemble members around the analysis forecast with a distribution drawn to sample the analysis covariance deterministically.
It has been show that the perturbed observation approach introduces additional sampling error that reduces the analysis covariance accuracy and increases the probability of of underestimating analysis error covariance \shortcite{hamill2001distance}.
For the examples given, we demonstrate both methods.

\begin{figure}[h!]
  \centering
  \includegraphics[width=0.79\textwidth]{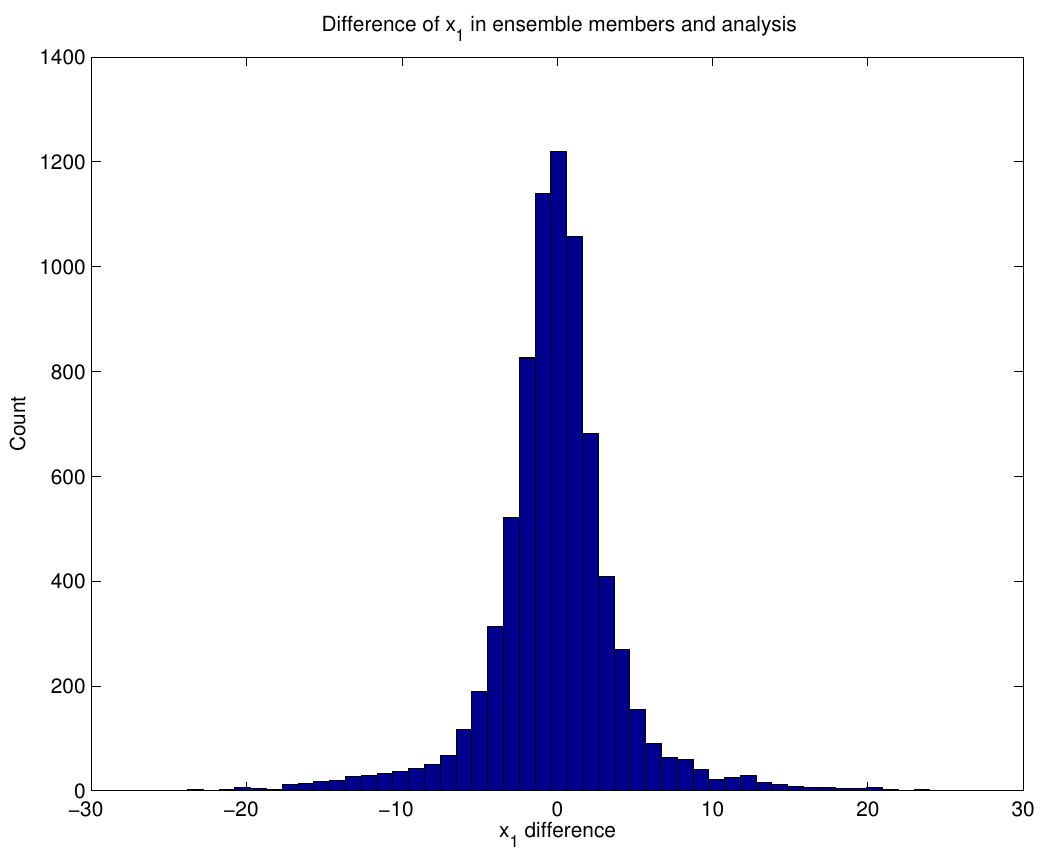}
  \caption[The difference of ensemble forecasts from the analysis is reported for 760 assimilation windows in one model run of length 200, with 10 ensemble members and an assimilation window of length 0.261]{
    The difference of ensemble forecasts from the analysis is reported for 760 assimilation windows in one model run of length 200, with 10 ensemble members and an assimilation window of length 0.261.
    This has the same shape of as the difference between ensemble forecasts and the mean of the forecasts (not shown).
    This spread of ensemble forecasts is what allows us to estimate the error covariance of the forecast model, and appears to be normally distributed.
  }
  \label{fig:EnKFhist}
\end{figure}

With a finite ensemble size this method is only an approximation and therefore in practice it often fails to capture the full spread of error.
To better capture the model variance, additive and multiplicative inflation factors are used to obtain a good estimate of the error covariance matrix (Section 2.6).
The spread of ensemble members in the $x_1$ variable of the Lorenz model, as distance from the analysis, can be seen in Figure \ref{fig:EnKFhist}.

The only difference between this approach and the EKF, in general, is that the forecast error covariance $\mbP^f$ is computed from the ensemble members, without the need for a tangent linear model:
\[ \mbP^f \approx \frac{1}{K-2} \sum _{k\neq l} \left ( \mbx_k ^f - \overline{\mbx} ^f _l \right ) \left (\mbx_k ^f - \overline{\mbx} ^f _l \right ) ^T .\]

In computing the error covariance $\mbP_f$ from the ensemble, we wish to add up the error covariance of each forecast with respect to the mean forecast.
But this would underestimate the error covariance, since the forecast we're comparing against was used in the ensemble average (to obtain the mean forecast).
Therefore, to compute the error covariance matrix for each forecast, that forecast itself is excluded from the ensemble average forecast.

\begin{figure}[h!]
  \centering
  \includegraphics[width=0.99\textwidth]{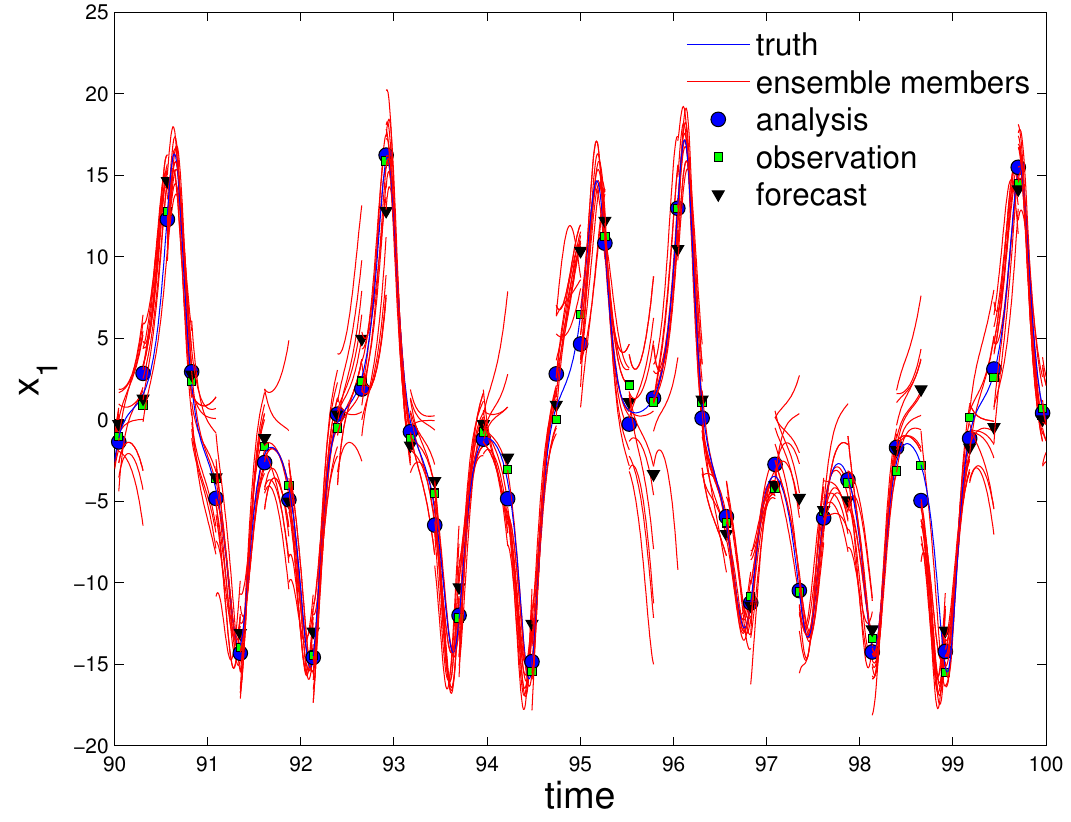}
  \caption[A sample time-series of the ensembles used in the EnKF]{
    A sample time-series of the ensembles used in the EnKF.
    In all tests, as seen here, 10 ensemble members are used.
    For this run, 384 assimilation cycles were performed with a window length of 0.26 model time units.
  }
  \label{fig:spaghetti}
\end{figure}

We can see the classic spaghetti of the ensemble with this filter implemented on Lorenz 63 in Figure \ref{fig:spaghetti}.

We denote the forecast within an ensemble filter as the average of the individual ensemble forecasts, and an explanation for this choice is substantiated by Burgers \shortcite{burgers1998analysis}.
The general EnKF which we use is most similar to that of Burgers.
Many algorithms based on the EnKF have been proposed and include the Ensemble Transform Kalman Filter (ETKF) \shortcite{ott2004local}, Ensemble Analysis Filter (EAF) \shortcite{anderson2001new}, Ensemble Square Root Filter (EnSRF) \shortcite{tippett2003ensemble}, Local Ensemble Kalman Filter (LEKF) \shortcite{ott2004local}, and the Local Ensemble Transform Kalman Filter (LETKF) \shortcite{hunt2007efficient}.
We explore some of these in the following sections.
A comprehensive overview through 2003 is provided by Evensen \shortcite{evensen2003ensemble}.

\subsection{Ensemble Transform Kalman Filter (ETKF)}

The ETKF introduced by Bishop is one type of square root filter, and we present it here to provide background for the formulation of the LETKF \shortcite{bishop2001adaptive}.
For a square root filter in general, we begin by writing the covariance matrices as the product of their matrix square roots.
Because $\mbP_a,\mbP_f$ are symmetric positive-definite (by definition) we can write
\begin{equation} \mbP_a = \mbZ_a \mbZ_a^T ~~,~~~ \mbP_f = \mbZ_f \mbZ_f^T \end{equation}
for $\mbZ_a,\mbZ_f$ the matrix square roots of $\mbP_a,\mbP_f$ respectively.
We are not concerned that this decomposition is not unique, and note that $\mbZ$ must have the same rank as $\mbP$ which will prove computationally advantageous.
The power of the SRF is now seen as we represent the columns of the matrix $\mbZ_f$ as the difference from the ensemble members from the ensemble mean, to avoid forming the full forecast covariance matrix $\mbP_f$.
The ensemble members are updated by applying the model $M$ to the states $\mbZ_f$ such that an update is performed by
\begin{equation} \mbZ_f = M \mbZ_a .\end{equation}

Per \shortcite{tippett2003ensemble}, the solution for $\mbZ_a$ can be solved by perturbed observations
\begin{equation} \mbZ_a = (\mathbf{I} -\mbK \mbH ) \mbZ _f + \mbK \mbW \end{equation}
where $\mbW$ is the observation perturbation of Gaussian random noise.
This gives the analysis with the expected statistics
\begin{align} \langle \mbZ_a \mbZ_a^T \rangle &= (\mathbf{I}-\mbK\mbH)\mbP_f(\mathbf{I}-\mbK\mbH)^T + \mbK\mbR\mbK^T\\
&=\mbP_a \end{align}
but as mentioned earlier, this has been shown to decrease the resulting covariance accuracy.

In the SRF approach, the ``Potter Method'' provides a deterministic update by rewriting
\begin{align} \mbP_a &= \mbZ_a \mbZ_a^T = \left(\mathbf{I}-\mbP_f\mbH^T(\mbR+\mbH\mbP_f\mbH^T)^{-1}\mbH\right)\mbP_f \\
&= \mbZ_f\left(\mathbf{I}-\mbZ_f\mbH^T(\mbH\mbZ_f\mbZ_f^T\mbH^T+\mbR)^{-1}\mbH\mbZ_f\right)\mbZ_f^T\\
&= \mbZ_f\left(\mathbf{I}-\mathbf{V}\mathbf{D}^{-1}\mathbf{Z}^T\right)\mbZ_f^T. \end{align}
We have defined $\mathbf{V}\equiv(\mbH\mbZ_f)^T$ and $\mathbf{D}\equiv\mathbf{V}^T\mathbf{V}+\mbR$.
Now for the ETFK step, we use the Sherman-Morrison-Woodbury identity to show
\begin{equation} \mathbf{I} - \mathbf{Z}\mathbf{D}^{-1}\mathbf{Z}^T = (\mathbf{I}+\mathbf{Z}_f\mathbf{H}^T\mathbf{R}^{-1}\mathbf{H}\mathbf{Z}_f)^{-1}.\end{equation}
The matrix on the right hand side is practical to compute because it is of dimension of the rank of the ensemble, assuming that the matrix inverse $\mbR^{-1}$ is available.
This is the case in our experiment, as the observation are uncorrelated.
The analysis update is thus
\begin{equation} \mbZ_a = \mbZ_f \mathbf{C}(\mathbf{\Gamma} +\mathbf{I})^{-1/2} ,\end{equation}
where $\mathbf{C}\mathbf{\Gamma}\mathbf{C}^T$ is the eigenvalue decomposition of $\mbZ_f^T\mbH^T\mbR^{-1}\mbH\mbZ_f$.
From \shortcite{tippett2003ensemble} we note that the matrix $\mathbf{C}$ of orthonormal vectors is not unique.

Then the analysis perturbation is calculated from
\begin{equation} \mbZ_a = \mbZ_f \mathbf{X} \mathbf{U}, \end{equation}
where $\mathbf{X}\mathbf{X}^T = (\mathbf{I}-\mathbf{Z}\mathbf{D}^{-1}\mathbf{V}^T)$ and $\mathbf{U}$ is an arbitrary orthogonal matrix.

In this filter, the analysis perturbations are assumed to be equal to the backgroud perturbations post-multiplied by a transformation matrix $\mathbf{T}$ so that that the analysis error covariance satisfies
$$ \mbP_a = (1- \mbK \mbH ) \mbP_f .$$

The analysis covariance is also written
$$ \mbP_f = \frac{1}{P-1} ( \mathbf{X}_a (\mathbf{X}_a) ^T = \mathbf{X} ^b \hat{\mathbf{A}} (\mathbf{X} ^b) ^T $$
where $\hat{\mathbf{A}} = [(P-1) \mathbf{I} + (\mathbf{H} \mathbf{X} ^b ) ^T \mathbf{R} ^{-1} (\mathbf{H} \mathbf{X} ^b ) ] ^{-1} $.
The analysis perturbations are $\mathbf{X} ^a = \mathbf{X} ^b \mathbf{T}$ where $\mathbf{T} = [(P-1)\hat{\mathbf{A}} ] ^{1/2} $.

To summarize, the steps for the ETKF are to (1) form $\mbZ_f^T\mbH^T\mbR^{-1}\mbH\mbZ_f$, assuming $\mathbf{R}^{-1}$ is easy, and (2) compute its eigenvalue decomposition, and apply it to $\mbZ_f$.

\subsection{Local Ensemble Kalman Filter (LEKF)}

The LEKF implements a strategy that becomes important for large simulations: localization.
Namely, the analysis is computed for each grid-point using only local observations, without the need to build matrices that represent the entire analysis space.
This localization removes long-distance correlations from $\mathbf{B}$ and allows greater flexibility in the global analysis by allowing different linear combinations of ensemble members at different spatial locations \shortcite{kalnay20074}.

Difficulties are presented however, in the form of computational grid dependence for the localization procedure.
The toroidal shape of our experiment, when discretized by OpenFOAM's mesher and grid numbering optimized for computational bandwidth for solving, does not require that adjacent grid points have adjacent indices.
In fact, they are decidedly non-adjacent for certain meshes, with every other index corresponding to opposing sides of the loops.
We overcome this difficulty by placing a more stringent requirement on the mesh of localization of indexing, and extracting the adjacency matrix for use in filtering.

The general formulation of the LEKF by Ott \shortcite{ott2004local} goes as follows, quoting directly:
\begin{enumerate}
\item Globally advance each ensemble member to the next analysis timestep. Steps 2-5 are performed for each grid point
\item Create local vectors from each ensemble member
\item Project that point's local vectors from each ensemble member into a low dimensional subspace as represented by perturbations from the mean
\item Perform the data assimilation step to obtain a local analysis mean and covariance
\item Generate local analysis ensemble of states
\item Form a new global analysis ensemble from all of the local analyses
\item Wash, rinse, and repeat.
\end{enumerate}

\begin{figure}[h!]
  \centering
  \includegraphics[width=0.79\textwidth]{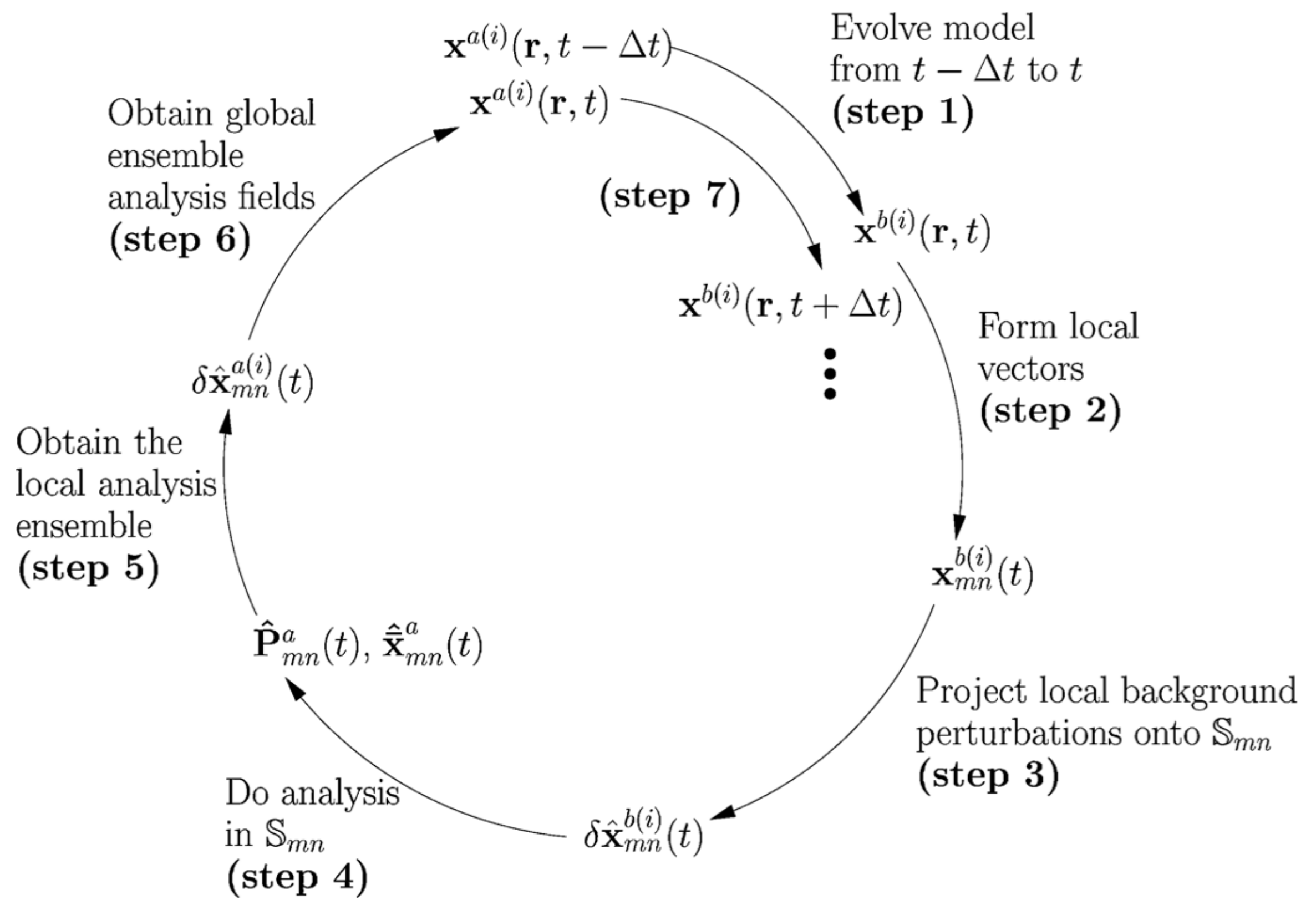}
  \caption[The algorithm for the LEKF as described by Ott (2004)]{
    The algorithm for the LEKF as described by Ott (2004).
    }
  \label{fig:ott2004algorithm}
\end{figure}

Figure \ref{fig:ott2004algorithm} depicts this cycle in greater detail, with similar notations to that used thus far \shortcite{ott2004local}.

In the reconstructed global analysis from local analyses, since each of the assimilations were performed locally, we cannot be sure that there is any smoothness, or physical balance, between adjacent points.
To ensure a unique solution of local analysis perturbations that respects this balance, the prior background forecasts are used to constrain the construction of local analysis perturbations to minimize the difference from the background.

\subsection{Local Ensemble Transform Kalman Filter (LETKF)}

Proposed by Hunt in 2007 with the stated objective of computational efficiency, the LETKF is named from its most similar algorithms from which it draws \shortcite{hunt2007efficient}.
With the formulation of the LEKF and the ETKF given, the LETKF can be  described as a synthesis of the advantages of both of these approaches.
This is the method that was sufficiently efficient for implementation on the full OpenFOAM CFD model of 240,000 model variables, and so we present it in more detail and follow the notation of Hunt et al (2007).
As in the LEKF, we explicitly perform the analysis for each grid point of the model.
The choice of observations to use for each grid point can be selected a priori, and tuned adaptively.
The simplest choice is to define a radius of influence for each point, and use observations within that.
In addition to this simple scheme, Hunt proposes weighing observations by their proximity to smoothly transition across them.
The goal is to avoid sharp boundaries between neighboring points.

Regarding the LETKF analysis problem as the construction of an analysis ensemble from the background ensemble, where we specify that the mean of the analysis ensemble minimizes a cost function, we state the equations in terms of a cost function.
The LETKF is thus
\begin{equation} \left\{ \mbx _{b(i)}\,:\, i=1..k \right\} \underrightarrow{~~~~~~~\text{LETKF}~~~~~~} \left\{ \mbx _{a(i)}\,:\,i=1..k\right\} \end{equation}
where
\begin{align*} J(\mbx) = &(\mbx - \overline{\mbx}_b) ^T (\mbP_b)^T(\mbx - \overline{\mbx}_b)\\
&+\left(\mby_o-H(\mbx)\right)^T \mbR^{-1} \left(\mby_o-H(\mbx)\right). \end{align*}
Since we have defined $\mbP_b = (1/(k-1))\mbX_b \mbX_b^T$ for the $m\times k$ matrix $\mbX_b = (\mbx_{b(1)},\ldots,\mbx_{b(k)})$, the rank of $\mbP_b$ is at most $k-1$, and the inverse of $\mbP$ is not defined.
To ameliorate this, we consider the analysis computation to take place in the subspace $S$ spanned by the columns of $\mbX_b$ (the ensemble states), and regard $\mbX_b$ as a linear transformation from a $k$-dimensional space $\tilde{S}$ into $S$.
Letting $\mbw$ be a vector in $\tilde{S}$, then $\mbX_b\mbw$ belongs to $S$ and $\mbx = \overline{\mbx}_b + \mbX_b \mbw$.
We can therefore write the cost function as
\begin{align*} \tilde{J} = &\frac{1}{k-1}\mbw^T\mbw + (\mby_o - H(\overline{\mbx}_b+\mbX_b\mbw))^T\\
&\times \mbR^{-1}\left(\mby_o - H(\overline{\mbx}_b+\mbX_b\mbw)\right)\end{align*}
and it has been shown that this is an equivalent formulation by \shortcite{hunt2007efficient}.

Finally, we linearize $H$ by applying it to each $\mbx_{b(i)}$ and interpolating.
This is accomplished by defining
\begin{equation} \mby _{b(i)} = H( \mbx_{b(i)}) \end{equation}
such that
\begin{equation} H(\overline{\mbx}_b - \mbX_b \mbw) = \overline{\mby}_b - \mbY _b \mbw \end{equation}
where $\mbY$ is the $l\times k$ matrix whose columns are $\mby_{b(i)} - \overline{\mby}_b$ for $\overline{\mby}_b$ the average of the $\mby_{b(i)}$ over $i$.
Now we rewrite the cost function for the linearize observation operator as
\begin{align*} \tilde{J}^* = &(k-1)\mbw^T \mbw + (\mby _o -\overline{\mby} _b - \mbY _b \mbw )^T\\
&\times\mbR^{-1} (\mby _o -\overline{\mby} _b - \mbY _b \mbw ) \end{align*}
This cost function minimum is the form of the Kalman Filter as we have written it before:
\begin{align*} \overline{\mbw}_a &= \tilde{\mbP}_a (\mbY _b )^T \mbR ^{-1} (\mby_o - \overline{\mby}_b),\\
\tilde{P}_a &= \left( (k-1)\mathbf{I}+(\mbY_b)^T\mbR^{-1}\mbY_b\right ) ^{-1} .\end{align*}
In model space these are
\begin{align*} \tilde{\mbx}_a &= \overline{\mbx} _b + \mbX _b \overline{\mbw}_a \\
\mbP _a &= \mbX _b \tilde{\mbP}_a (\mbX _b)^T . \end{align*}

We move directly to a detailed description of the computational implementation of this filter.
Starting with a collection of background forecast vectors $\{ \mbx_{b(i)} \,:\,i=1..k \}$, we perform steps 1 and 2 in a global variable space, then steps 3-8 for each grid point:

\begin{enumerate}
\item apply $H$ to $\mbx_{b(i)}$ to form $\mby_{b(i)}$, average the $\mby_b$ for $\overline{\mby_b}$, and form $\mbY b$.
\item similarly form $\mbX_b$. (now for each grid point)
\item form the local vectors
\item compute $\mathbf{C}=(\mbY_b)^T\mbR^{-1}$ (perhaps by solving $\mbR \mathbf{C}^T = \mbY_b$
\item compute $\tilde{\mbP}_a = \left( (k-1)\mathbf{I} / \rho + \mathbf{C} \mbY _b \right ) ^{-1}$ where $\rho > 1$ is a tunable covariance inflation factor
\item compute $\mbW_a = \left ( (k-1) \tilde{\mbP} _a \right ) ^{1/2}$
\item compute $\overline{\mbw} _a  = \tilde{\mbP}_a \mathbf{C} \left ( \mby_o - \overline{\mby} _b \right )$ and add it to the column of $\mbW_a$
\item multiply $\mbX_b$ by each $\mbw_{a(i)}$ and add $\tilde{\mbx}_b$ to get $\left\{ \mbx_{a(i)}\,:\,i=1..k\right \}$ to complete each grid point
\item combine all of the local analysis into the global analysis
\end{enumerate}

The greatest implementation difficulty in the OpenFOAM model comes back to the spatial discretization and defining the local vectors for each grid point.

\subsubsection{Four Dimensional Local Ensemble Transform Kalman Filter (4D-LETKF)}

We expect that extending the LETKF to incorporate observations from within the assimilation time window will become useful in the context of our target experiment high-frequency temperature sensors.
This is implemented without great difficulty on top of the LETKF, by projecting observations as they are made into the space $\tilde{S}$ with the observation time-specific matrix $H_\tau$ and concatenating each successive observation onto $\mby_o$.
Then with the background forecast at the time $\tau$, we compute $( \overline{\mby}_b )_\tau$ and $ (\mbY _b ) _\tau$, concatenating onto $\overline{\mby}_b$ and $\mbY_b$, respectively.

\section{4D-Var}

The formulation of 4D-Var aims to incorporate observations at the time they are observed.
This inclusion of observations has resulted in significant improvements in forecast skill, cite.
The next generation of supercomputers that will be used in numerical weather prediction will allow for increased model resolution.
The main issues with the 4D-Var approach include attempting to parallelize what is otherwise a sequential algorithm, and the computational difficulty of the strong-constraint formulation with many TLM integrations, with work being done by Fisher and Andersson (2001).
The formulation of 4D-Var is not so different from 3D-Var, with the inclusion of a TLM and adjoint model to assimilation observations made within the assimilation window.
At the ECMWF, 4D-Var is run on a long-window cycle (24 and 48 hours) so that many observations are made within the time window.
Generally, there are two main 4D-Var formulations in use: strong-constraint and weak-constaint.
The strong-constraint formulation uses the forecast model as a strong constraint, and is therefore easier to solve.
From Kalnay (2007) we have the standard weak-constaint formulation with superscript denoting the timestep within the assimilation window (say there are $0,\ldots,N$ observation time steps within the assimilation window):
\begin{align*} J(\mbx^0) = &\frac{1}{2} (\mbx^0 - \mbx^0_b) ^T \mbB ^{-1} (\mbx^0 - \mbx_b^0 \frac{1}{2} \sum _{i=0}^N \left( H^i(\mbx^i)-\mby^i\right)^T\\
&(\mbR^i) ^{-1} \left( H^i(\mbx^i)-\mby^i \right ) .\end{align*}

The weak-constraint formulation has only recently been implemented operationally.
This is \shortcite{kalnay20074}:
\begin{align*} J(\mbx^0 , \beta) = &\frac{1}{2} (\mbx^0 - \mbx^0_b)^T \mbB^{-1} (\mbx^0 - \mbx^0_b) + \frac{1}{2} \sum _{i=1} ^N\\
& \left(H^i(\mbx^i+\beta)-\mby^i\right)^T (\mbR^i)^{-1} \left(H^i(\mbx^i+\beta)-\mby^i\right) \\
&+\frac{1}{2}\beta^T\mathbf{Q}^{-1}\beta = J_b +J_o+J_\beta.\end{align*}

Future directions of operational data assimilation combine what is best about 4D-Var and ensemble-based algorithms, and these methods are reffered to as ensemble-variational hybrids.
We have already discussed a few such related approaches: the use of the Kalman gain matrix $\mbK$ used within the 3D-Var cost function, the ``NCEP'' method for 3D-Var, and inflating the ETKF covariance with the 3D-Var static (climatological) covariance $\mbB$.

\subsection{Ensemble-variational hybrids}

Ensemble-variational hybrid approaches benefit from flow dependent covariance estimates obtained from an ensemble method, and the computational efficiency of iteratively solving a cost function to obtain the analysis, are combined.
An ensemble of forecasts is used to update the static error covariance $\mbB$ in the formulation of 4D-Var (this could be done in 3D-Var as well) with some weighting from the climatological covariance $\mbB$ and the ensemble covariance.
This amounts to a formulation of
\begin{equation} \mbB = \alpha \mbP_f + (1-\alpha) \mbB \end{equation}
for some $\alpha$ between 0 and 1.
Combining the static error covariance $\mbB$ and the flow dependent error covariance $\mbP_f$ is an active area of research \shortcite{clayton2012operational}.

\section{Performance with Lorenz 63}

\begin{figure}[h!]
  \centering
  \includegraphics[width=0.79\textwidth]{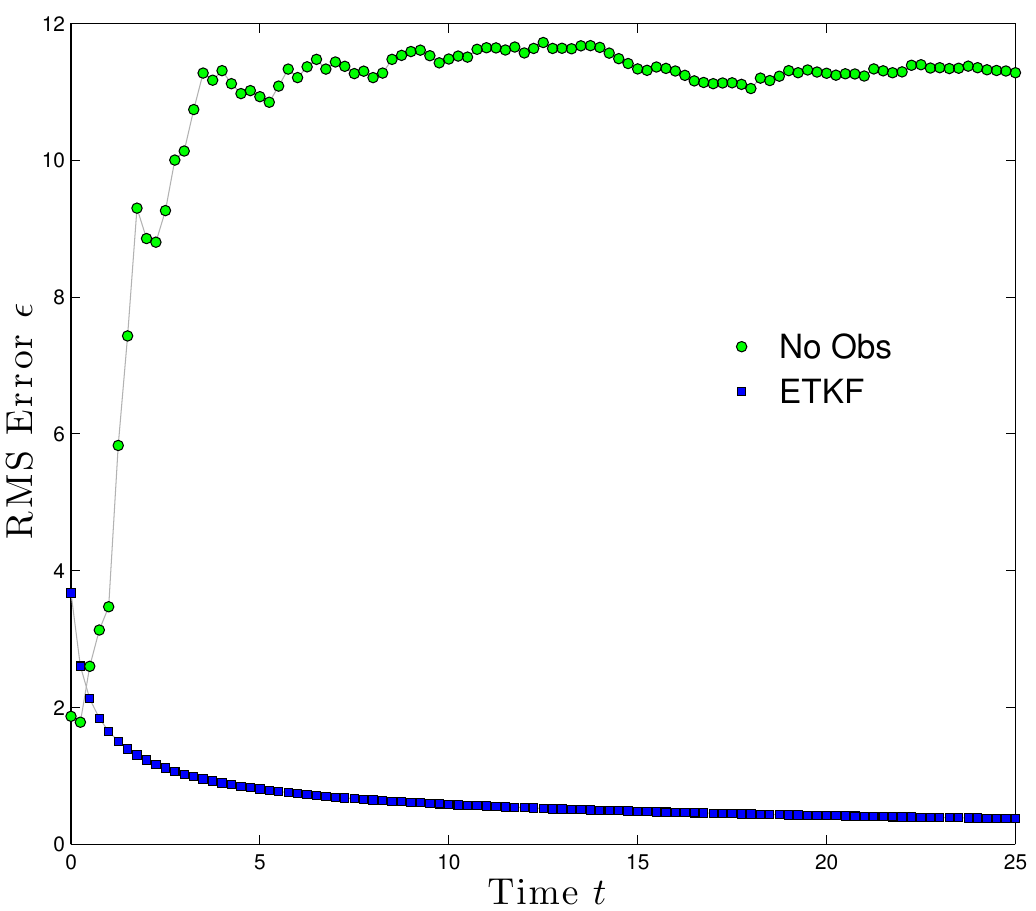}
  \caption[For near-perfect initial knowledge of the atmospheric state, the running average RMS error of a prediction using DA is compared against a prediction that does not]{
    For near-perfect initial knowledge of the atmospheric state, the running average RMS error of a prediction using DA is compared against a prediction that does not.
    Initial error for DA is normal with mean 0 and variance 0.01, and the RMS errors reported are averages of 3 runs.
    Observational noise is normally distributed with mean 0 and variance 0.05 (approximately 5\% of climatological variation), and an assimilation window length of 0.25 is used.
    As would be expected, in the long run DA greatly improves forecast skill, as forecasts without DA saturate to climatological error.
    The initial period where a forecast exceeds the skill of DA is potentially due to the spin-up time required by the filter to dial in the analysis error covariance.
  }
  \label{fig:DA_test}
\end{figure}

With the EKF, EnKF, ETKF, and EnSRF algorithms coded as above, we present some basic results on their performance.
First, we see that using a data assimilation algorithm is highly advantageous in general (Figure \ref{fig:DA_test}).
Of course we could think to use the observations directly when DA is not used, but since we only observing the $x_1$ variable (and in real situations measurements are not so direct), this is not prudent.

Next, we test the performance of each filter against increasing assimilation window length.
As we would expect, longer windows make forecasting more challenging for both algorithms.
The ensemble filter begins to out-perform the extended Kalman filter for longer window lengths because the ensemble more accurately captures the nonlinearity of the model.

\begin{figure}[h!]
  \centering
  \includegraphics[width=0.79\textwidth]{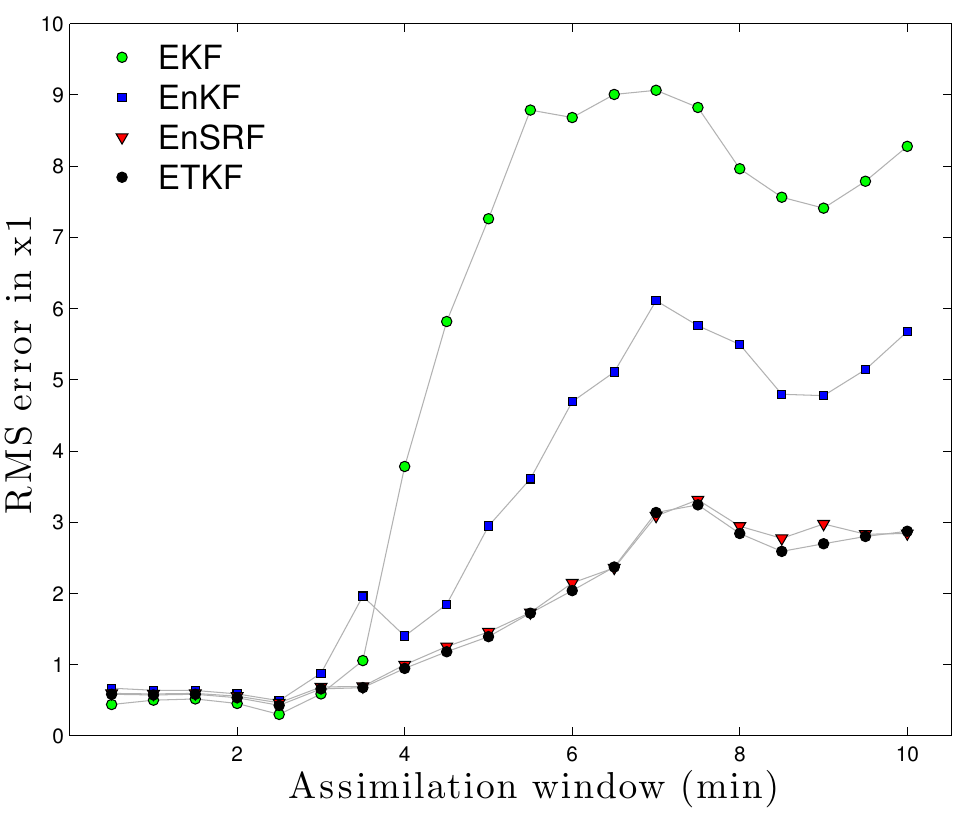}
  \caption[The RMS error is reported for our EKF and EnKF filters]{
    The RMS error (not scaled by climatology) is reported for our EKF and EnKF filters, measured as the difference between forecast and truth at the end of an assimiliation window for the latter 2500 assimiliation windows in a 3000 assimilation window model run.
    Error is measured in the only observed variable, $x_1$.
    Increasing the assimilation window led to an decrease in predictive skill, as expected.
    Additive and multiplicative covariance inflation is the same as Harris et. al (2011).
  }
  \label{fig:window_test}
\end{figure}

\begin{figure}[h!]
  \centering
  \includegraphics[width=0.79\textwidth]{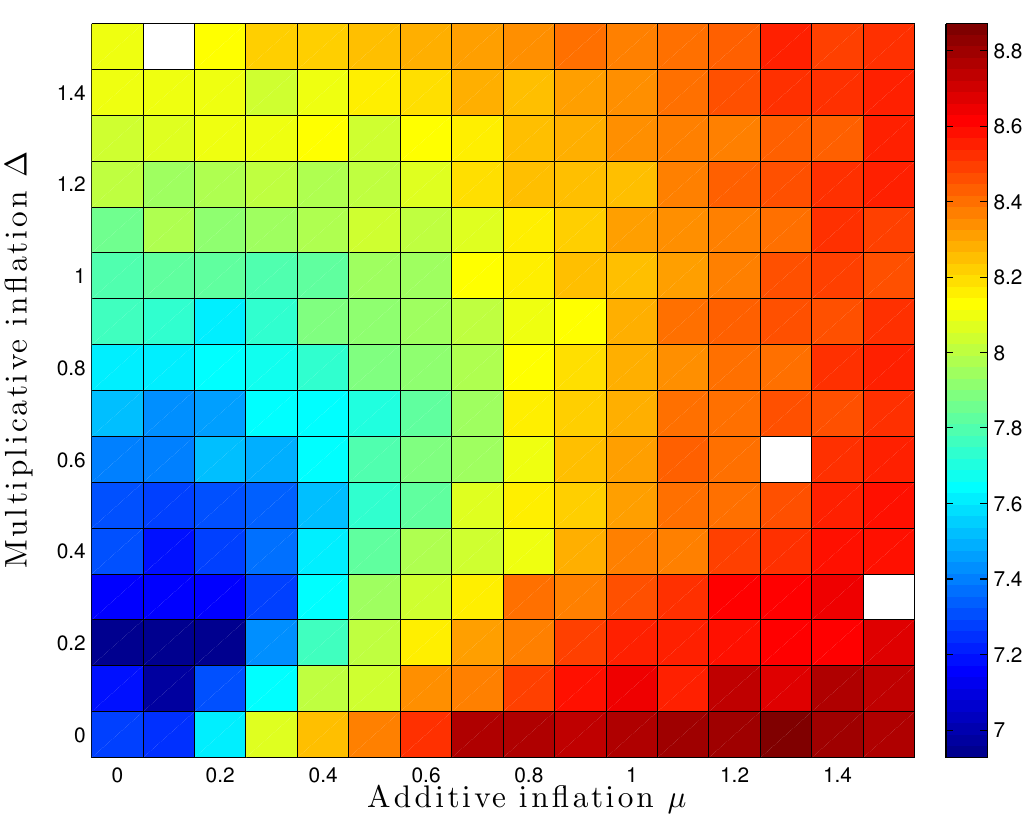}
  \caption[The RMS error averaged over 100 model runs of length 1000 windows is reported for the ETKF for varying additive and multiplicative inflation factors]{
    The RMS error averaged over 100 model runs of length 1000 windows is reported for the ETKF for varying additive and multiplicative inflation factors $\Delta$ and $\mu$.
    Each of the 100 model runs starts with a random IC, and the analysis forecast starts randomly.
    The window length here is 390 seconds.
    The filter performance RMS is computed as the RMS value of the difference between forecast and truth at the assimilation window for the latter 500 windows, allowing a spin-up of 500 windows.
  }
  \label{fig:ETKF_cov_tuning_390s}

\end{figure}

Each of the filters is sensitive to the size of the covariance matrix, and in some cases to maintain the numerical stability of the assimilation it is necessary to use additive and/or multiplicative covariance inflation.
In all filters multiplicative error is performed with inflation factor $\Delta$ by
\begin{equation} \mbK \leftarrow (1+\Delta) \mbK .\end{equation}
Additive inflation is performed for the EKF as in Yang et al (2006) and Harris et al (2011) where uniformly distributed random numbers in $[0,\mu]$ are added to the diagonal.
This is accomplished in MATLAB by drawing a vector  $\nu$ of length $n$ of uniformly distributed numbers in $[0,1]$, multiplying it by $\mu$, and adding it to the diagonal:
\begin{equation} \mbK \leftarrow \mbK + \text{diag}(\mu\cdot \nu).\end{equation}
For the EnKF we perform additive covariance inflation to the analysis state itself for each ensemble:
\begin{equation} \mbx \leftarrow \mbx + \nu .\end{equation}
The result of searching over these parameters for the optimal values to minimize the analysis RMS error differs for each window length, so here here we present the result for a window of length 390 seconds for the ETKF.
In Figure \ref{fig:ETKF_cov_tuning_390s} we see that there is a clear region of parameter space for which the filter performs with lower RMS, although this region appears to extend into values of $\Delta > 3$.

\chapter{Thermal Convection Loop Experiment}
\chaptermark{Thermosyphon}
\begin{quote}
The thermosyphon, a type of natural convection loop or non-mechanical heat pump, can be likened to a toy model of climate \shortcite{harris2011predicting}.
In this chapter we present both physical and computationally simulated thermosyphons, the equations governing their behavior, and results of the computational experiment.
\end{quote}

\section{Physical Thermal Convection Loops}

\begin{figure}[t!]
  \centering
  \includegraphics[width=0.79\textwidth]{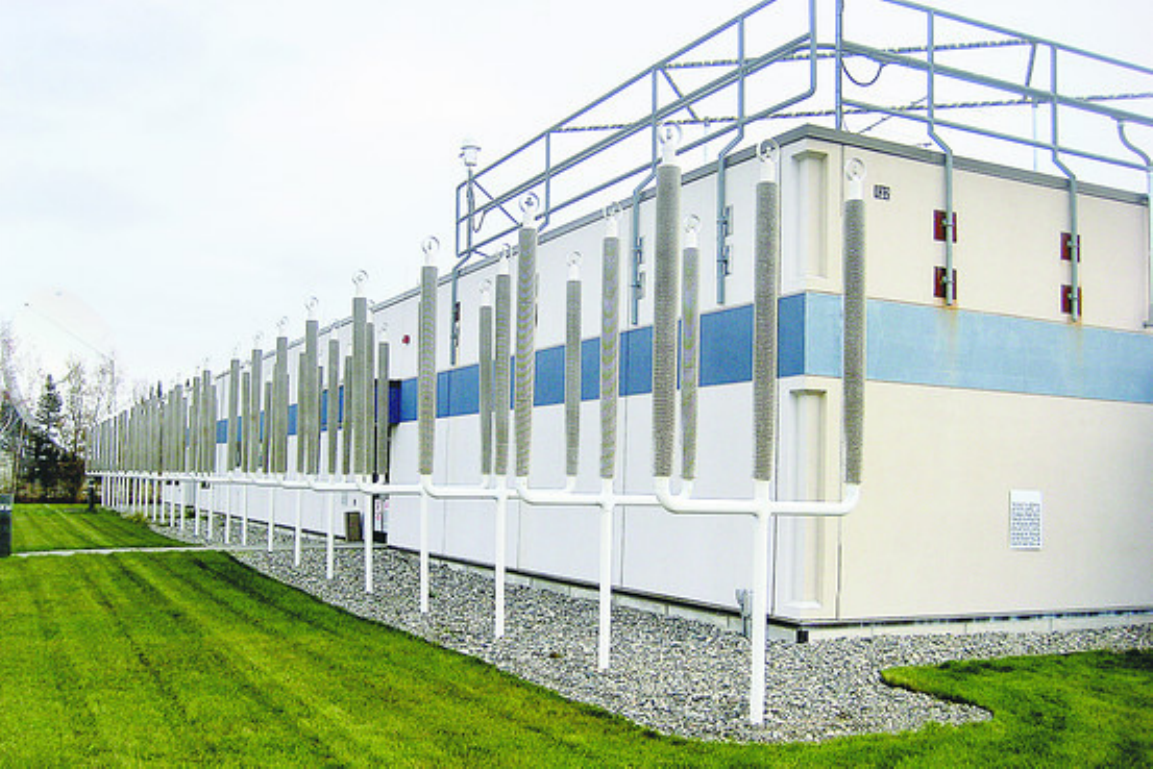}
  \caption[Thermosyphons flank a Federal Aviation Administration building in Fairbanks, Alaska]{
  Thermosyphons flank a Federal Aviation Administration building in Fairbanks, Alaska.
  The devices help keep the permafrost frozen beneath buildings and infrastructure by transferring heat out of the ground.
  Photo credit Jim Carlton {\protect\shortcite{WSJthermo}}.
  }
  \label{fig:fairbanks}
\end{figure}

Under certain conditions, thermal convection loops operate as non-mechanical heat pumps, making them useful in many applications.
In Figures \ref{fig:thermosyphons} and \ref{fig:fairbanks} we see four such applications.
Solar hot water heaters warm water using the sun's energy, resulting in a substantial energy savings in solar-rich areas \shortcite{WSJthermo}.
Geothermal systems heat and cool buildings with the Earth's warmth deep underground, and due to low temperature differences are often augmented with a physical pump.
CPU coolers allow modern computers to operate at high clock rates, while reducing the noise and energy consumption of fans \shortcite{CPUThermo}.
Ground thermosyphons transfer heat away from buildings built on fragile permafrost, which would otherwise melt in response to a warming climate, and destroy the building \shortcite{Xu2008experimental}.
The author notes that these systems rely on natural convection, and when augmented with pumps the convection is mixed or forced, and would necessitate different models.

\begin{figure}[t!]
  \centering
  \includegraphics[width=0.89\textwidth]{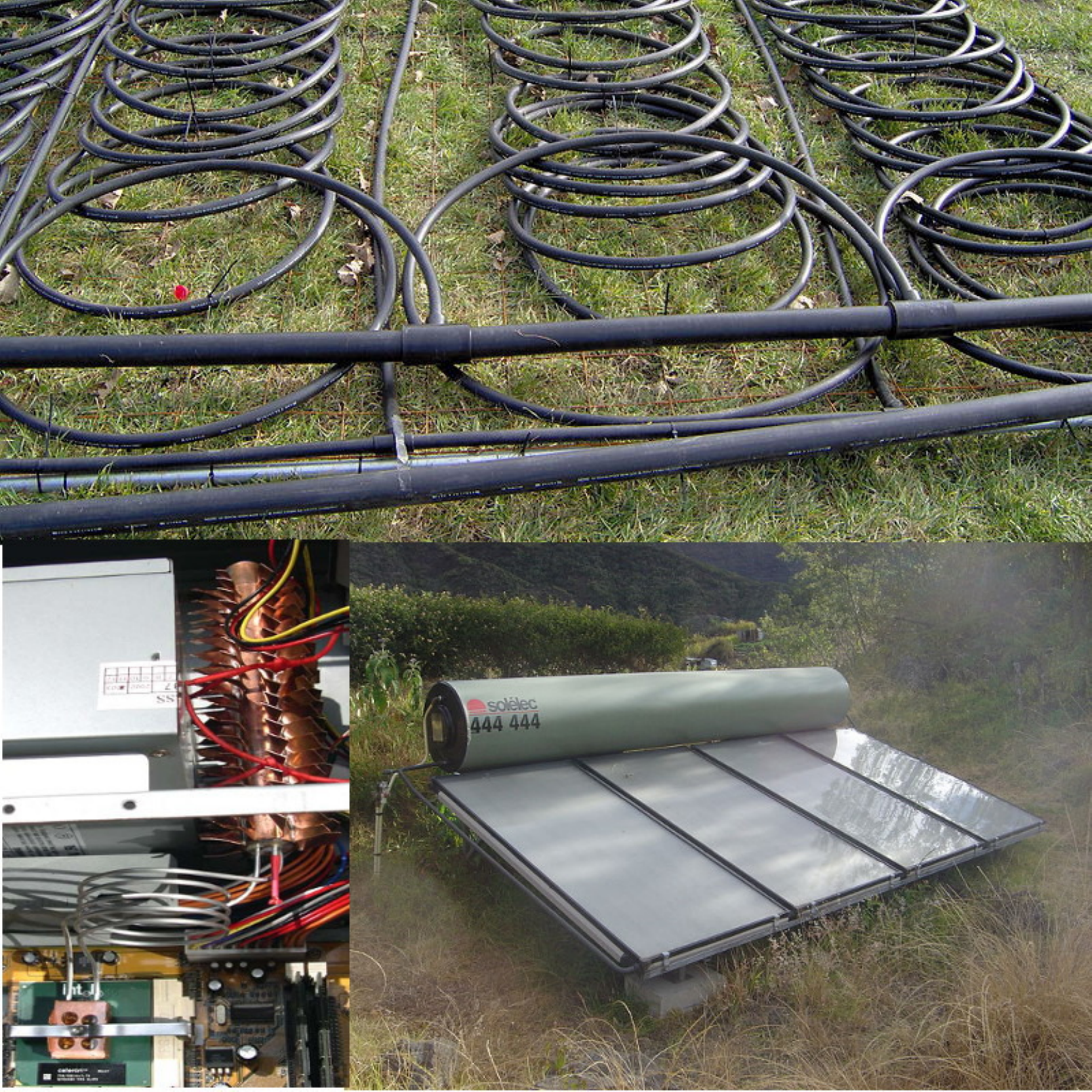}
  \caption[Here we see many different thermosyphon applications]{
    Here we see many different thermosyphon applications.
    Clockwise starting at the top, we have a geothermal heat exchanger, solar hot water heater, and computer cooler.
    Images credit to BTF Solar and Kuemel Bernhard {\protect\shortcite{CPUThermo,WSJthermo}}.
  }
  \label{fig:thermosyphons}
\end{figure}

\subsection{Ehrhard-M\"{u}ller Equations}

\begin{figure}[t!]
  \centering
  \includegraphics[width=0.79\textwidth]{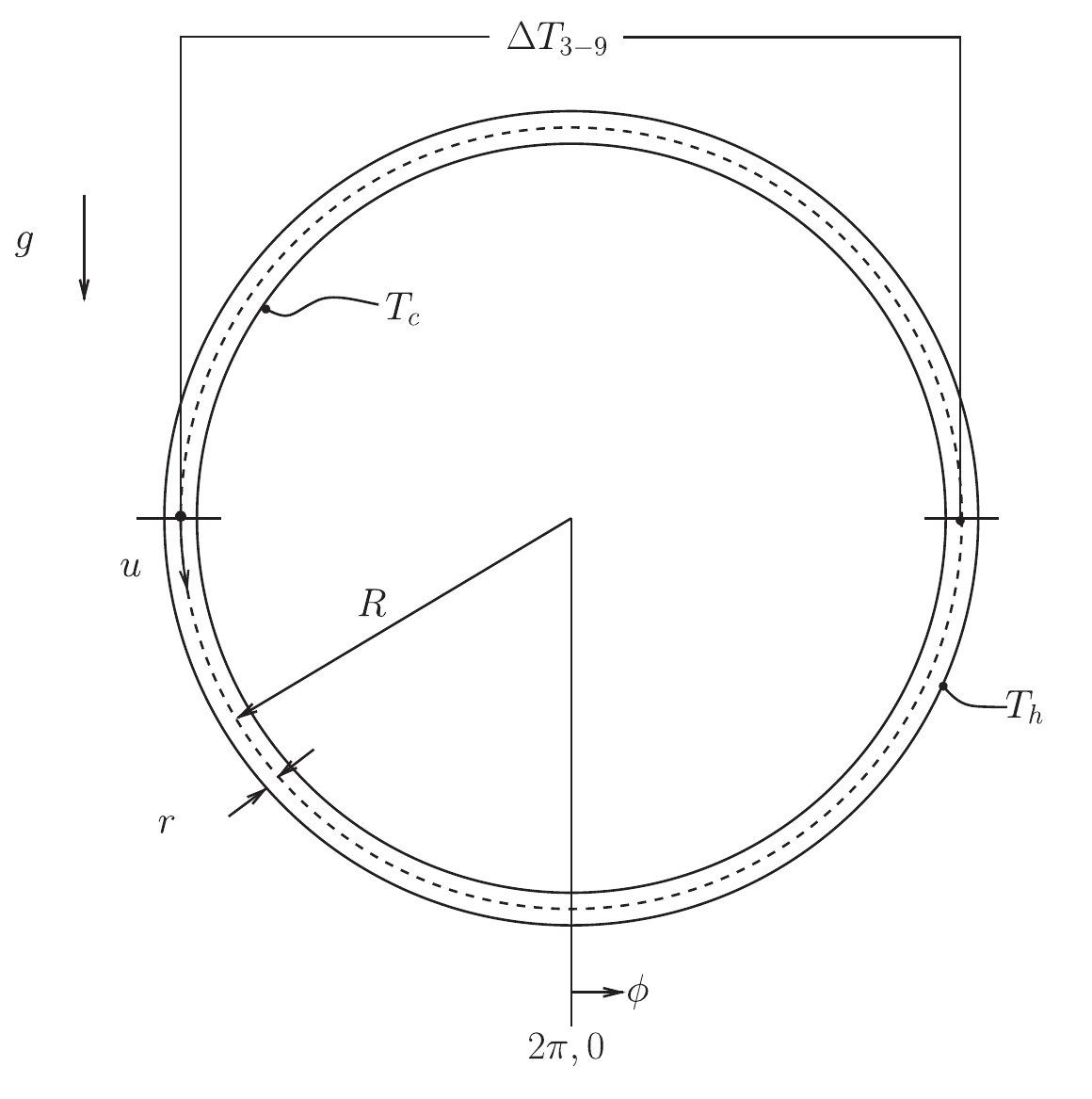}
  \caption[Schematic of the experimental, and computational, setup from Harris et al (2012)]{
    Schematic of the experimental, and computational, setup from Harris et al (2012).
    The loop radius is given by $R$ and inner radius by $r$.
    The top temperature is labeled $T_c$ and bottom temperature $T_h$, gravity $g$ is defined downward, the angle $\phi$ is prescribed from the 6 o'clock position, and temperature difference $\Delta T_{3-9}$ is labeled.
  }
  \label{fig:thermosyphons}
\end{figure}

The reduced order system describing a thermal convection loop was originally derived by Gorman \shortcite{gorman1986} and Ehrhard and M\"{u}ller \shortcite{ehrhard1990dynamical}.
Here we present this three dimensional system in non-dimensionalized form.
In Appendix B we present a more complete derivation of these equations, following the derivation of Harris \shortcite{harris2011predicting}.
For $x_{\{1,2,3\}}$ the mean fluid velocity, temperature different $\Delta T_{3-9}$, and deviation from conductive temperature profile, respectively, these equations are:
\begin{align}
& \diff{x_1}{t} = \alpha (x_2 - x_1),\\
& \diff{x_2}{t} = \beta x_1 - x_2 (1 + Kh(|x_1|)) - x_1x_3,\\
& \diff{x_3}{t} = x_1x_2 - x_3 (1 + Kh(|x_1|)) .\end{align}
The parameters $\alpha$ and $\beta$, along with scaling factors for time and each model variable can be fit to data using standard parameter estimation techniques.

Tests of the data assimilation algorithms in Chapter 2 were performed with the Lorenz 63 system, which is analogous to the above equations with Lorenz's $\beta = 1$, and $K = 0$.
From these tests, we have found the optimal data assimilation parameters (inflation factors) for predicting time-series with this system.
We focus our efforts on making prediction using computational fluid dynamics models.
For a more complete attempt to predict thermosyphon flow reversals using this reduced order system, see the work of Harris et al (2012).

\subsection{Experimental Setup}

\begin{figure}[h!]
  \centering
  \includegraphics[width=0.95\textwidth]{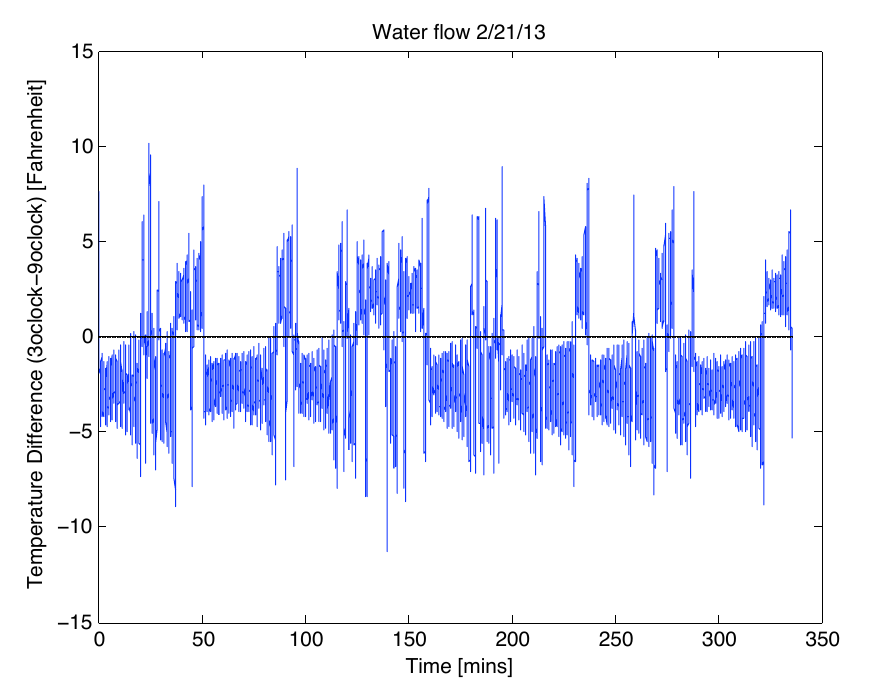}
  \caption[A time series of the physical thermosyphon, from the Undergraduate Honor's Thesis of Darcy Glenn {\protect \shortcite{glenn2013}}]{
    A time series of the physical thermosyphon, from the Undergraduate Honor's Thesis of Darcy Glenn {\protect \shortcite{glenn2013}}.
    The temperature difference (plotted) is taken as the difference between temperature sensors in the 3 and 9 o'clock positions.
    The sign of the temperature difference indicates the flow direction, where positive values are clockwise flow.
  }
  \label{fig:chrisLoop}
\end{figure}

\begin{figure}[h!]
  \centering
  \includegraphics[width=0.79\textwidth]{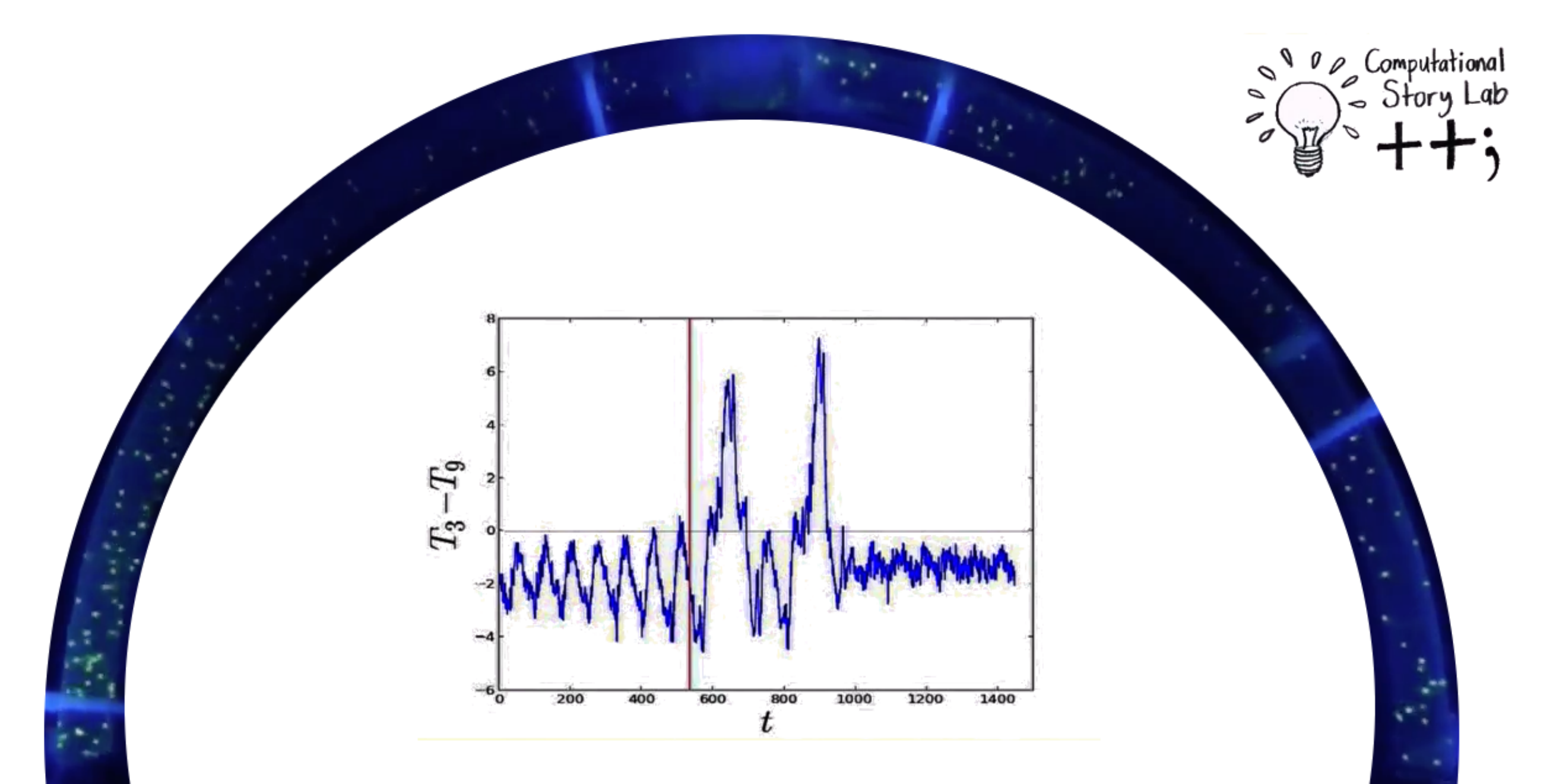}
  \caption[The Thermal Convection Loop experiment]
  {
    The Thermal Convection Loop experiment, inlaid with a time series of the temperature difference $\Delta T_{3-9}$ whose sign indicates the flow direction.
    Fluorescent balls are lit by a black-light and move with the flow.
    As the timeseries indicates, the flow changes direction aperiodically, giving rise to the chaotic behavior that we model.
    A video of the loop, from which this is a snapshot, can be viewed at \url{http://www.youtube.com/watch?feature=player_embedded&v=Vbni-7veJ-c}.
	  }
  \label{fig:chrisLoop}
\end{figure}

Originally built (and subsequently melted) by Dr. Chris Danforth at Bates College, the University of Vermont maintains experimental thermosyphons.
Operated by Dave Hammond, UVM's Scientific Electronics Technician, the experimental thermosyphons access the chaotic regime of state space found in the principled governing equations.
We quote the detailed setup from Darcy Glenn's undergraduate thesis \shortcite{glenn2013}:
\begin{quote}
The [thermosyphon] was a bent semi-flexible plastic tube with a 10-foot heating rope wrapped around the bottom half of the upright circle.
The tubing used was light-transmitting clear THV from McMaster-Carr, with an inner diameter of 7/8 inch, a wall thickness of 1/16 inch, and a maximum operating temperature of 200F.
The outer diameter of the circular thermosyphon was 32.25 inches.
This produced a ratio of about 1:36 inner tubing radius to outside thermosyphon radius (Danforth 2001).
There were 1 inch 'windows' when the heating cable was coiled in a helix pattern around the outside of the tube, so the heating is not exactly uniform.
The bottom half was then insulated using aluminum foil, which allowed fluid in the bottom half to reach 176F.
A forcing of 57 V, or 105 Watts, was required for the heating cable so that chaotic motion was observed.
Temperature was measured at the 3 o'clock and 9 o'clock positions using unsheathed copper thermocouples from Omega.
\end{quote}

We first test our ability to predict this experimental thermosyphon using synthetic data, and describe the data in more detail in Chapter 4.
Synthetic data tests are the first step, since the data being predicted is generated by the predictive model.
We consider the potential of up to 32 synthetic temperature sensors in the loop, and velocity reconstruction from particle tracking, to gain insight into which observations will make prediction possible.

\section{Computational Experiment}
In this section, we first present the governing equations for the flow in our thermal convection loop experiment.
A spatial and temporal discretization of the governing equations is then necessary so that they may be solved numerically.
After discretization, we must specify the boundary conditions.
With the mesh and boundary conditions in place, we can then simulate the flow with a computational fluid dynamics solver.

We now discuss the equations, mesh, boundary conditions, and solver in more detail.
With these considerations, we present our simulations of the thermosyphon.

\subsection{Governing Equations}
We consider the incompressible Navier-Stokes equations with the Boussinesq approximation to model the flow of water inside a thermal convection loop.
Here we present the main equations that are solved numerically, noting the assumptions that were necessary in their derivation.
In Appendix C a more complete derivation of the equations governing the fluid flow, and the numerical for our problem, is presented.
In standard notation, for $u,v,w$ the velocity in the $x,y,z$ direction, respectively, the continuity equation for an incompressible fluid is
\begin{equation} \frac{\partial u}{\partial x} + \frac{\partial v}{\partial y} + \frac{\partial w}{\partial z} = 0. \label{eq:NScontIco} \end{equation}

The momentum equations, presented compactly in tensor notation with bars representing averaged quantities, are given by
\begin{equation} \rho_\text{ref} \left ( \frac{\partial \bar{u}_i}{\partial t} + \frac{\partial}{\partial x_j} \left( \bar{u}_j \bar{u}_i \right) \right )
= -\frac{\partial \bar{p}} {\partial{x_i}} +  \mu \frac{\partial \bar{u}_i}{\partial x_j^2} + \bar{\rho} g_i \end{equation}
for $\rho_\text{ref}$ the reference density, $\rho$ the density from the Boussinesq approximation, $p$ the pressure, $\mu$ the viscocity, and $g_i$ gravity in the $i$-direction.
Note that $g_i = 0$ for $i \in \{ x,y\}$ since gravity is assumed to be the $z$ direction.
The energy equation is given by
\begin{equation} \frac{\partial }{\partial t} \left ( \rho \overline{e}\right ) + \frac{\partial}{\partial x_j} \left ( \rho \overline{e} \overline{u}_j \right )
=
- \frac{\partial q_k^*}{\partial x_k}
- \frac{\partial \overline{q}_k}{\partial x_k}
\end{equation}
for $e$ the total internal energy and $q$ the flux (where $q = \overline{q} + q^*$ is the averaging notation).

\begin{figure}[t!]
  \centering
  \includegraphics[width=0.79\textwidth]{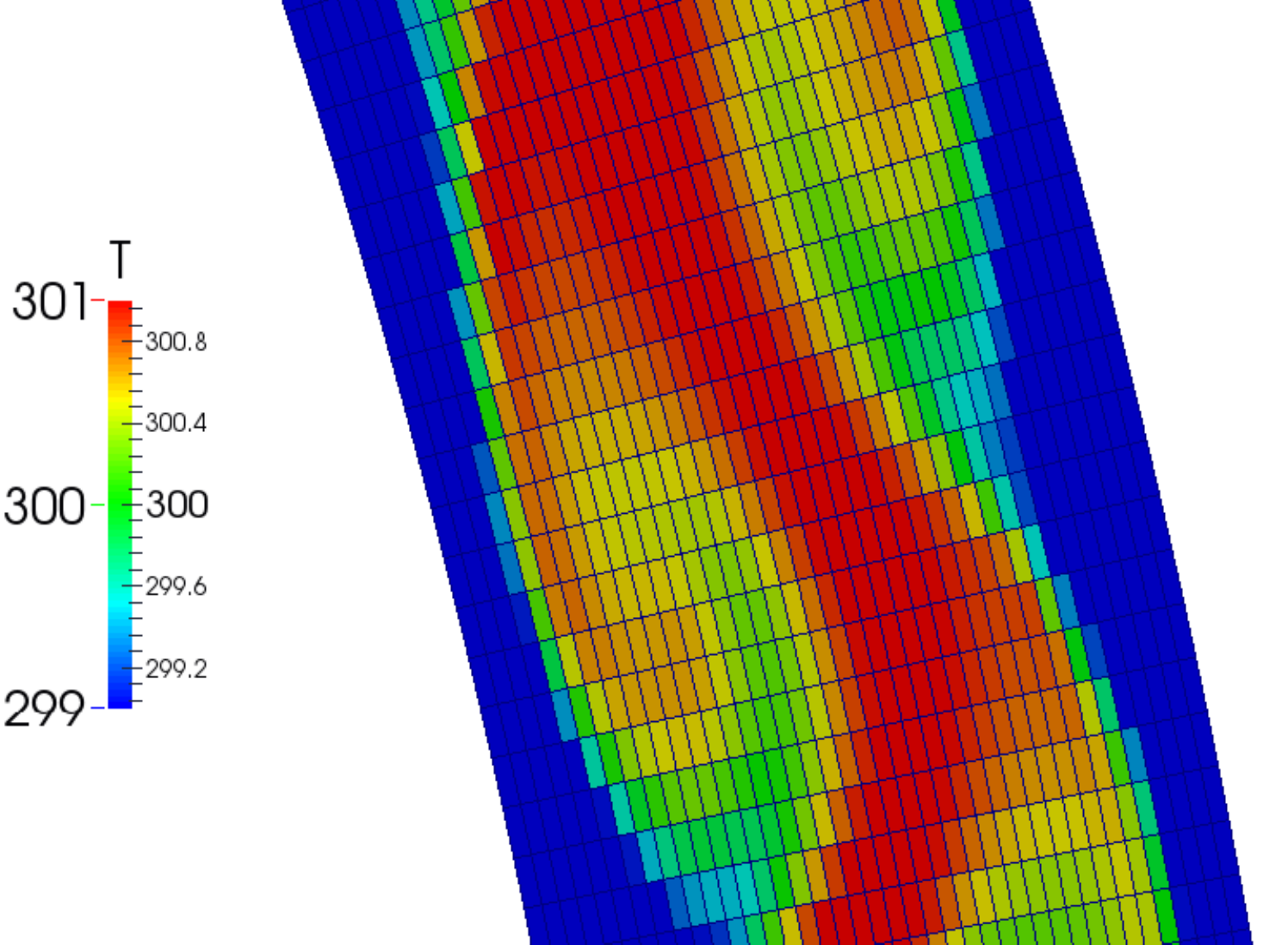}
  \caption[A snapshot of the mesh used for CFD simulations]{
  A snapshot of the mesh used for CFD simulations.
  Shown is an initial stage of heating for a fixed value boundary condition, 2D, laminar simulation with a mesh of 40000 cells including wall refinement with walls heated at 340K on the bottom half and cooled to 290K on the top half.
  }
  \label{fig:CFDmesh1}
\end{figure}

\subsection{Mesh}

Creating a suitable mesh is necessary to make an accurate simulation.
Both 2-dimensional and 3-dimensional meshes were created using OpenFOAM's native meshing utility ``blockMesh.''
After creating a mesh, we consider refining the mesh near the walls to capture boundary layer phenomena and renumbering the mesh for solving speed.
To refine the mesh near walls, we used the ``refineWallMesh'' utility,
Lastly, the mesh is renumbered using the ``renumberMesh'' utility, which implements the Cuthill-McKee algorithm to minimize the adjacency matrix bandwidth, defined as the maximum distance from diagonal of nonzero entry.
The 2D mesh contains 40,000 points, and is imposed with fixed value boundary conditions on the top and bottom of the loop.

\begin{figure}[t!]
  \centering
  \includegraphics[width=0.79\textwidth]{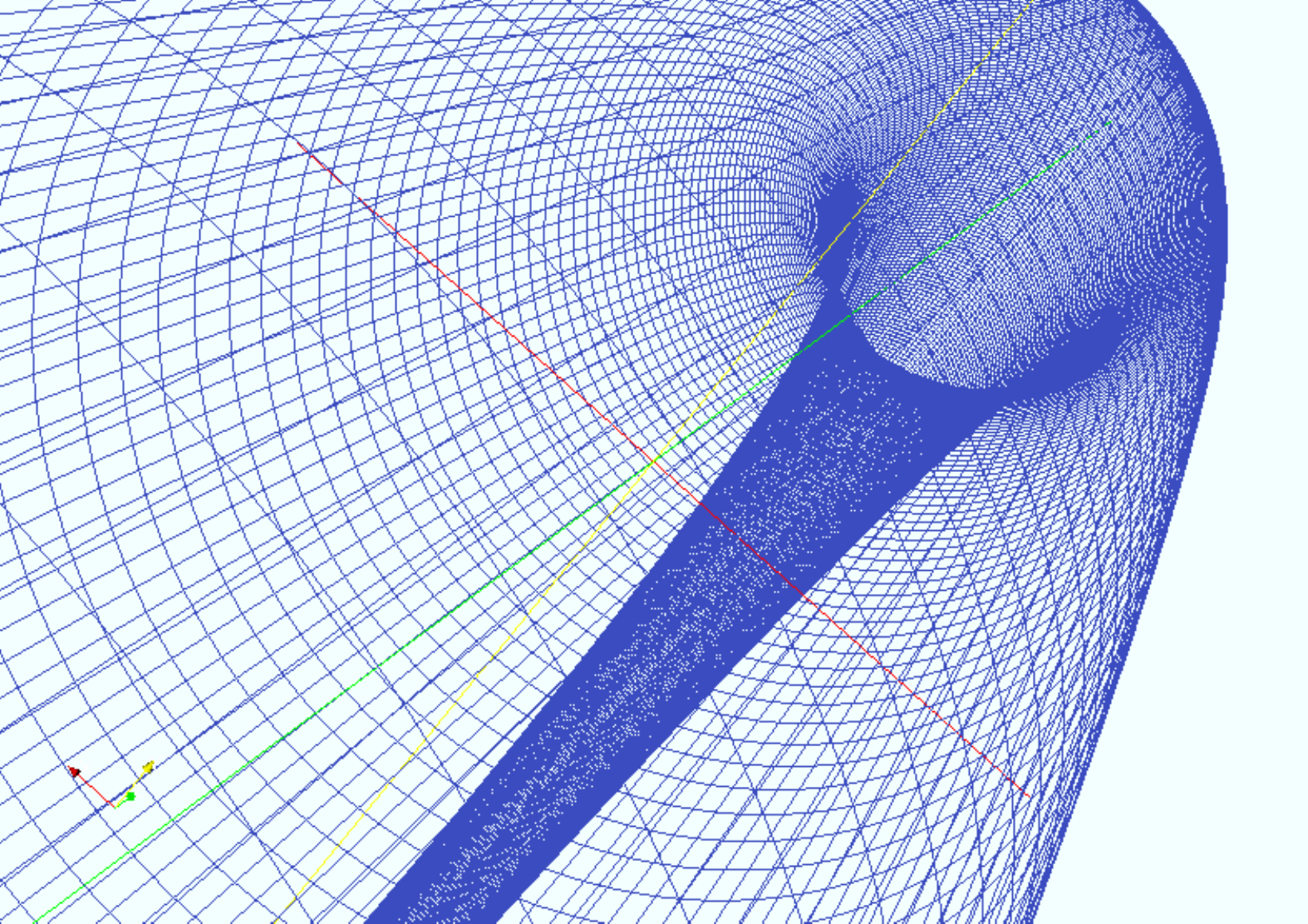}
  \caption[The 3D mesh viewed as a wireframe from within]{
  The 3D mesh viewed as a wireframe from within.
  Here there are 900 cells in each slice (not shown), for a total mesh size of 81,000 cells.
  Simulations using this computational mesh are prohibitively expensive for use in a real time ensemble forecasting system, but are possible offline.
  }
  \label{fig:CFDmesh2}
\end{figure}

\begin{figure}[t!]
  \centering
  \includegraphics[width=0.79\textwidth]{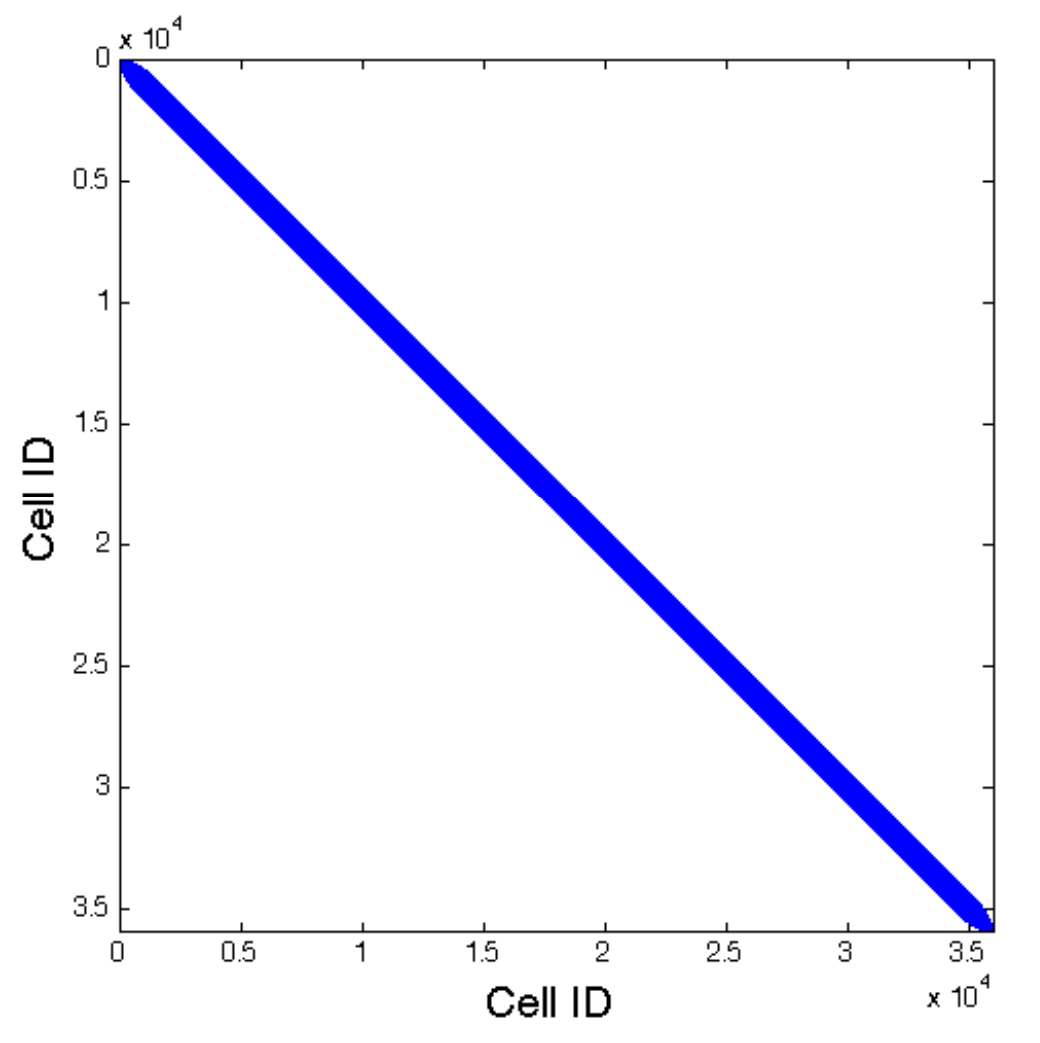}
  \caption[The 2D mesh Adjacency matrix for a local radius of 1.0cm]{
  The 2D mesh Adjacency matrix for a local radius of 1.0cm.
  In general, we see that cells have neighbors with nearby cell ID.
  With 5,260,080 nonzero entries, it is difficult to see structure looking at the entire matrix.
  This figure was bitmapped for inclusion on the arXiv, for a full version see {\protect \url{http://andyreagan.com}}.
  }
  \label{fig:CFDmeshAdj1}
\end{figure}

\begin{figure}[h!]
  \centering
  \includegraphics[width=0.79\textwidth]{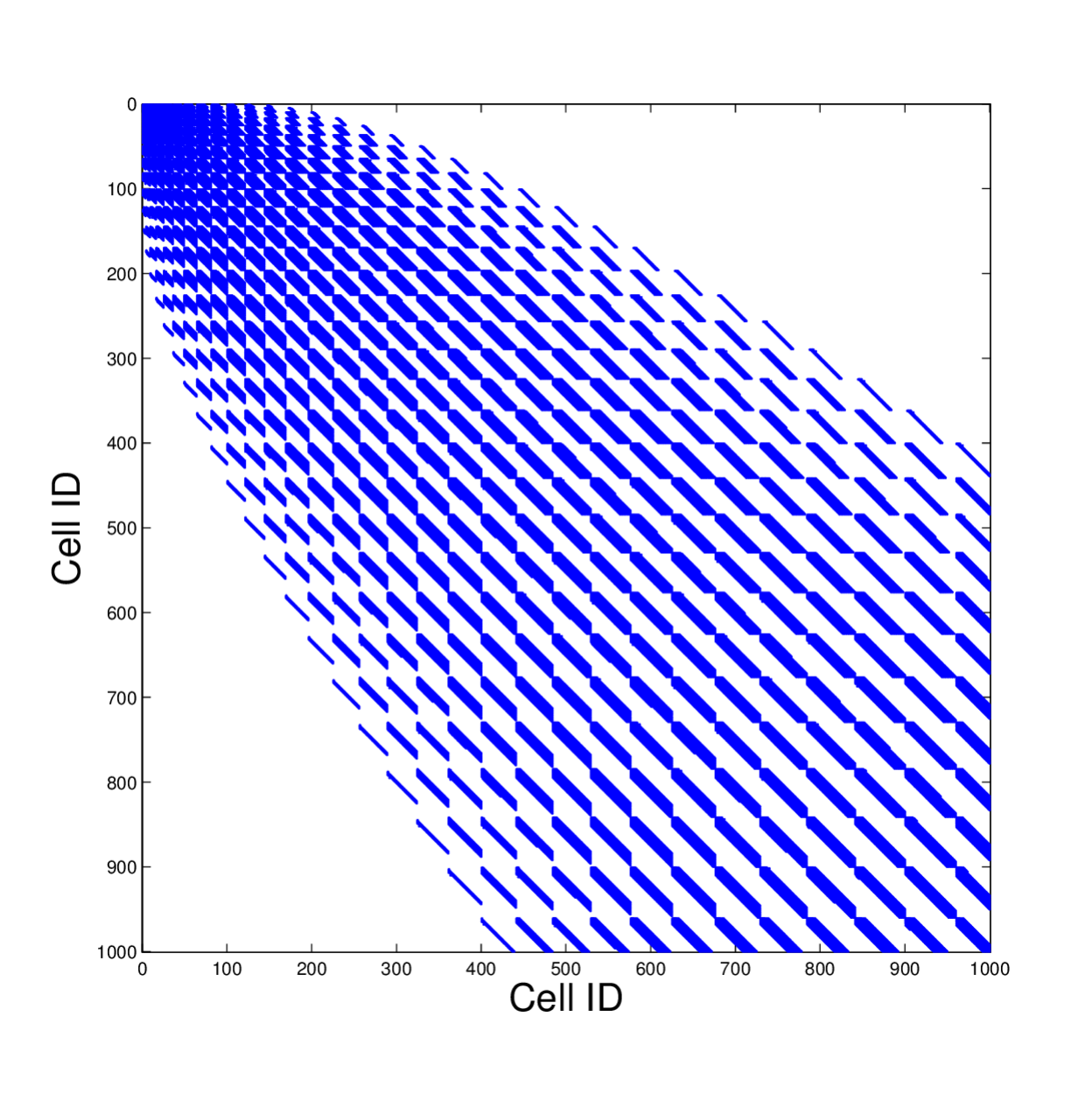}
  \caption[The 2D mesh adjacency matrix, zoomed in]{
  The 2D mesh adjacency matrix, zoomed in.
  The Cuthill-McKee renumbering algorithm was used, and we see 23 distinct regions of local numbering, reflective of the geometry and Cuthill-McKee algorithm.
  The mesh shown here has 30 cells in each perpendicular slice of the loop, and they are numbered iteratively along the edge of triangular front emanating from cell 0 (chose to be the left-most cell).
  This diagonalization in both space and ID number effectively reduces the bandwidth of the mesh from 947 to 18.
  This figure was bitmapped for inclusion on the arXiv, for a full version see {\protect \url{http://andyreagan.com}}.
    }
  \label{fig:CFDmeshAdj2}
\end{figure}

\subsection{Boundary Conditions}

Available boundary conditions (BCs) that were found to be stable in OpenFOAM's solver were constant gradient, fixed value conditions, and turbulent heat flux.
Simulations with a fixed flux BC is implemented through the externalWallHeatFluxTemperature library were unstable and resulted in physically unrealistic results.
Constant gradient simulations were stable, but the behavior was empirically different from our physical system.
Results from attempts at implementing these various boundary conditions are presented in Table \ref{tb:ofbc}.
To most accurately model our experiment, we choose a fixed value BC.
This is due to the thermal diffusivity and thickness of the walls of the experimental setup.
For $Q$ the flow of heat from the hot water bath into the loop's outer material, and $A$ the area of a small section of boundary we can write the heat flux $q$ into the fluid inside the loop as a function of the thermal conductivity $k$, temperature difference $\Delta T$ and wall thickness $\Delta x$:
\begin{equation} q = \frac{Q}{A} = k \frac{\Delta T}{\Delta x} . \end{equation}
For our experimental setup, the thermal conductivity of a PVC wall is approximately $1.1 \, \left [ \text{BTU in}/(h\cdot \text{ft}^2 F) \right ]  = 0.15 \, \left [W/(m\cdot K)\right ] $.
With a wall thickness of 1cm, and $\Delta T = 10 \left [K \right ]$ the heat flux is $q=79 \left [ W/m^2 \right ]$.
The numerical implementation of this boundary condition is considered further in Appendix C.
In Table \ref{tb:ofbc} we list the boundaries conditions with which we experimented, and the quality of the resulting simulation.

\begin{table}
\begin{tabular}{llll}
\hline
Boundary Condition & Type & Values Considered & Result\\
\hline
\hline
fixedValue & value & $\pm 5 \to \pm 30$ & realistic\\
\hline
fixedGradient & gradient & $\pm 10 \to \pm 100000$ & plausible\\
\hline
wallHeatTransfer & enthalpy & - & cannot compile\\
\hline
externalWallHeatFluxTemperature & flux & - & cannot compile\\
\hline
turbulentHeatFluxTemperature & flux & $\pm .001 \to \pm10$ & unrealistic\\
\hline
\end{tabular}
\label{tb:ofbc}
\caption[Boundary conditions considered for the OpenFOAM CFD model are labeled with their type, values considered, and empirical result]{
  Boundary conditions considered for the OpenFOAM CFD model are labeled with their type, values considered, and empirical result.
  We find the most physically realistic simulations using the simplest ``fixedValue'' condition.
  The most complex conditions, ``wallHeatTransfer,'' ``externalWallHeatFluxTemperature,'' and ``turbulentHeatFluxTemperature'' either did not compile or produced divergent results.
  }
\end{table}

\subsection{Solver}

In general, there are three approaches to the numerical solution of PDEs: finite differences, finite elements, and finite volumes.
We perform our computational fluid dynamics simulation of the thermosyphon using the open-source software OpenFOAM.
OpenFOAM has the capacity for each scheme, and for our problem with an unstructured mesh, the finite volume method is the simplest.
The solver in OpenFOAM that we use, with some modification, is ``buoyantBoussinesqPimpleFoam.''
Solving is accomplished by the Pressure-Implicit Split Operator (PISO) algorithm \shortcite{issa1986solution}.
Modification of the code was necessary for laminar operation.
In Table \ref{tb:ofsolution}, we note the choices of schemes for time-stepping, numerical methods and interpolation.

The Boussinesq approximation makes the assumption that density changes due to temperature differences are sufficiently small so that they can be neglected without effecting the flow, except when multiplied by gravity.
The non-dimensionalized density difference $\Delta \rho / \rho$ depends on the value of heat flux into the loop, but remains small even for large values of flux.
Below, we approximate the validity of the Boussinesq approximation, which is generally valid for values of $\beta(T-T_\text{ref})\rhoref \ll 1$:
\begin{equation} \frac{\beta(T - T_\text{ref})}{\rhoref} \simeq 7.6 \cdot 10^{-6} \ll 1 \label{eq:bouss} \end{equation}

\begin{table*}
\begin{center}
\begin{tabular}{ll}
\hline
Problem & Scheme\\
\hline
\hline
Main solver & buoyantBoussinesqPimpleFoam\\
\hline
Solving algorithm & PISO\\
\hline
Grid & Semi-staggered\\
\hline
Turbulence model & None (laminar)\\
\hline
Time stepping & Euler\\
\hline
Gradient & Gauss linear\\
\hline
Divergence & Gauss upwind\\
\hline
Laplacian & Gauss linear\\
\hline
Interpolation & Linear\\
\hline
Matrix preconditioner & Diagonal LU\\
\hline
Matrix solver & Preconditioned Bi-linear Conjugate Gradient\\
\hline
\end{tabular}
\caption[Numerical schemes used for solving the finite volume problem]{
  Numerical schemes used for solving the finite volume problem.
  }
\label{tb:ofsolution}
\end{center}
\end{table*}

\subsection{Simulation results}
With the mesh, BCs, and solver chosen, we now simulate the flow.
From the data of $T,\phi,u,v,w$ and $p$ that are saved at each timestep, we extract the mass flow rate and average temperature at the $12,3,6$ and $9$ o'clock positions on the loop.
Since $\phi$ is saved as a face-value flux, we compute the mass flow rate over the cells $i$ of top (12 o'clock) slice as
\begin{equation} \sum _i\phi_{f(i)} \cdot v_i \cdot \rho_i\end{equation}
where $f(i)$ corresponds the face perpendicular to the loop angle at cell $i$ and $\rho$ is reconstructed from the Boussinesq approximation $\rho = \rhoref (1-\beta(T-T_\text{ref}))$.

Here is a picture of the thermosyphon colored by temperature.
\begin{figure}[h!]
  \centering
  \includegraphics[width=0.95\textwidth]{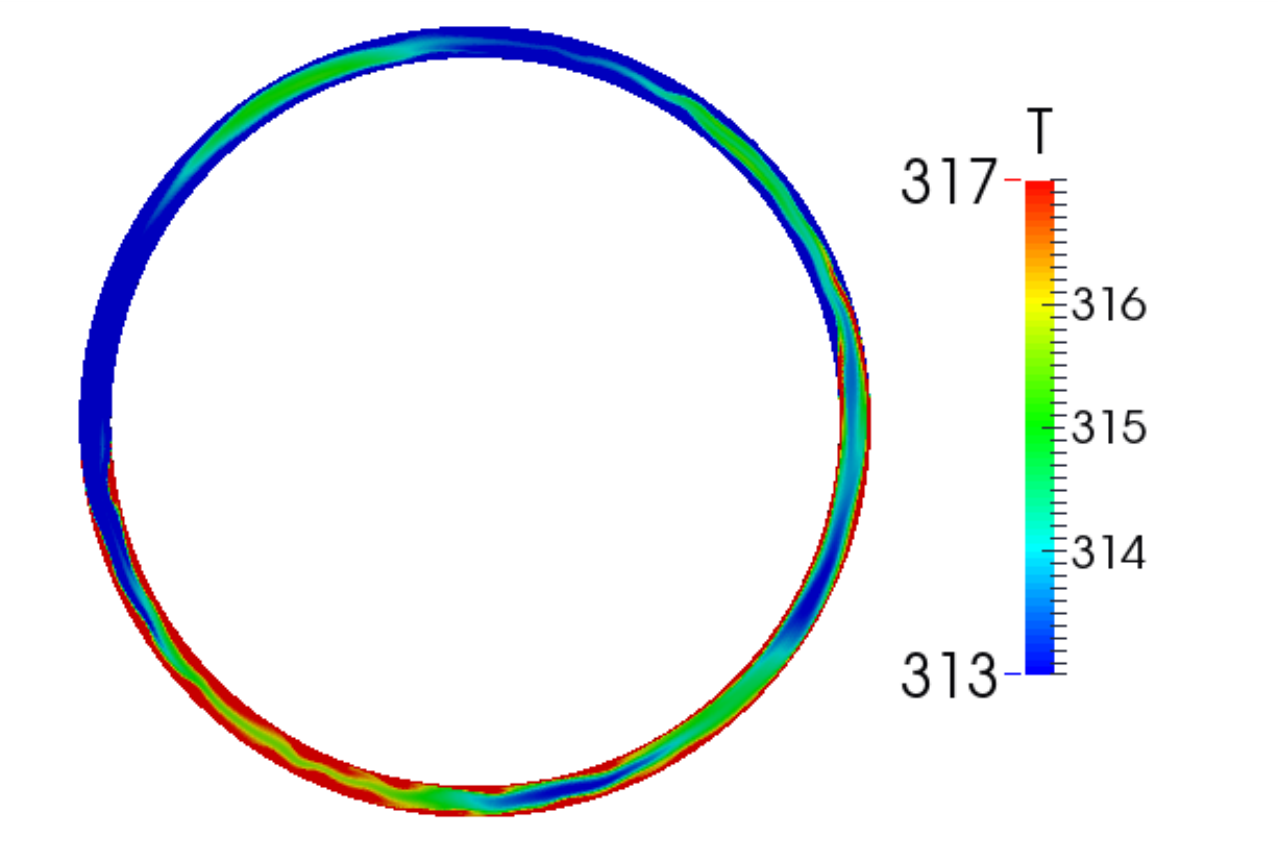}
  \caption[A screenshot of the whole simulated loop]{
  A screenshot of the whole computationally simulated loop.
  Shown is an initial stage of heating for a fixed value boundary condition, 2D, laminar simulation with a mesh of 40000 cells including wall refinement with walls heated at 340K on the bottom half and cooled to 290K on the top half.
  Note that the colorbar limits are truncated to allow for visualization of the dynamics far from the boundary.
  }
  \label{fig:CFDloopSS}
\end{figure}

\chapter{Results}
\chaptermark{}

\section{Computational Model and DA Framework}

\subsection{Core framework}

The first output of this work is a general data assimilation framework for MATLAB.
By utilizing an object-oriented (OO) design, the model and data assimilation algorithm code are separate and can be changed independently.
The principal advantage of this approach is the ease of incorporation of new models and DA techniques.

To add a new model to this framework, it is only necessary to write a MATLAB class wrapper for the model, and the existing Lorenz 63 and OpenFOAM models are templates.
Testing a new technique for performing the DA step is as simple as writing a MATLAB function for the assimilation steps, with the model forecast and observations as input.

The core of this framework is the MATLAB script modelDAinterface.m which wraps observations, the model class, and any DA method together.
A description of the main algorithm is provided in Appendix \ref{app:foamlab}.

For large models, it is necessary to perform the assimilation in a local space, and this consideration is accounted for within the framework.
This procedure, referred to as localization, eliminates spurious long-distance correlations and parallelizes the assimilation.

\begin{figure}[t!]
  \centering
  \includegraphics[width=0.79\textwidth]{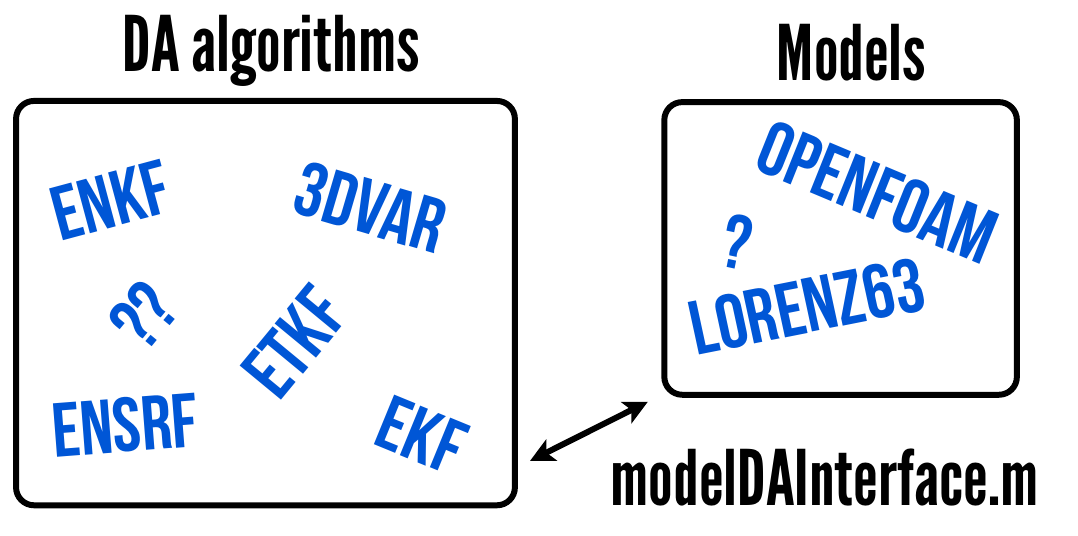}
  \caption[A schematic of the modelDAinterface.m functionality]{
    A schematic of the modelDAinterface.m functionality.
    By linking arbitrary data assimilation algorithms and models we have the ability to test prediction of any given timeseries with a choice of model, algorithm, and localization scheme.
    This functionality was necessary to perform consistent tests for many different model parameterizations and observation choices.
  }
  \label{fig:modelDAschem}
\end{figure}

\subsection{Tuned OpenFOAM model}

In addition to wrapping the OpenFOAM model into MATLAB, defining a physically realistic simulation was necessary for accurate results.
A non-uniform sampling of the OpenFOAM model parameter space across BC type, BC values, and turbulence model was performed.
The amount of output data constitutes what is today ``big data,'' and for this reason it was difficult to find a suitable parameter regime systematically.
Specifically, each model run generates 5MB of data per time save, and with 10000 timesteps this amounts to roughly 50GB of data for each model run.
We reduce the data to time-series for a subset of model variables, namely timeseries of the reconstructed mass flow rate and temperature difference $\Delta T_{3-9}$.
From these timeseries, we then select the BC, BC values, and turbulence model empirically based on the flow behavior observed.

\begin{figure}[h!]
  \centering
    \includegraphics[width=0.79\textwidth]{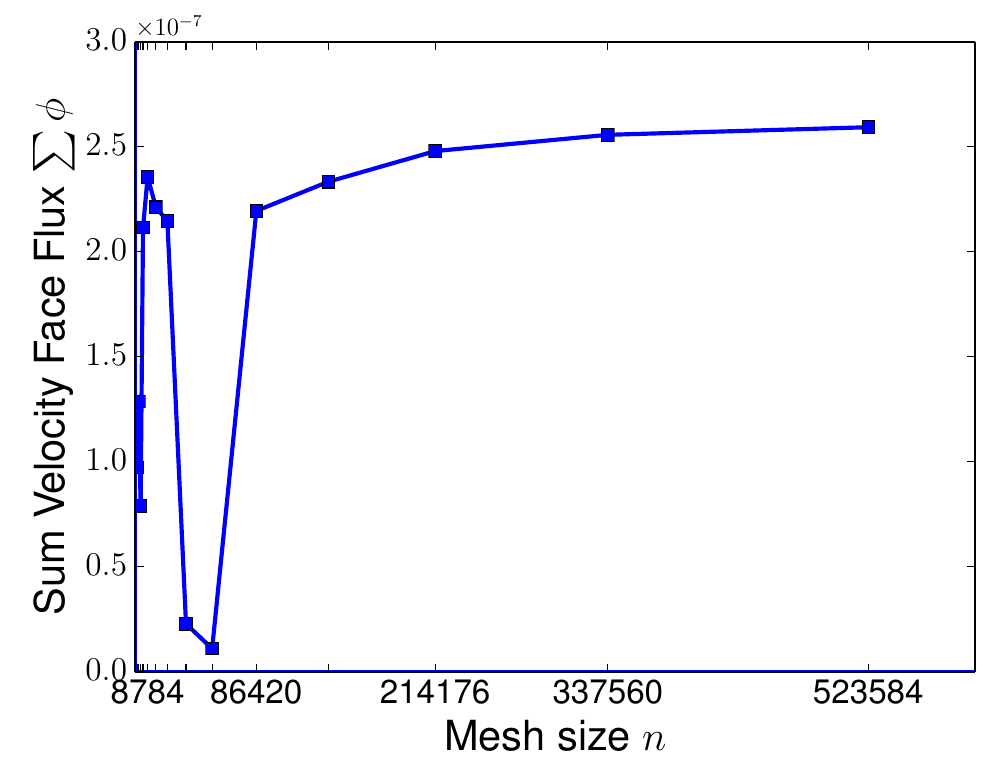}
  \caption[Long-term behavior for different meshes]
  {
    With a fixed choice of solver, boundary conditions, and initial conditions that lead to a stable convective state, we present the long-term behavior of the velocity face flux at the top slice for different meshes.
    The face flux is reported as the average for the last 20 times saves for which the velocity flux is summed across a slice perpendicular to the loop, here we show the top slice.
    We choose a fixed time step of 0.005 for each simulation, and run the solver for 60 hours on 8 cores.
    The computational limit of mesh creation was a memory limit at 818280 cells, so we present results for meshes starting at 1600 cells and cells decreasing in size by a factor of 1.25 in both $y$ and $z$ up to a mesh of 523584 cells.
    For meshes with more than 80,000 cells we see that the solutions are very similar.
    The smaller meshes generate increasing unstable flow behavior, leading to oscillations of flux and then flow reversals for the smallest meshes of size 2500 and 1600 cells.
  }
  \label{fig:meshverification}
\end{figure}

\section{Prediction Skill for Synthetic Data}

We first present the results pertaining to the accuracy of forecasts for synthetic data.
There are many possible experiments given the choice of assimilation window, data assimilation algorithm, localization scheme, model resolution, observational density, observed variables, and observation quality.
We focus on considering the effect of observations and observational locations on the resulting forecast skill.

\subsubsection{Covariance Localization}

\begin{figure}[h!]
  \centering
  \includegraphics[width=0.79\textwidth]{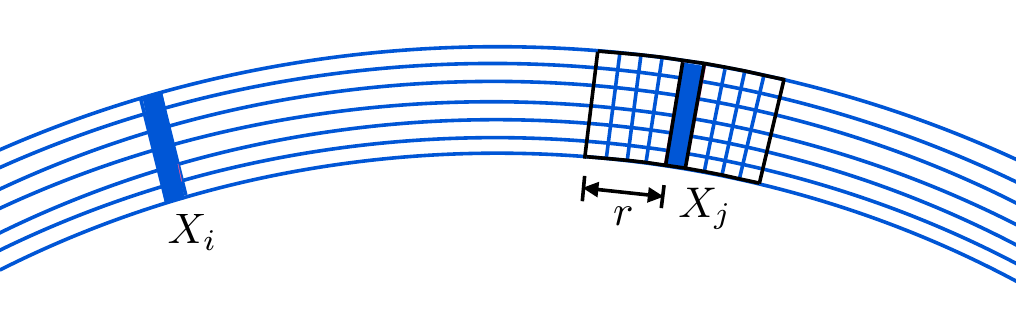}
  \caption[A schematic of covariance localized by slice is shown for example slices $X_{\{i,j\}}$ and localization radius $r$ where $i,j$ are indexed within the number of slices]{
    A schematic of covariance localized by slice is shown for example slices $X_{\{i,j\}}$ and localization radius $r$ where $i,j$ are indexed within the number of slices.
    Note that for this scheme, $r$ is chosen as an integer.
    Of the localization schemes tested, slice localization requires the least individual assimilations by considering groups of the size of cells in each slice.
  }
  \label{fig:localcovzones}
\end{figure}

\begin{figure}[h!]
  \centering
  \includegraphics[width=0.79\textwidth]{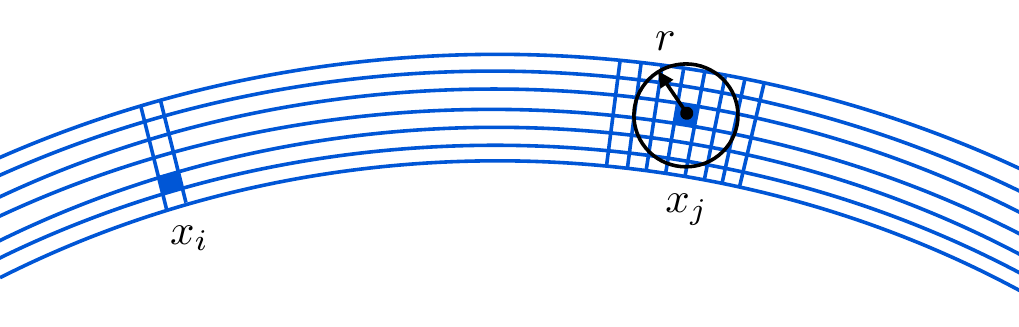}
  \caption[A schematic of covariance localized by a radial distance from cell centers is shown for cells $x_{\{i,j\}}$ where $i,j$ are indexed by the number of cells and $r$ is any float]{
    A schematic of covariance localized by a radial distance from cell centers is shown for cells $x_{\{i,j\}}$ where $i,j$ are indexed by the number of cells and $r$ is any float.
    Cells are local to one another if the Euclidean distance between their cell centers is less than $r$.
    Assimilation is performed for each cell, so there are as many individual assimilations are cells in the mesh.
    Since observations are assumed to be independent, this assimilation can be performed in parallel for efficiency.
  }
  \label{fig:localcovradii}
\end{figure}

The need for covariance localization was discussed in Chapter 2, and here we implement two different schemes.
First, we consider cells within each perpendicular slice of the loop to be a cohesive group, and their neighbors those cells within $r$ slices (where here, $r$ is an integer).
Computationally this scheme is the simplest to implement, and results in the fewest covariance localizations.

Next, we consider localization by distance from center.
Data assimilation is performed for each grid point, considering neighboring cells and observations with cell center within a radius $r$ from each point (e.g. $r = 0.5$cm).
The difficulty of this approach is the localization itself which is mesh dependent and nontrivial, see Figure \ref{fig:CFDmeshAdj2}.
For our simulations, we perform the assimilations sequentially, and note that there many be gains in efficiency from parallelization.

\subsubsection{Generation of Synthetic Data}

Synthetic data is generated from one model run, with a modified ``buoyantBoussinesqPimpleFoam'' solver applied to fixed value boundary conditions of 290K top and 340K bottom, 80,000 grid points, a time step of .01, and laminar turbulence model.
Data is saved every half second, and is subsampled for 2,4,8,16 and 32 synthetic temperature sensors spaced evenly around the loop.
The temperature at these synthetic sensors is reported as the average for cells with centers within 0.5cm of the observation center.
Going beyond 32 temperature observations, which we find to be necessary, samples are performed by observed cell-center temperature values.
The position of the chosen cell depends on the number of observations being made, and within that number there are options that are tested.
For less than 1000 observations, (1000 being the number of angular slices) observations are placed at slice centers, for slices spaced maximally apart.
For 1000 observations both slice center (midpoint going all the way around) and staggered (closer to one wall, alternating) are implemented.
For $n$ observations with $n>1000$, assuming that $n/1000$ is an integer, for each slice we space the $n/1000$ observations maximally far apart within each slice.

In addition to temperature sensing, the velocity in $y$ and $z$ is sampled at 50, 100, and 300 points chosen randomly at each time step.
These velocity observations are designed to simulate a reconstructed velocity measurement based on video particle tracking.
The following results do not include the use of these velocity observations.
To obtain realistic sampling intervals, this data was aggregated and is then subsampled at the assimilation window length.

For both the Lorenz 63 model and OpenFOAM, error from the observed $T_{3-9}$ is computed as RMSE, such that the models can be compared.
In the EM model, a maximum of 2 temperature sensors are considered for observations at $T_3$ and $T_9$, and no mean velocity reconstruction from subsampled velocity is attempted.

For any of the predictions to work (even with complete observation coverage), it was first necessary to scale the model variables to the same order of magnitude.
This is accomplished here by dividing through by the mean climatological values from the observations and model results before assimilation (e.g. multiplying temperature by $1/300$).
While other, more complex, strategies may result in better performance there are no general rules to determine the best scaling factors and good scaling is often problem dependent \shortcite{navon1992variational,Mitchell2012PhD}.

\subsection{Observational Network}

In general, we see that increasing observational density leads to improved forecast accuracy.
With too few observations, the data assimilation is unable to recover the underlying dynamics.

\begin{figure}[h!]
  \centering
  \includegraphics[width=0.79\textwidth]{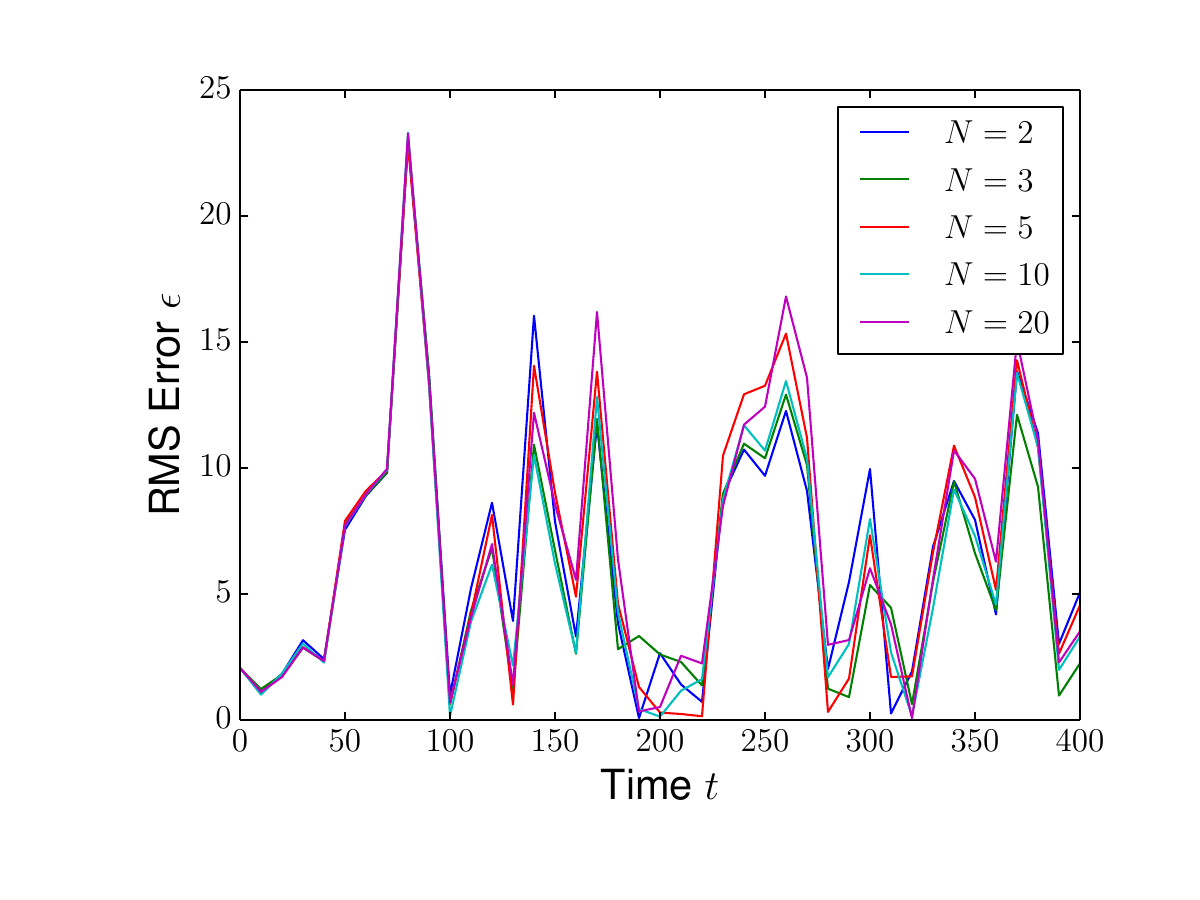}
  \caption[Prediction skill over time for synthetic temperature sensors]{
    Prediction skill over time for synthetic temperature sensors.
    With these observations, our DA scheme is unable to predict the model run at all.
    Compared with operational NWP, we are observing far less model variables in this situation, with $10^5$ times more model variables than observations, even for the 32 sensors.
  }
  \label{fig:TsensorExp}
\end{figure}

\begin{figure}[h!]
  \centering
  \includegraphics[width=0.79\textwidth]{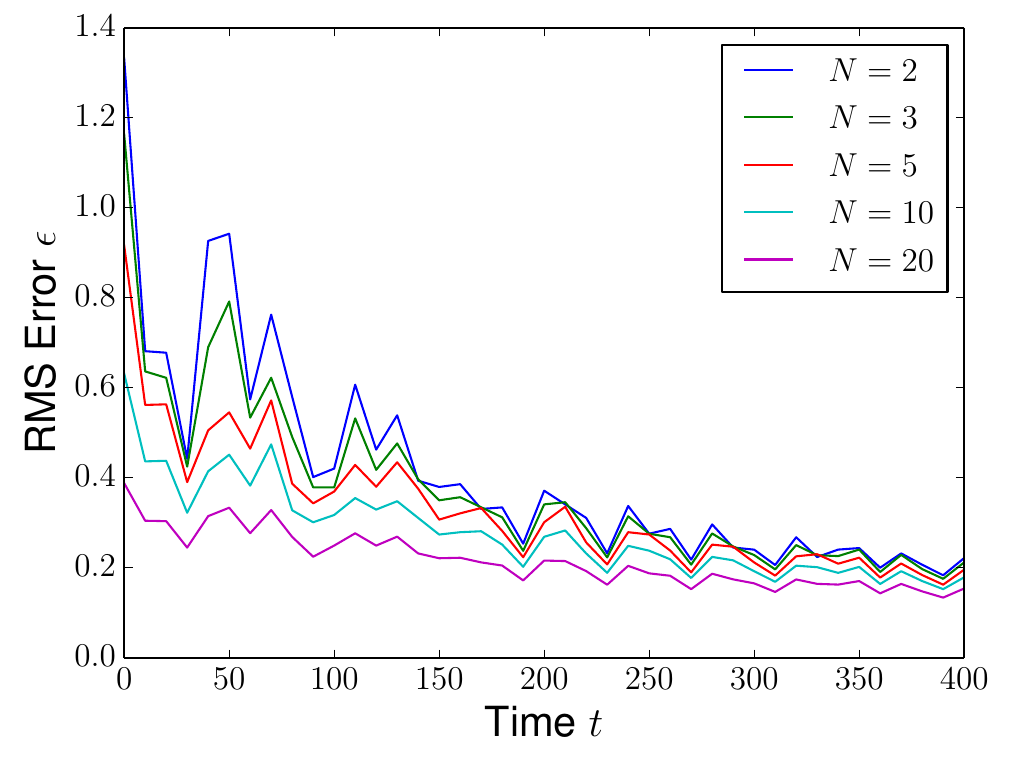}
  \caption[Prediction skill over time for complete observations of temperature]{
    Prediction skill over time for complete observations of temperature.
    With these observations, our DA scheme is able to predict the model run well.
    Increasing the number of ensemble members $N$ leads to a decrease in the forecast error, as one would hope with a more complete sampling of model uncertainty.
    Note that although we may full observations of temperature, it is only one of five model variables assimilated.
  }
  \label{fig:fullObsExp}
\end{figure}

While the inability to predict the flow direction with these synthetic temperature sensors may cast doubt on the possibility of predicting a physical experiment using CFD, we discuss improvements to our prediction system that may be able to overcome these initial difficulties.
Regardless, in Figure \ref{fig:fullObsExp} we verify directly that our assimilation procedure is working adequately with a sufficient number of observations: temperature at all cells.

With success at predicting flow direction given full temperature observations, we now ask the question: just how many observations are necessary?
The aforementioned generation of observational data considers the different possible observational network approaches for given number of observations $N_\text{obs}$.
From these approaches, we select the best method available to assimilation $N_\text{obs}$ for each $N$.
Our results initially indicate that the averaging of temperature over many cells for the synthetic temperature sensors up to 32 sensors degraded prediction skill.

\begin{figure}[h!]
  \centering
  \includegraphics[width=0.79\textwidth]{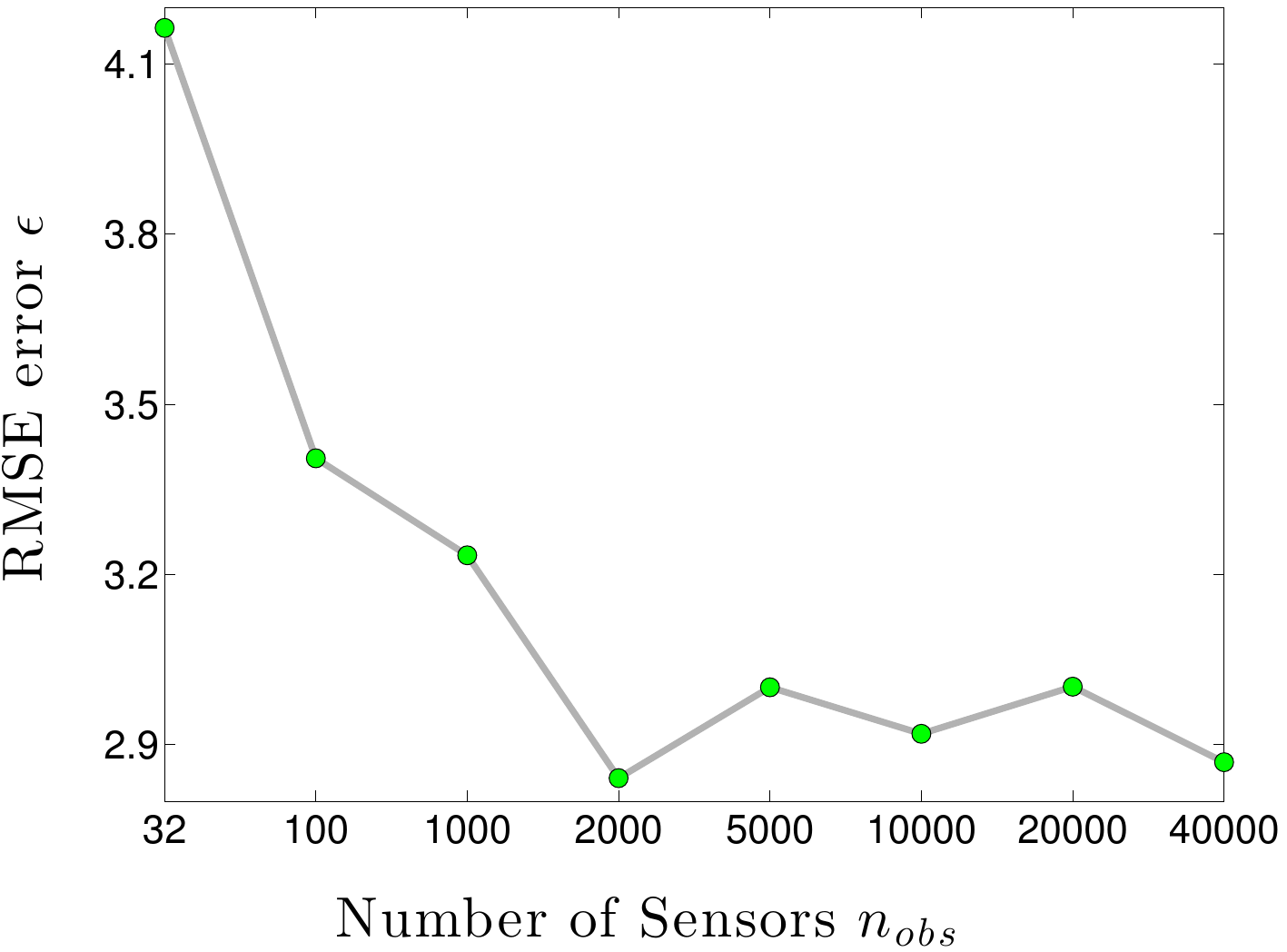}
  \caption[The RMSE of forecasts at the end of 40 assimilation cycles for differing observational density is plotted]
  {
    The RMSE of forecasts at the end of 40 assimilation cycles for differing observational density is plotted.
    For each, we use an assimilation window length of 10, observations with errors drawn from a normal distribution with standard deviation 0.05, and 20 ensemble members within an ETKF.
    No improvement in skill is witnessed when using more than 2000 observations (5\% of the model resolution).
  }
  \label{fig:obsexp}
\end{figure}

\section{Conclusion}

Throughout the process of building a data assimilation scheme from scratch, as well as deriving and tuning a computational model, we have learned a great deal about the difficulties and successes of modeling.
We have found that data assimilation algorithms and CFD solvers are sensitive to parameter choices, and cannot be used successfully without an approach that is specific to the problem at hand.
Such an approach has been developed, and proved to be successful.

For the future of data assimilation, the numerical schemes will need to capitalize on new computing resources.
The trend toward massively parallel computation puts the spotlight on effective localization of the assimilation step, which may be difficult in schemes that do not use ensembles, such as the 4D-Var.
Maintaining TLM and adjoint code is difficult and time-consuming, and with future increases in model resolution and atmospheric understanding, TLMs of full GCM models may become impossible.
Again, this points to the potential use of local ensemble based filters such as the LETKF as implemented here.

\subsection{Future directions}

The immediate next step for this work is to test prediction skill in real time.
All of the pieces are there: an observation-based data assimilation framework and models that run at near real-time speed.
Currently on the VACC, the OpenFOAM simulation with grid size 36K runs with approximately half of real-time speed on one core.
This allows the possibility of running $N$ ensemble methods using $4N$ cores, allowing time for the assimilation.

\begin{figure}[h!]
  \centering
  \includegraphics[width=0.79\textwidth]{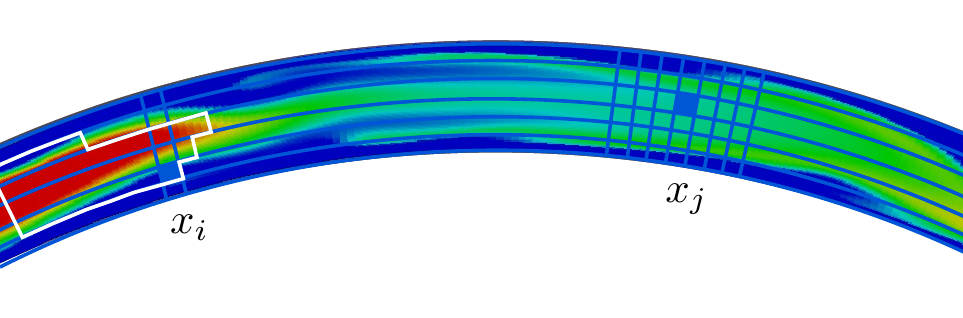}
  \caption[A snapshot of the thermosyphon flow with a sample grid superimposed]{
    A snapshot of the thermosyphon flow with a sample grid superimposed.
    The cells $x_i$ and $x_j$ are labeled and we consider the localization of covariance around cell $x_i$.
    This snapshot was chosen to demonstrate a situation where the local radius and slice zone approaches would include spurious correlations.
    By choosing the localization appropriately, the data assimilation can be improved.
  }
  \label{fig:localcovadaptive}
\end{figure}

Beyond real-time prediction, it is now possible to design and test new methods of data assimilation for large scale systems.
New methods specific to CFD, possibly using principle mode decomposition of flow fields have the potential to greatly improve the skill of real-time CFD applications.
In Figure \ref{fig:localcovadaptive}, flow structures present in our computational experiment demonstrate the potential for adaptive covariance localization.

The numerical coupling of CFD to experiment by DA should be generally useful to improve the skill of CFD predictions in even data-poor experiments, which can provide better knowledge of unobservable quantities of interest in fluid flow.

\bibliographystyle{chicago}

\appendix
\addappheadtotoc

\titleformat{\chapter}[hang]{\normalfont\huge\bfseries}{\chaptertitlename\ \thechapter:}{1em}{}

\chapter{Parameters}
\label{params}

In Appendix A we present the parameters used for the OpenFOAM and Lorenz63 models, as well as the Data Assimilation covariance inflation parameters.
The inflation parameters shown here are those which minimize the RMS forecast error of predictions, using 100 runs of 1000 assimilation windows, for each window length on each filter.

\begin{table}[h]
\caption[Specific constants used for thermal properties of water]{
  Specific constants used for thermal properties of water.
  For temperature dependent quantities, the value is evaluated at the reference temperature $T_\text{ref}  = 300K$.
  Derivation of the Boussineq condition is presented in Equation \ref{eq:bouss}.
  }
\begin{center}
\begin{tabular}{ll}
\hline
Variable & Value (at 300K, if applicable)\\
\hline
\hline
Rayleigh number & $[1.5\cdot 10^6, 4\cdot 10^7]$\\
\hline
Reynolds number & $\sim 416$\\
\hline
Gravity $g_z$ (fixed) & $9.8 \text{m}/\text{s}^2$\\
\hline
Specific heat $c_p$ & $4.187 \text{kj}/\text{kg}$\\
\hline
Thermal expansion $\beta$ & $0.303\cdot 10^{-3} \,1/\text{K}$\\
\hline
Density $\rhoref$ & $995.65 \text{kg}/\text{m}^3$\\
\hline
Kinematic Viscoscity $\nu$ & $0.801\cdot 10^{-6} \text{m}^2/\text{s}$\\
\hline
Dynamic Viscoscity $\mu$ & $0.798$\\
\hline
Laminar Prandtl number & $5.43$\\
\hline
Boussinesq condition & $7.6\cdot 10^{-6}$\\
\hline
\end{tabular}
\label{tb:flowvariables}
\end{center}
\end{table}

\begin{table}[h!]
\caption[EnKF covariance inflation test parameters]{EnKF covariance inflation test parameters.}
\begin{center}
\begin{tabular}{ccc}
\hline
{\bf Window Length ~~~~} & {\bf Additive Inflation~~} & {\bf Multiplicative Inflation~~}\\
\hline
\hline
~30s ~& 0.10 & 0.00 \\ \hline
~60s ~& 0.20 & 0.00 \\ \hline
~90s ~& 0.40 & 0.00 \\ \hline
~120s ~& 0.30 & 0.00 \\ \hline
~150s ~& 0.70 & 0.00 \\ \hline
~180s ~& 0.70 & 0.00 \\ \hline
~210s ~& 1.00 & 0.00 \\ \hline
~240s ~& 1.00 & 0.00 \\ \hline
~270s ~& 0.90 & 0.00 \\ \hline
~300s ~& 1.00 & 0.00 \\ \hline
~330s ~& 1.00 & 0.00 \\ \hline
~360s ~& 0.80 & 0.00 \\ \hline
~390s ~& 0.30 & 0.00 \\ \hline
~420s ~& 0.40 & 0.00 \\ \hline
~450s ~& 1.40 & 0.10 \\ \hline
~480s ~& 1.10 & 0.00 \\ \hline
~510s ~& 0.70 & 0.10 \\ \hline
~540s ~& 1.50 & 1.40 \\ \hline
~570s ~& 1.40 & 1.20 \\ \hline
~600s ~& 1.50 & 0.70 \\
\hline
\end{tabular}
\end{center}
\label{table:EnKFCovInfl}
\end{table}

\begin{table}[h!]
\caption[EKF covariance inflation test parameters]{EKF covariance inflation test parameters.}
\begin{center}
\begin{tabular}{ccc}
\hline
{\bf Window Length ~~~~} & {\bf Additive Inflation~~} & {\bf Multiplicative Inflation~~}\\
\hline
\hline
~30s ~& 0.10 & 0.00 \\ \hline
~60s ~& 0.00 & 1.50 \\ \hline
~90s ~& 0.00 & 1.50 \\ \hline
~120s ~& 0.00 & 1.50 \\ \hline
~150s ~& 0.00 & 1.50 \\ \hline
~180s ~& 0.20 & 1.50 \\ \hline
~210s ~& 0.60 & 1.50 \\ \hline
~240s ~& 0.80 & 1.50 \\ \hline
~270s ~& 0.90 & 1.50 \\ \hline
~300s ~& 1.30 & 0.60 \\ \hline
~330s ~& 0.50 & 0.70 \\ \hline
~360s ~& 0.50 & 0.90 \\ \hline
~390s ~& 0.90 & 0.40 \\ \hline
~420s ~& 0.60 & 0.30 \\ \hline
~450s ~& 0.80 & 0.10 \\ \hline
~480s ~& 0.60 & 0.40 \\ \hline
~510s ~& 0.20 & 0.30 \\ \hline
~540s ~& 0.00 & 0.40 \\ \hline
~570s ~& 0.00 & 0.30 \\ \hline
~600s ~& 1.20 & 0.20 \\
\hline
\end{tabular}
\end{center}
\label{table:EKFCovInfl}
\end{table}

\begin{table}[h!]
\caption[EnSRF covariance inflation test parameters]{EnSRF covariance inflation test parameters.}
\begin{center}
\begin{tabular}{ccc}
\hline
{\bf Window Length ~~~~} & {\bf Additive Inflation~~} & {\bf Multiplicative Inflation~~}\\
\hline
\hline
~30s ~& 0.00 & 0.00 \\ \hline
~60s ~& 1.50 & 0.80 \\ \hline
~90s ~& 1.40 & 0.40 \\ \hline
~120s ~& 1.50 & 0.90 \\ \hline
~150s ~& 0.10 & 1.30 \\ \hline
~180s ~& 0.00 & 0.80 \\ \hline
~210s ~& 0.10 & 1.30 \\ \hline
~240s ~& 0.00 & 1.50 \\ \hline
~270s ~& 0.00 & 1.40 \\ \hline
~300s ~& 0.80 & 1.50 \\ \hline
~330s ~& 1.20 & 0.50 \\ \hline
~360s ~& 1.40 & 1.00 \\ \hline
~390s ~& 1.40 & 1.30 \\ \hline
~420s ~& 1.50 & 1.30 \\ \hline
~450s ~& 1.50 & 1.40 \\ \hline
~480s ~& 1.50 & 0.30 \\ \hline
~510s ~& 0.00 & 1.40 \\ \hline
~540s ~& 0.00 & 1.40 \\ \hline
~570s ~& 0.00 & 1.00 \\ \hline
~600s ~& 0.00 & 0.80 \\
\hline
\end{tabular}
\end{center}
\label{table:EnSRFCovInfl}
\end{table}

\begin{table}[h!]
\caption[ETKF covariance inflation test parameters]{ETKF covariance inflation test parameters.}
\begin{center}
\begin{tabular}{ccc}
\hline
{\bf Window Length ~~~~} & {\bf Additive Inflation~~} & {\bf Multiplicative Inflation~~}\\
\hline
\hline
~30s ~& 0.10 & 0.00 \\ \hline
~60s ~& 0.20 & 0.00 \\ \hline
~90s ~& 0.30 & 0.00 \\ \hline
~120s ~& 0.30 & 0.00 \\ \hline
~150s ~& 0.20 & 0.00 \\ \hline
~180s ~& 0.20 & 0.00 \\ \hline
~210s ~& 0.50 & 0.00 \\ \hline
~240s ~& 0.50 & 0.00 \\ \hline
~270s ~& 0.50 & 0.00 \\ \hline
~300s ~& 0.20 & 0.00 \\ \hline
~330s ~& 0.40 & 0.10 \\ \hline
~360s ~& 0.20 & 0.10 \\ \hline
~390s ~& 0.20 & 0.10 \\ \hline
~420s ~& 0.60 & 0.10 \\ \hline
~450s ~& 1.50 & 0.20 \\ \hline
~480s ~& 1.20 & 0.40 \\ \hline
~510s ~& 0.30 & 0.10 \\ \hline
~540s ~& 0.90 & 0.50 \\ \hline
~570s ~& 1.50 & 0.40 \\ \hline
~600s ~& 1.50 & 0.30 \\
\hline
\end{tabular}
\end{center}
\label{table:ETKFCovInfl}
\end{table}

\begin{table}[h]
\caption[Summary of Lorenz '63 parameters]
  {
    Summary of Lorenz '63 parameters.
    With $\rho=28$ this system is strongly nonlinear, and prediction for long windows is difficult.
  }
\label{table:lorenz63params}
\begin{center}
\begin{tabular}{ll} \hline
Parameter & Value \\ \hline
\hline
$\sigma$ & 10 \\\hline
$\beta$ & 8/3 \\\hline
$\rho$ & 28 \\\hline
\end{tabular}
\end{center}
\end{table}

\chapter{Derivation of Ehrhard-M\"{u}ller Equations}

Following the derivation by Harris \shortcite{harris2011predicting}, itself a representation of the derivation of Gorman \shortcite{gorman1986} and namesakes Ehrhard and M\"{u}ller \shortcite{ehrhard1990dynamical}, we derive the equations governing a closed loop thermosyphon.

Similar to the derivation of the governing equations of computational fluid dynamics in Appendix C, we start with a small but finite volume inside the loop.
Here, however, the volume is described by $\pi r^2 R \text{d} \phi$ for $r$ the interior loop size (such that $\pi r^2$ is the area of a slice) and $R\text{d}\phi$ the arc length (width) of the slice.
Newton's second law states that momentum is conserved, such that the sum of the forces acting upon our finite volume is equal to the change in momentum of this volume.
Therefore we have the basic starting point for forces $\sum F$ and velocity $u$ as
\begin{equation} \sum F = \rho \pi r^2 R \text{d}\phi \diff{u}{t} .\end{equation}
The sum of the forces is $\sum F = F_{\{p,f,g\}}$ for net pressure, fluid shear, and gravity, respectively.
We write these as
\begin{align} & F_p = -\pi r^2 \text{d} \phi \pdiff{p}{\phi}\\
& F_w = -\rho \pi r^2 \text{d} \phi f_w\\
& F_g = -\rho \pi r^2 \text{d} \phi g \sin (\phi)\end{align}
where $\partial p /\partial \phi$ is the pressure gradient, $f_w$ is the wall friction force, and $g \sin (\phi)$ is the vertical component of gravity acting on the volume.

We now introduce the Boussinesq approximation which states that both variations in fluid density are linear in temperature $T$ and density variation is insignificant except when multiplied by gravity.
The consideration manifests as
\begin{equation*} \rho = \rho (T) \simeq \rho _\text{ref} (1 - \beta (T - T_\text{ref}) \end{equation*}
where $\rho _0$ is the reference density and $T_\text{ref}$ is the reference temperature, and $\beta$ is the thermal expansion coefficient.
The second consideration of the Boussinesq approximation allows us to replace $\rho$ with this $\rhoref$ in all terms except for $F_g$.
We now write momentum equation as
\begin{equation} -\pi r^2 \dphi \pdiff{p}{\phi} - \rhoref \phi r^2 R \dphi f_w
- \rhoref (1 - \rho (T- T_\text{ref}) ) \pi r^2 R \dphi g \sin (\phi) = \rhoref \pi r^2 R \dphi \diff{u}{t}. \end{equation}
Canceling the common $\pi r^2$, dividing by $R$, and pulling out $\dphi$ on the LHS we have
\begin{equation} -\dphi \left ( \pdiff{p}{\phi}  \frac{1}{R} - \rhoref f_w - \rhoref (1 - \rho (T- T_\text{ref}) ) g \sin (\phi) \right ) = \rhoref \dphi \diff{u}{t}. \label{eq:EM07} \end{equation}
We integrate this equation over $\phi$ to eliminate many of the terms, specifically we have
\begin{align*}
& \int _{0} ^{2\pi} -\dphi \pdiff{p}{\phi} \frac{1}{R} \rightarrow 0\\
& \int _{0} ^{2\pi} -\dphi \rhoref g \sin (\phi) \rightarrow 0\\
& \int _{0} ^{2\pi} -\dphi \rhoref \beta T_\text{ref} g \sin (\phi) \rightarrow 0.\end{align*}
Since $u$ (and hence $\diff{u}{\phi}$) and $f_w$ do not depend on $\phi$, we can pull these outside an integral over $\phi$ and therefore the momentum equation is now
\begin{equation*} 2\pi f_w \rho _0 + \int _{0} ^{2\pi} \dphi \rhoref \beta T g \sin (\phi) = 2\pi \diff{u}{\phi} \rhoref .\end{equation*}
Diving out $2\pi$ and pull constants out of the integral we have our final form of the momentum equation
\begin{equation} f_w \rhoref + \frac{\rhoref \beta g }{2 \pi} \int _{0} ^{2\pi} \dphi T \sin (\phi) = \diff{u}{\phi} \rhoref \label{eq:EM10}.\end{equation}
Now considering the conservation of energy within the thermosyphon, the energy change within a finite volume must be balanced by transfer within the thermosyphon and to the walls.
The internal energy change is given by
\begin{equation} \rhoref \pi r^2 R \dphi \left ( \pdiff{T}{t} + \frac{u}{R}\pdiff{T}{\phi} \right ) \label{eq:EMeg1}\end{equation}
which must equal the energy transfer through the wall, which is, for $T_w$ the wall temperature:
\begin{equation} \dot{q} = -\pi r^2 R \dphi h_w (T - T_w) . \label{eq:EMeg2} \end{equation}
Combining Equations \ref{eq:EMeg1} and \ref{eq:EMeg2} (and canceling terms) we have the energy equation:
\begin{equation} \left ( \pdiff{T}{t} + \frac{u}{R}\pdiff{T}{\phi} \right ) = \frac{-h_w}{\rhoref c_p} \left( T - T_w \right ) \label{eq:EMeq}.\end{equation}
The $f_w$ which we have yet to define and $h_w$ are fluid-wall coefficients and can be described by \shortcite{ehrhard1990dynamical}:
\begin{align*} & h_w = h_{w_0} \left ( 1 + K h(|x_1|) \right ) \\
& f_w = \frac{1}{2} \rhoref f_{w_0} u .\end{align*}
We have introduced an additional function $h$ to describe the behavior of the dimensionless velocity $x_1 \alpha u$.
This function is defined piece-wise as
\begin{equation*} h (x) = \left \{ \begin{array}{ll} x^{1/3} & ~~\text{when} ~x \geq 1\\ p (x) & ~~\text{when} ~ x <1 \end{array} \right. \end{equation*}
where $p(x)$ can be defined as $p(x) = \left( 44x^2 -55 x^3 + 20x^4 \right ) /9$ such that $p$ is analytic at 0 \shortcite{harris2011predicting}.

Taking the lowest modes of a Fourier expansion for $T$ for an approximate solution, we consider:
\begin{equation} T(\phi , t) = C_0 (t) + S(t) \sin (\phi ) + C(t) \cos (\phi) . \end{equation}
By substituting this form into Equations \ref{eq:EM10} and \ref{eq:EMeq} and integrating, we obtain a system of three equations for our solution.
We then follow the particular nondimensionalization choice of Harris et al such that we obtain the following ODE system, which we refer to as the Ehrhard-M\"{u}ller equations:
\begin{align}
& \diff{x_1}{t'} = \alpha (x_2 - x_1),\\
& \diff{x_2}{t'} = \beta x_1 - x_2 (1 + Kh(|x_1|)) - x_1x_3,\\
& \diff{x_3}{t'} = x_1x_2 - x_3 (1 + Kh(|x_1|)) .\end{align}
The nondimensionalization is given by the change of variables
\begin{align}
& t' = \frac{h_{w_0}}{\rhoref c_p}t,\\
& x_1 = \frac{\rhoref c_p }{R h_{w_0}} u, \\
& x_2 = \frac{1}{2} \frac{\rhoref c_p \beta g}{ R h_{w_0} f_{w_0}} \Delta T_{3-9}, \\
& x_3 = \frac{1}{2} \frac{\rhoref c_p \beta g}{ R h_{w_0} f_{w_0}} \left ( \frac{4}{\pi} \Delta T_w - \Delta T_{6-12} \right )
\end{align}
and
\begin{align}
& \alpha = \frac{1}{2} R c_p f_{w_0} / h_{w_0} ,\\
& \gamma = \frac{2}{\pi} \frac{\rhoref c_p \beta g}{Rh_{w_0} f_{w_0}} \Delta T_w. \end{align}

Through careful consideration of these non-dimensional variable transformations we verify that $x_1$ is representative of the mean fluid velocity, $x_2$ of the temperature difference between the 3 and 9 o'clock positions on the thermosyphon, and $x_3$ the deviation from the vertical temperature profile in a conduction state \shortcite{harris2011predicting}.

\chapter{Derivation of CFD Solving Code}

Closely following the original derivation of Stokes \shortcite{stokes}, I present a full derivation of the main equations governing the flow of water at temperatures near 300K (an incompressible, heat-conducting, Newtonian fluid), and detail their implementation in a finite volume numerical solver.
Detailed equations follow from my notes based on the reference work of {\em Computational Fluid Dynamics} written by Anderson et al for the von Karman Institute lectures \shortcite{anderson1995computational}.
From the available approaches, I consider a fixed finite volume to derive the equations used in OpenFOAM.

\section{Continuity Equation}

Considering a finite volume element with side lengths $\Delta x, \Delta y$ and $\Delta z$.
The amount of mass that enters any given face on this volume is a function of the fluid density $\rho$, the velocity tangent to this face, and the area of the face.
For the a side with edges specified by both $\Delta x$ and $\Delta y$, let the velocity tangent to this face be $u$ and this mass flow in this side is equal to
$$ \rho u \Delta x \Delta y .$$
Assuming $u$ is positive in this direction, the mass flux out of the opposite side needs to reflect possible changes in velocity $u$ and density $\rho$ through the volume and can be written
$$ -(\rho + \Delta \rho) (u+\Delta u) \Delta x \Delta y.$$
Similarly for the other two directions tangent to our volume, we assign velocities $v$ and $w$ with incoming mass
\begin{align*} &\rho v \Delta x \Delta z\\
&\rho w \Delta y \Delta z\end{align*}
and outgoing mass
\begin{align*} &-(\rho + \Delta \rho )( v+ \Delta v) \Delta x \Delta z\\
&-(\rho+\Delta \rho) (w+\Delta w) \Delta y \Delta z.\end{align*}
The rate of mass accumulation in the volume,
$$ \Delta \rho (\Delta x \Delta y \Delta z) / \Delta t, $$
must be equal to the sum of the mass entering and leaving the volume (given by the six equations above).
Equating these, with a little cancellation and dividing by $(\Delta x \Delta y \Delta z)$ we are left with
\begin{equation} \frac{\Delta \rho}{\Delta t} + \frac{\Delta (\rho u)}{\Delta x} + \frac{\Delta (\rho v)}{\Delta y} + \frac{\Delta (\rho w)}{\Delta z} = 0. \end{equation}
As $\Delta t \to 0$, this can be written in terms of the partial derivatives
\begin{equation} \frac{\partial \rho}{\partial t} + \frac{\partial (\rho u)}{\partial x} + \frac{\partial (\rho v)}{\partial y} + \frac{\partial (\rho w)}{\partial z} = 0. \label{eq:NScont} \end{equation}

For an incompressible fluid, we have that $\partial \rho / \partial \{ t,x,y,z\} = 0$ such the continuity equation becomes
\begin{equation} \frac{\partial u}{\partial x} + \frac{\partial v}{\partial y} + \frac{\partial w}{\partial z} = 0. \label{eq:NScontIco} \end{equation}

\section{Momentum Equation}

We now consider the second conservation law: the conservation of momentum.
This states that the rate of change of momentum must equal the net momentum flux into the control volume in addition to any external forces (e.g. gravity) on the control volume.
Again consider a small but finite volume element with side lengths $\Delta x, \Delta y$ and $\Delta z$.
Since the conservation of momentum applies in $x,y$ and $z$ directions, without loss of generality I consider only the $x$ direction.

The rate of change of momentum with respect to time in our volume element is
\begin{equation*} \frac{\partial }{\partial t} \left ( \rho u \right) \Delta x \Delta y \Delta z .\end{equation*}

The momentum flux into the volume in the $x$ direction is the mass flux times $x$-velocity, since momentum is mass$\times$velocity.
Recall the mass flux is given by $\rho u \times A$ for $A$ the area orthogonal to $x$, in this case $\Delta y \Delta z$.
Therefore the momentum flux is
\begin{equation*} u \times \rho u (\Delta y \Delta z) .\end{equation*}
Similarly, with the additional consideration of the change within the volume, the momentum flux leaving through the opposite side is
\begin{equation*} - \left ( u\rho u + \frac{\partial}{\partial x} \left( u\rho u\right ) \Delta x \right ) \Delta y \Delta z .\end{equation*}

The $y$ and $z$ direction momemtum flux (of $x$ direction momentum) is found in the same way and we have the fluxes in those directions as
\begin{align*} u \times \rho v (\Delta x \Delta z) & ~~~\& ~~~ - \left ( u\rho v + \frac{\partial}{\partial y} \left( u\rho v\right ) \Delta y \right ) \Delta x \Delta z \\
u \times \rho w (\Delta x \Delta y) & ~~~\& ~~~ - \left ( u\rho w + \frac{\partial}{\partial z} \left( u\rho w\right ) \Delta z \right ) \Delta x \Delta y . \end{align*}

Equating the change in internal momentum with the sum of the fluxes, and the sum of the external forces in the $x$ direction $\mathbf{F}_1$, we have the most basic form of the momentum equation:

\begin{align*} \frac{\partial }{\partial t} \left ( \rho u \right) \Delta x \Delta y \Delta z &= \frac{\partial }{\partial t} \left ( \rho u \right) \Delta x \Delta y \Delta z  - \left ( u\rho u + \frac{\partial}{\partial x} \left( u\rho u\right ) \Delta x \right ) \Delta y \Delta z +  u \rho v (\Delta x \Delta z) \\
& - \left ( u\rho v + \frac{\partial}{\partial y} \left( u\rho v\right ) \Delta y \right ) \Delta x \Delta z + u \rho w (\Delta x \Delta y)\\
& - \left ( u\rho w + \frac{\partial}{\partial z} \left( u\rho w\right ) \Delta z \right ) \Delta x \Delta y + \sum \mathbf{F}_1.\end{align*}

Canceling the first terms from each flux we are left
\begin{align*} \frac{\partial }{\partial t} \left ( \rho u \right) \Delta x \Delta y \Delta z &=  -  \frac{\partial}{\partial x} \left( u\rho u\right ) \Delta x  \Delta y \Delta z - \frac{\partial}{\partial y} \left( u\rho v\right ) \Delta y \Delta x \Delta z\\
& - \frac{\partial}{\partial z} \left( u\rho w\right ) \Delta z \Delta x \Delta y + \sum \mathbf{F}_1.\end{align*}

Pulling out the $\Delta x \Delta y \Delta z$ and moving just the forces to the right we have

\begin{align*} \Delta x \Delta y \Delta z \left (\frac{\partial }{\partial t} \left ( \rho u \right)  +  \frac{\partial}{\partial x} \left( u\rho u\right )  + \frac{\partial}{\partial y} \left( u\rho v\right ) +\frac{\partial}{\partial z} \left( u\rho w\right ) \right ) &= \sum \mathbf{F}_1.\end{align*}

Applying the product rule to the partial derivatives with respect to $x,y$ and $z$ we find the continuity equation which we know to be zero, and are left
\begin{align*} \Delta x \Delta y \Delta z \left ( \rho \frac{\partial u }{\partial t} + \rho u \frac{\partial u }{\partial x} + \rho v \frac{\partial u}{\partial y} + \rho w \frac{\partial u}{\partial z} \right ) &= \sum \mathbf{F}_1 .\end{align*}

The sum of the forces in the $x$ direction includes the force of gravity that acts on the mass of the entire volume ($g_1 \times \rho \Delta x\Delta y \Delta z$) and the surface stress.
We assume that gravity acts only in the $z$ direction where the term $g_3 \times \rho \Delta x \Delta y \Delta z$ shows up in the sum of the forces and is otherwise zero.
The $x$-direction force of the stresses acting on the volume is the product of the stress and the area on which it acts, and in Cartesian coordinates this is either direct or shear stress.
The $x$-normal stress, which we denote $s_{xx}$ has forces on both sides tangent to $x$ given by
\begin{equation*} s_{xx} \Delta y \Delta z ~~~\& ~~~ \left ( s_{xx} + \frac{\partial s_{xx} }{\partial x} \Delta x \right ) \Delta y \Delta z \end{equation*}
for which the sum is
\begin{equation*} \frac{\partial s_{xx} }{\partial x} \Delta x \Delta y \Delta z .\end{equation*}
The forces of the shear stress on the faces planar to $x$ are similarly
\begin{equation*} \frac{\partial s_{yx} }{\partial x} \Delta x \Delta y \Delta z ~~~\&~~~\frac{\partial s_{zx} }{\partial x} \Delta x \Delta y \Delta z .\end{equation*}

Since the pressure $p$ does not apply shear stress, $p$ is only a part of $s_{xx}$ and not $s_{yx},s_{zx}$.
In each direction we denote the viscous stress as $\tau _{ij}$ for $i$ and $j$ the directions $x,y,$ and $z$.
Thus for general fluids the sum of the forces in the $x$ direction is
\begin{equation*} \left ( \frac{\partial p }{\partial x} + \frac{\partial \tau _{xx}}{\partial x} + \frac{\partial \tau _{yx}}{\partial y} + \frac{\partial \tau _{zx}}{\partial z} \right ) \Delta x \Delta y \Delta z . \end{equation*}

Since the thermosyphon is filled with water, not blood, we can happily assume that the fluid is Newtonian, and the rapture has been avoided.
This means that we can relate the viscous stress to the local rate of deformation linearly, by the viscocity $\mu$:
\begin{align*} \tau _{xx} &= 2 \mu \frac{\partial u}{\partial x} \\
\tau _{yx} &= \mu \left ( \frac{\partial v}{\partial x} + \frac{\partial u}{\partial y} \right )\\
\tau _{zx} &= \mu \left ( \frac{\partial w}{\partial x} + \frac{\partial u}{\partial z} \right )\end{align*}

Putting this together for our Newtonian fluid we have by the equality of mixed partials and the continuity equation we have the sum of all the forces given by
\begin{align*} & \left (-\frac{\partial p}{\partial x} + 2 \mu \frac{\partial u}{\partial x^2} + \mu \frac{\partial}{\partial y} \frac{\partial v}{\partial x} + \mu \frac{\partial }{\partial y} \frac{\partial u}{\partial y} + \mu \frac{\partial }{\partial z} \frac{\partial u}{\partial z} + \mu \frac{\partial }{\partial z} \frac{\partial w}{\partial x}\right ) \Delta x \Delta y \Delta z\\
& = \left ( -\frac{\partial p}{\partial x} + \mu \frac{\partial u}{\partial x^2} + \mu\frac{\partial u}{\partial y^2} + \mu \frac{\partial u}{\partial z^2}  + \mu \frac{\partial u}{\partial x^2} + \mu \frac{\partial}{\partial x} \frac{\partial v}{\partial y} + \mu \frac{\partial }{\partial x} \frac{\partial w}{\partial z} \right ) \Delta x \Delta y \Delta z\\
& = \left (-\frac{\partial p}{\partial x} + \mu \frac{\partial u}{\partial x^2} + \mu\frac{\partial u}{\partial y^2} + \mu \frac{\partial u}{\partial z^2}  + \mu \frac{\partial }{\partial x}\left ( \frac{\partial u}{\partial x} + \frac{\partial v}{\partial y} + \frac{\partial w}{\partial z} \right ) \right ) \Delta x \Delta y \Delta z\\
& = \left ( -\frac{\partial p}{\partial x} + \mu \frac{\partial u}{\partial x^2} + \mu\frac{\partial u}{\partial y^2} + \mu \frac{\partial u}{\partial z^2}\right ) \Delta x \Delta y \Delta z\end{align*}

Thus the momemtum equation for $x$ is
\begin{align}  \rho \left ( \frac{\partial u }{\partial t} + u \frac{\partial u }{\partial x} + v \frac{\partial u}{\partial y} + w \frac{\partial u}{\partial z} \right ) &= -\frac{\partial p}{\partial x} + \mu \frac{\partial u}{\partial x^2} + \mu\frac{\partial u}{\partial y^2} + \mu \frac{\partial u}{\partial z^2}\end{align}

In tensor notation this is often written for all three equations (recall, the above is only for $x$):

\begin{equation} \rho \frac{\partial u_i}{\partial t} + \frac{\partial}{\partial x_j} \left( u_j u_i \right) = -\frac{\partial p}{\partial x_i}  + \mu \frac{\partial u_i}{\partial x_j^2}  + \rho g_i . \end{equation}

Working our way towards OpenFOAM's implementation of this equation, write the filtered equation for averaged quantities of pressure $\bar{p}$, density $\bar{\rho}$ and velocity $\bar{u}$, and assume that the density is constant except for when multiplied by gravitational effects, writing $\bar{\rho} = \overline{\rho} / \rho _0$:
\begin{equation} \rho _0 \left ( \frac{\partial \bar{u}_i}{\partial t} + \frac{\partial}{\partial x_j} \left( \bar{u}_j \bar{u}_i \right) \right )
= -\frac{\partial \bar{p}} {\partial{x_i}} +  \mu \frac{\partial \bar{u}_i}{\partial x_j^2} + \bar{\rho} g_i. \end{equation}

Since we will be relying on the use of a turbulence model, we reintroduce the stress tensor $\tau _{ij}$ and split $\tau_{ij}$ into the mean stress tensor $\tau _{ij}$ and the turbulent stress tensor $\tau ^* _{ij}$.
In tensor notation the mean stress tensor can be written
\begin{equation*} \tau _{ij} = \mu \left ( \left ( \frac{\partial \bar{u} _i}{\partial x_j} + \frac{\partial \bar{u} _j }{\partial x_i} \right) - \frac{2}{3} \frac{\partial \bar{u}_k} {\partial x_k} \delta _{ij} \right ) \end{equation*}
for $\delta _{ij} = 0$ if and only if $i=j$.
Assume constant viscocity $\mu = \mu _0$, divide by $\rho _0$ and incorporate the latter for the stress tensor and we now have the momentum equation as
\begin{equation*} \frac{\partial \bar{u}_i}{\partial t} + \frac{\partial}{\partial x_j} \left( \bar{u}_j \bar{u}_i \right)
= -\frac{\partial } {\partial{x_i}} \frac{\bar{p}}{\rho_0} + \frac{\partial }{\partial x_j} \left ( \nu _0 \left ( \left ( \frac{\partial \bar{u}_i} {\partial x_j} + \frac{\partial \bar{u} _j}{\partial x_i} \right ) - \frac{2}{3}\left ( \frac{\partial \bar{u} _k}{\partial x_k }\right ) \delta _{ij}\right ) - \tau ^* _{ij} \right ) + \bar{\rho} g_i \end{equation*}
where $\nu _0 = \mu _0 /\rho _0$.
We can then decompose the turbulent stress tensor $\tau ^* _{ij}$ further using $R_{ij}$ to denote the Reynolds stress tensor
\begin{equation*} \tau ^* _{ij} = R_{ij} = R^D _{ij} + \frac{2}{3} k \delta _{ij} \end{equation*}
for $R_{ij}^D$ the deviatoric part and $k = R_{ii}/2$ the turbulent kinetic energy.
For LES this becomes the sub-grid kinetic energy.
Note that OpenFOAM outputs the resolved kinematic pressure
\begin{equation*} \tilde{p} = \frac{\bar{p} }{\rho _0} + \frac{2}{3}k \end{equation*}
and uses the Boussinesq approximation for the last term
\begin{equation*} g_i\frac{\bar{\rho}}{\rho _0} = g_i \left ( 1 + \frac{\bar{\rho} - \rho _0}{\rho _0} \right ) = g_i \rho _k  = g_i \left (1 - \beta \left (\overline{T} - T _0\right ) \right ) .\end{equation*}

\section{Energy Equation}

Again following the derivation of Anderson, we derive the energy equation into the form used in OpenFOAM \shortcite{anderson2009governing}.
Since our problem deals with an incompressible fluid, we are concerned mainly with the temperature.

The first law of thermodynamics applied to a finite volume says that the rate of change of energy inside the fluid element is equal to the net flux of heat into the element plus the rate of work done by external forces on the volume.
Writing the absolute internal thermal energy as $e$ and the kinematic energy per unit mass as $v^2/2$ we have the rate of change of total internal energy as
\begin{equation} \frac{\partial}{\partial t} \left ( \rho \left ( e + \frac{v^2}{2} \right ) \right ) \label{eq:energy1} \end{equation}
and the net transfer of energy through the control volume is
\begin{equation} \frac{\partial}{\partial x} \left ( u \rho \left ( e + \frac{v^2}{2} \right ) \right ) + \frac{\partial}{\partial y} \left ( v \rho \left ( e + \frac{v^2}{2} \right ) \right ) + \frac{\partial}{\partial z} \left ( w \rho \left ( e + \frac{v^2}{2} \right ) \right ) \label{eq:energy2} .\end{equation}

Heat flux $\dot{q}$ is defined to be positive for flux leaving the control volume, so we write the net heat flux entering the volume as
\begin{equation} - \frac{\partial q_x}{\partial x} - \frac{\partial q_y}{\partial y} - \frac{\partial q_z}{\partial z} \label{eq:energy3}.\end{equation}

Since heat flux by thermal conduction is determined by Fourier's Law of Heat Conduction as
\begin{equation} q_{\{x,y,z\}} = -k \frac{\partial T}{\partial \{x,y,z\}} \label{eq:energy6} \end{equation}
we write Equation \ref{eq:energy3} as
\begin{equation} k\left ( \frac{\partial T}{\partial x^2} - \frac{\partial T}{\partial y^2} - \frac{\partial T}{\partial z^2} \label{eq:energy5} \right ).\end{equation}

As before, the shear stress in the $i$ direction on face $j$ is denoted $s_{ij}$ such that the net rate of work done by these stresses is the sum of their components
\begin{equation} - \frac{\partial }{\partial x} \left ( u s_{xx} + v s_{xy} + w s_{xz} \right) - \frac{\partial }{\partial y} \left (u s_{yx} + v s_{yy} + w s_{yz} \right) - \frac{\partial }{\partial z}\left (u s_{zx} + v s_{zy} + w s_{zz} \right) \label{eq:energy6}.\end{equation}

Putting this all together in tensor notation, the temperature (energy) equation is thus
\begin{equation} \frac{\partial }{\partial t} \left ( \rho e\right ) + \frac{\partial}{\partial x_j} \left ( \rho e u_j \right )
=
-k \frac{\partial T}{\partial x_k^2}.
\end{equation}

Averaging this equation and separating the heat flux into the average and turbulent parts $q = \overline{q} + q^*$ we have
\begin{equation} \frac{\partial }{\partial t} \left ( \rho \overline{e}\right ) + \frac{\partial}{\partial x_j} \left ( \rho \overline{e} \overline{u}_j \right )
=
- \frac{\partial q_k^*}{\partial x_k}
- \frac{\partial \overline{q}_k}{\partial x_k}
\end{equation}

\section{Implementation}

The PISO (Pressure-Implicit with Splitting of Operators) algorithm derives from the work of \shortcite{issa1986solution}, and is complementary to the SIMPLE (Semi-Implicit Method for Pressure-Linked Equations) \shortcite{patankar1972calculation} iterative method.
The main difference of the PISO and SIMPLE algorithms is that in the PISO, no under-relaxation is applied and the momentum corrector step is performed more than once \shortcite{ferziger1996computational}.
They sum up the algorithm in nine steps:
\begin{itemize}
\item Set the boundary conditions
\item Solve the discretized momentum equation to compute an intermediate velocity field
\item Compute the mass fluxes at the cell faces
\item Solve the pressure equation
\item Correct the mass fluxes at the cell faces
\item Correct the velocity with respect to the new pressure field
\item Update the boundary conditions
\item Repeat from step \#3 for the prescribed number of times
\item Repeat (with increased time step).
\end{itemize}

The solver itself has 647 dependencies, of which I present only a fraction.
The main code is straight forward, relying on include statements to load the libraries and equations to be solved.
\lstset{language=C++,
	basicstyle=\ttfamily\scriptsize\singlespacing,
	keywordstyle=\color{blue},
	stringstyle=\color{red},
	commentstyle=\color{green},
	morecomment=[l][\color{magenta}]{\#},
	frame=L,
	xleftmargin=\parindent,
				numbersep=5pt,
	breaklines=true,		breakatwhitespace=false,	    escapeinside={\%*}{*)}
}

\begin{lstlisting}[language=C++,firstline=48,lastline=53]
#include "fvCFD.H"
#include "singlePhaseTransportModel.H"
#include "RASModel.H" // AJR edited 2013-10-14
#include "radiationModel.H"
#include "fvIOoptionList.H"
#include "pimpleControl.H"
\end{lstlisting}

The main function is then

\begin{lstlisting}[language=C++,firstline=57,lastline=70]
int main(int argc, char *argv[])
{
    #include "setRootCase.H"
    #include "createTime.H"
    #include "createMesh.H"
    #include "readGravitationalAcceleration.H"
    #include "createFields.H"
    #include "createIncompressibleRadiationModel.H"
    #include "createFvOptions.H"
    #include "initContinuityErrs.H"
    #include "readTimeControls.H"
    #include "CourantNo.H"
    #include "setInitialDeltaT.H"
    pimpleControl pimple(mesh);
\end{lstlisting}

We then enter the main loop.
This is computed for each time step, prescribed before the solver is applied.
Note that the capacity is available for adaptive time steps, choosing to keep the Courant number below some threshold, but I do not use this.
For the distributed ensemble of model runs, it is important that each model complete in nearly the same time, so that the analysis is not waiting on one model and therefore under-utilizing the available resources.

\begin{lstlisting}[language=C++,firstline=71,lastline=88]
    while (runTime.loop())
    {
        #include "readTimeControls.H"
        #include "CourantNo.H"
        #include "setDeltaT.H"
        while (pimple.loop())
        {
            #include "UEqn.H"
            #include "TEqn.H"
            while (pimple.correct())
            {
                #include "pEqn.H"
            }
        }
	if (pimple.turbCorr())
	{
	    turbulence->correct();
	}
\end{lstlisting}

Opening up the equation for $U$ we see that Equation

\begin{lstlisting}[language=C++]
    // Solve the momentum equation
    fvVectorMatrix UEqn
    (
        fvm::ddt(U)
      + fvm::div(phi, U)
      + turbulence->divDevReff(U)
     ==
        fvOptions(U)
    );
    UEqn.relax();
    fvOptions.constrain(UEqn);
    if (pimple.momentumPredictor())
    {
        solve
        (
            UEqn
         ==
            fvc::reconstruct
            (
                (
                  - ghf*fvc::snGrad(rhok)
                  - fvc::snGrad(p_rgh)
                )*mesh.magSf()
            )
        );
        fvOptions.correct(U);
    }
\end{lstlisting}

Solving for $T$ is

\begin{lstlisting}[language=C++]
{
    alphat = turbulence->nut()/Prt;
    alphat.correctBoundaryConditions();
    volScalarField alphaEff("alphaEff", turbulence->nu()/Pr + alphat);
    fvScalarMatrix TEqn
    (
        fvm::ddt(T)
      + fvm::div(phi, T)
      - fvm::laplacian(alphaEff, T)
     ==
        radiation->ST(rhoCpRef, T)
      + fvOptions(T)
    );
    TEqn.relax();
    fvOptions.constrain(TEqn);
    TEqn.solve();
    radiation->correct();
    fvOptions.correct(T);
    rhok = 1.0 - beta*(T - TRef); // Boussinesq approximation
}
\end{lstlisting}

Finally, we solve for the pressure $p$ in ``pEqn.H'':

\begin{lstlisting}[language=C++]
{
    volScalarField rAU("rAU", 1.0/UEqn.A());
    surfaceScalarField rAUf("Dp", fvc::interpolate(rAU));
    volVectorField HbyA("HbyA", U);
    HbyA = rAU*UEqn.H();
    surfaceScalarField phig(-rAUf*ghf*fvc::snGrad(rhok)*mesh.magSf());
    surfaceScalarField phiHbyA
    (
        "phiHbyA",
        (fvc::interpolate(HbyA) & mesh.Sf())
      + fvc::ddtPhiCorr(rAU, U, phi)
      + phig
    );
    while (pimple.correctNonOrthogonal())
    {
        fvScalarMatrix p_rghEqn
        (
            fvm::laplacian(rAUf, p_rgh) == fvc::div(phiHbyA)
        );
        p_rghEqn.setReference(pRefCell, getRefCellValue(p_rgh, pRefCell));
        p_rghEqn.solve(mesh.solver(p_rgh.select(pimple.finalInnerIter())));
        if (pimple.finalNonOrthogonalIter())
        {
            // Calculate the conservative fluxes
            phi = phiHbyA - p_rghEqn.flux();
            // Explicitly relax pressure for momentum corrector
            p_rgh.relax();
            // Correct the momentum source with the pressure gradient flux
            // calculated from the relaxed pressure
            U = HbyA + rAU*fvc::reconstruct((phig - p_rghEqn.flux())/rAUf);
            U.correctBoundaryConditions();
        }
    }
    #include "continuityErrs.H"
    p = p_rgh + rhok*gh;
    if (p_rgh.needReference())
    {
        p += dimensionedScalar
        (
            "p",
            p.dimensions(),
            pRefValue - getRefCellValue(p, pRefCell)
        );
        p_rgh = p - rhok*gh;
    }
}
\end{lstlisting}

The final operation being the conversion of pressure to hydrostatic pressure,
\begin{equation*} p _\text{rgh} = p - \rho _k g_h . \end{equation*}
This ``pEqn.H'' is then re-run until convergence is achieved, and the PISO loop begins again.

\chapter{foamLab Framework}
\label{app:foamlab}

The core components of the code, used to make all results herein, are as follows.
The main link is the MATLAB function modelDAinterface.m, which uses models as a standard class and the DA algorithms as function files.
For the use of the OpenFOAM, the two main files are the MATLAB class OpenFOAM.m and the complementary shell script foamLab.sh with which it directly interacts.
First I present the full modelDAinterface.m function, then one of the assimilation algorithms, and the class definition for OpenFOAM.
Only the non-localized modelDAinterface.m is included, as the local version is much longer.

Performing an experiment with a given model is then performed by inputing an observation timeseries into modelDAinterface.m, as well as the model to use for prediction, the filter to use to assimilate, and other details about the observations.

{\scriptsize 
\begin{lstlisting}[language=Matlab]
function [x_f_all,x_a_all] = modelDAinterface(modelname,filter,y,y_t,H,R,window,delta,mu,varargin)
% modelDAinterface.m
%
% written by Andy Reagan
% 2013-10-14
%
% combine abstract model class and DA code for experimentation
% things should be general enough for this to work both lorenz and foam
%
% INPUT:
%   modelname       the model: try either 'lorenz63' or 'OpenFOAM'
%   filter          filtering algorithm: 'EKF' or 'EnKF' for starters
%   y               matrix of observations [y_0,y_1,...,y_N]
%   y_t             vector of observation time [0,1,2,...,100]
%   H               observation operator
%   R               obs error covariance (block-diagonal,sym)
%   window          length of assimilation window
%   delta           multiplicative covariance influation
%   mu              additive covariance inflation
% (optional, see more below)
%   writecontrol,t  specify density of forecast output:
%    		    default is only at the windows
%                       
%   numens,N        specify number of ensemble members
%                   e.g. 'numens',20
%                   default is 10
%
%   params,paramCell
%                   set the model parameters for choosing
%
% OUTPUT:
%   x_f_all       matrix store of the forecasts
%   x_a_all       matrix store of the analyses

%% defaults
save_int = window; % output times
N = 10; % number of ens members
modelnamestr = func2str(modelname);
silent = 1;
setparams = 0;
np = 1;
seed = sum(100*clock);

%% parse input
for i=1:2:nargin-9
    switch varargin{i}
        case 'writecontrol'
            save_int = varargin{i+1};
        case 'numens'
            N =   varargin{i+1};
        case 'params'
            setparams = 1; params = varargin{i+1};
        case 'silent'
            silent = varargin{i+1};
        case 'OFmodel'
            OFCase = varargin{i+1};
        case 'randseed'
            seed = varargin{i+1};
    end
end

%% reset the random seed
%% rng('shuffle','twister'); %% for 2013a
RandStream.setDefaultStream(RandStream.create('mt19937ar','seed',seed)); %% for 2010b
randStr = sprintf('%06d',randi(1000000,1));

%%%%%%%%%%%%%%%%%%%%%%%%%%%%%%%%%%%%%%%%%%%%%%%%%%%%%%%%%%%%%
%% INITIALIZE
%%%%%%%%%%%%%%%%%%%%%%%%%%%%%%%%%%%%%%%%%%%%%%%%%%%%%%%%%%%%%

%% set up for having a EKF
switch filter
    case 'EKF'
        model = modelname();
        if setparams
            model.params = params;
        end
        model.init(sprintf('%s/model',randStr));
        model.window = window;
        x_a = model.x;
        cov_init_max = 10; % start the covariance large, it will calm down
        %% passing this by value instead of via disk
        %% save('P.mat',cov_init_max*eye(model.dim));
        p_a = cov_init_max*eye(model.dim);
        %% make the dimension local
        dim = model.dim;
        
    case {'direct','none'}
        %% set up for having no DA, either direct input or none
        model = modelname();
        if setparams
            model.params = params;
        end
        model.init(sprintf('%s/model',randStr));
        model.window = window;
        dim = model.dim*np;
        x_a = model.x;
        
    case {'EnKF','ETKF','EnSRF'}
        %% set up ensemble if we're doing an ensemble filter
        ensemble = cell(N,1);
        ensemble{1} = modelname();
        if setparams
            ensemble{1}.params = params;
        end
        ensemble{1}.init(sprintf('%s/ens%02d',randStr,1));
        %% start the first model
        ensemble{1}.window = window;
        dim = ensemble{1}.dim*np;
        
        X_f = ones(dim,N);
        X_a = ones(dim,N);
        X_a(:,1) = ensemble{1}.x;
        %% make the rest of the ensemble
        for i=2:N
            ensemble{i} = modelname();
            if setparams
                ensemble{i}.params = params;
            end
            ensemble{i}.init(sprintf('%s/ens%02d',ranStr,i));
            ensemble{i}.window = window;
            X_a(:,i) = ensemble{i}.x;
        end %% ensemble initialization
end %% switch

num_windows = length(y(1,:));
x_f_all = ones(dim,num_windows);
x_a_all = ones(dim,num_windows);


%% predict the given timeseries
i=0;
for time=0:window:max(y_t)
    i=i+1; %% loop counter
    
    %%%%%%%%%%%%%%%%%%%%%%%%%%%%%%%%%%%%%%%%%%%%%%%%
    %% FORECAST: RUN MODEL
    %%%%%%%%%%%%%%%%%%%%%%%%%%%%%%%%%%%%%%%%%%%%%%%%
    
    switch filter
        case 'EKF'
            %% set the IC
            model.x = x_a;
            %% run both models
            model.runTLM(p_a);
            model.run();
            %% store the forecast
            x_f_all(:,i) = model.x;
        case {'direct','none'}
            %% run the forward model for no DA
            %% set the IC
            model.x = x_a;
            %% run model
            model.run();
            %% store the forecast
            x_f_all(:,i) = model.x;
        case {'EnKF','EnSRF','ETKF'}
            %% run the ens if we're doing the EnKF
            %% run the whole ensemble forward
            for j=1:N
                % set the IC
                ensemble{j}.x = X_a(:,j);
                % run
                ensemble{j}.run();
                % record the model output from run
                X_f(:,j) = ensemble{j}.x;
            end
            %% store the forecast
            x_f_all(:,i) = mean(X_f,2);
    end %% switch
    
    %%%%%%%%%%%%%%%%%%%%%%%%%%%%%%%%%%%%%%%%%%%%%%%%%
    %% ASSIMILATE
    %%%%%%%%%%%%%%%%%%%%%%%%%%%%%%%%%%%%%%%%%%%%%%%%%
    
    switch filter
        case 'EKF'
            %% pre-analysis multiplicative inflation
            model.p_f = (1+delta).*model.p_f;
            [x_a,p_a] = EKF(model.x,y(:,i),H,R,model.p_f,'silent');
            %% post analysis additive inflation
            p_a = p_a + mu*diag(rand(size(x_a)));
            x_a_all(:,i) = x_a;
            
        case 'direct'
            x_a = direct(x_f_all(:,i),y(:,i),H);
            x_a_all(:,i) = x_a;
            
        case 'none'
            x_a = none(x_f_all(:,i));
            x_a_all(:,i) = x_a;
            
        case 'EnKF'
            X_a = EnKF(X_f,y(:,i),H,R,delta,'silent');
            %% post analysis additive inflation
            X_a = X_a+mu*rand(size(X_a));
            x_a_all(:,i) = mean(X_a,2);
            
        case 'EnSRF'
            X_a = EnSRF(X_f,y(:,i),H,R,delta,'silent');
            %% post analysis additive inflation
            X_a = X_a+mu*rand(size(X_a));
            x_a_all(:,i) = mean(X_a,2);
                     
        case 'ETKF'
            X_a = ETKF(X_f,y(:,i),H,R,delta,'silent');
            %% post analysis additive inflation
            X_a = X_a+mu*rand(size(X_a));
            x_a_all(:,i) = mean(X_a,2);

    end %% switch
    %% wash, rinse, repeat

end %% for loop
end %% function




\end{lstlisting}}

Now we take a look at how simply a data assimilation algoritm can be defined within this framework.
The code for the Ensemble Transform Kalman Filter is contained in ETKF.m and is:

{\scriptsize 
\begin{lstlisting}[language=Matlab]
function X_a = ETKF(X_f,y_o,H,R,delta,varargin)

% ETKF.m
%
% Ensemble Transform Kalman Filter: compute the new analysis
% 
% INPUTS
%  X_f:     the whole forecast, columns are model state
%
% OUTPUTS
%  X_a:     the analysis matrix, columns are model state analysis
%
% written by Andy Reagan
% 2013-10-18

silent = 0;

i=1; 
while i<=length(varargin), 
  argok = 1; 
  if ischar(varargin{i}), 
    switch varargin{i},
        case 'silent',       silent  = 1;
        otherwise, argok=0; 
    end
  end
  if ~argok, 
    fprintf('invalid argument %s\n',varargin{1})
  end
  i = i+1; 
end

N = length(X_f(1,:));

%% let x_f now be the average
x_f = mean(X_f,2);

X_f_diff = X_f - repmat(x_f,1,N);

%% multiplicative inflation
X_f_diff = sqrt(1+delta).*X_f_diff;

p_f = 1/(N-1)*(X_f_diff*X_f_diff');

K = (p_f*H')/(R+H*p_f*H');

d = y_o - H*x_f;

x_a = x_f + K*d;

%% perform the transform
p_a = inv((N-1)*eye(N)+(H*X_f_diff)'/(R)*(H*X_f_diff));

T = sqrtm((N-1)*p_a);

X_a_diff = X_f_diff*T;

%% now compute the analysis for each ensemble member
X_a=repmat(x_a,1,N) + X_a_diff;

if ~silent
  fprintf('K is\n')
  disp(K);
  fprintf('p_f is\n')
  disp(p_f);
end


\end{lstlisting}}

Finally we include the OpenFOAM class definition.
This definition relies entirely on a shell script foamLab.sh, which for brevity's sake we do not include.

{\scriptsize 
\begin{lstlisting}[language=Matlab]
classdef OpenFOAM<handle
    %OpenFOAM: class for an OpenFOAM run
    %   Goal is to implement a class for which MATLAB can
    %   control and run an OpenFOAM instance
    %
    %   USAGE:
    %       -initialize an OpenFOAM directory, read to be run, at 'testCase'
    %        ens1 = OpenFOAM('2D-foamLabBase-mesh40000');
    %        ens1.init('testCase');
    %       -run OpenFOAM on this case
    %        ens1.run();
    
    properties
        %% basics
        x
        dim
        tstep = 1;
        time = 0;
        windowLen = 20;
        params
        
        foamLab = '/users/a/r/areagan/work/2013/data-assimilation/OpenFOAM/foamLab.sh';
        TIME=0;
        ETIME=20;
        DIR='testCase'
        FLUX=10000; 
        BCbottom = '340';
        BCtop = '290';
        BC = 'fixedValue';
        TURB = 'off';
        BASE = '2D-foamLabSmall';
        TURBMODEL = 'kEpsilon';
        T
        U
        p
        p_rgh
        Tstep = 1;
        WRITEp = 6;
        VALUE = [];
        np = 6;
    end %% properties
    
    methods
        function self = OpenFOAM(varargin)
            %% intialize the class
            if nargin > 0
                self.BASE = varargin{1};
            end
            switch self.BASE
                case '2D-foamLabSmall'
                    self.dim = 600;
                case '2D-foamLabBase'
                    self.dim = 36000;
                case '2D-foamLabBase-mesh40000'
                    self.dim = 40000;
                otherwise
                    warning('chosen case does not have a predifined dimension. edit class to include it');
            end
            self.T = 300*ones(self.dim,1);
            self.U = zeros(self.dim,3); %%1e-6*ones(40832,3);
            self.p = zeros(self.dim,1); %%1e-6*ones(40832,1);
            self.p_rgh = zeros(self.dim,1); %%1e-6*ones(40832,1);
            self.x = ones(self.dim*self.np,1);
            self.vectorize('in');
        end %% constructor
        function init(self,varargin)
            if nargin > 1
                self.DIR = varargin{1};
            end
            %% intialize a a clean dir
            command = sprintf('%s -i -d %s -B %s -h %s -g %s -b %s -q %s -Q %s',self.foamLab,self.DIR,self.BASE,self.BCtop,self.BCbottom,self.BC,self.TURB,self.TURBMODEL);
            fprintf('%s\n',command);
            system(command);
        end %% init
        function read(self,time,varargin)
            %% read in a specific time variable
            %% first: have foamLab write the csv
            caseDir = sprintf('/users/a/r/areagan/OpenFOAM/areagan-2.2.1/run/%s',self.DIR);
            system(sprintf('%s -r %g -d %s -D %d',self.foamLab,time,caseDir,self.dim));
            self.T = csvread(sprintf('%s/%g/T.csv',caseDir,time));
            self.U = csvread(sprintf('%s/%g/U.csv',caseDir,time));
            self.p = csvread(sprintf('%s/%g/p.csv',caseDir,time));
            self.p_rgh = csvread(sprintf('%s/%g/p_rgh.csv',caseDir,time));
            self.vectorize('in');
        end %% read
        function write(self,time,varargin)
            %% write out variable
            self.vectorize('out');
            caseDir = sprintf('/users/a/r/areagan/OpenFOAM/areagan-2.2.1/run/%s',self.DIR);
            csvwrite(sprintf('%s/%g/T.csv',caseDir,time),self.T);
            csvwrite(sprintf('%s/%g/U.csv',caseDir,time),self.U);
            csvwrite(sprintf('%s/%g/p.csv',caseDir,time),self.p);
            csvwrite(sprintf('%s/%g/p_rgh.csv',caseDir,time),self.p_rgh);
            system(sprintf('%s -W %g -d %s',self.foamLab,time,caseDir));
        end %% write
         function run(self,varargin)
            caseDir = sprintf('/users/a/r/areagan/OpenFOAM/areagan-2.2.1/run/%s',self.DIR);
            
            %% write out the current x before running
            self.write(self.time)
            
            tmpcommand = sprintf('%s -x -d %s -t %g -e %g -l %g -w %g -c %g -B %s -D %d',self.foamLab,caseDir,self.time,self.time+self.windowLen,self.tstep,self.WRITEp,self.windowLen,self.BASE,self.dim);
            fprintf('%s\n',tmpcommand);
            system(tmpcommand);

            %% update time
            self.time = self.time+self.windowLen;
            self.read(self.time);
        end %% run
        function vectorize(self,direction)
            %% self.np is the number of variables used, order a permutation on their storage
	    switch np
                case {1,6}
		     order = 0:5;
		case 3 
		     order = [0 0 1 2 0 0];
	    end
            switch direction
                case 'in' %% into the datavec
                    self.x((1:self.np:(self.dim*self.np+1-self.np))+order(1)) = self.T;
                    self.x((1:self.np:(self.dim*self.np+1-self.np))+order(2)) = self.U(:,1);
                    self.x((1:self.np:(self.dim*self.np+1-self.np))+order(3)) = self.U(:,2);
                    self.x((1:self.np:(self.dim*self.np+1-self.np))+order(4)) = self.U(:,3);
                    self.x((1:self.np:(self.dim*self.np+1-self.np))+order(5)) = self.p;
                    self.x((1:self.np:(self.dim*self.np+1-self.np))+order(6)) = self.p_rgh;
                case 'out' %% out of x into T,U,p,p_rgh
                    self.T = self.x((1:self.np:(self.dim*self.np+1-self.np))+order(1));
                    self.U(:,1) = self.x((1:self.np:self.dim*self.np)+order(2));
                    self.U(:,2) = self.x((1:self.np:(self.dim*self.np+1-self.np))+order(3));
                    self.U(:,3) = self.x((1:self.np:(self.dim*self.np+1-self.np))+order(4));
                    self.p = self.x((1:self.np:self.dim*self.np)+order(5));
                    self.p_rgh = self.x((1:self.np:self.dim*self.np)+order(6));
            end %% switch
        end %% vectorize
        function destruct(self)
            %% destroy the folder, for sanity's sake
            system(sprintf('\\rm -rf /users/a/r/areagan/OpenFOAM/areagan-2.2.1/run/%s',self.DIR));
        end
    end %% methods
end %% classdef






\end{lstlisting}}

For more complete code documentation, and to download the source code, visit the GitHub repository: \url{https://github.com/andyreagan}.

\end{document}